\RequirePackage{etoolbox}
\csdef{input@path}{%
 {./}%
 {figures/}%
}%
\csgdef{bibdir}{./}%

\documentclass[reqno,letterpaper]{amsart}

\RequirePackage[OT1]{fontenc}

\usepackage{import, subfiles}
\usepackage{url}
\usepackage{amsmath, amsthm, bm, bbm, amssymb, amsfonts}
\usepackage{mathtools}

\newenvironment{proofof}[1]{\renewcommand{\proofname}{Proof of #1}\proof}{\endproof}
\usepackage{algpseudocode, algorithm}
\usepackage{rotating}
\usepackage{graphicx}
\usepackage{multirow}
\usepackage{color}
\usepackage{wrapfig}
\usepackage{graphicx} %
\usepackage{lipsum}
\usepackage{leftidx}
\usepackage{natbib}
\definecolor{DarkBlue}{rgb}{0, 0, .4}
\usepackage[linktocpage=true,colorlinks,linkcolor=blue,citecolor=blue,urlcolor=blue,pdfborder={0 0 0}]{hyperref}
\usepackage[capitalize,sort]{cleveref}
\usepackage{autonum}

\usepackage[arxiv]{optional}

\makeatletter
\renewcommand\paragraph{\@startsection{paragraph}{4}{\z@}%
  {\medskipamount}%
  {-1em}%
  {\bfseries}%
}
\makeatother

\usepackage{caption}
\usepackage{subcaption}

\crefname{nlem}{Lemma}{Lemmas}
\crefname{nprop}{Proposition}{Propositions}
\crefname{ncor}{Corollary}{Corollaries}
\crefname{nthm}{Theorem}{Theorems}
\crefname{exa}{Example}{Examples}
\crefname{assumption}{Assumption}{Assumptions}

\usepackage{jhh-misc2}

\newcommand{\BDrep}{\ensuremath{\operatorname{B-Rep}}}
\newcommand{\isBDrep}{\gets \BDrep}
\newcommand{\ILrep}{\ensuremath{\operatorname{IL-Rep}}}
\newcommand{\isILrep}{\gets \ILrep}
\newcommand{\Trep}{\ensuremath{\operatorname{T-Rep}}}
\newcommand{\isTrep}{\gets \Trep}
\newcommand{\Rrep}{\ensuremath{\operatorname{R-Rep}}}
\newcommand{\isRrep}{\gets \Rrep}

\newcommand{\DBFrep}{\ensuremath{\operatorname{DB-Rep}}}
\newcommand{\isDBFrep}{\gets \DBFrep}
\newcommand{\PLFrep}{\ensuremath{\operatorname{PL-Rep}}}
\newcommand{\isPLFrep}{\gets \PLFrep}
\newcommand{\Vrep}{\ensuremath{\operatorname{SB-Rep}}}
\newcommand{\isVrep}{\gets \Vrep}

\newcommand{\nuinv}{\nu^{\gets}}

\newcommand{\muinv}{\mu^{\gets}}

\allowdisplaybreaks[3]

\numberwithin{equation}{section}

\begin{document}
\sloppy

\title{Truncated Random Measures}

\thanks{$^{\star}$First authorship is shared jointly by T.~Campbell and J.~H.~Huggins.}

\author[T.~Campbell]{Trevor Campbell$^{\star}$}
\address{Computer Science and Artificial Intelligence Laboratory (CSAIL) \\ Massachusetts Institute of Technology}
\urladdr{http://www.trevorcampbell.me/}
\email{tdjc@mit.edu}

\author[J.~H.~Huggins]{Jonathan H.~Huggins$^{\star}$}
\address{Computer Science and Artificial Intelligence Laboratory (CSAIL) \\ Massachusetts Institute of Technology}
\urladdr{http://www.jhhuggins.org/}
\email{jhuggins@mit.edu}

\author[J.~P.~How]{Jonathan How}
\address{Laboratory for Information and Decision Systems (LIDS) \\ Massachusetts Institute of Technology}
\urladdr{http://www.mit.edu/~jhow/}
\email{jhow@mit.edu}

\author[T.~Broderick]{Tamara Broderick}
\address{Computer Science and Artificial Intelligence Laboratory (CSAIL) \\ Massachusetts Institute of Technology}
\urladdr{http://www.tamarabroderick.com}
\email{tbroderick@csail.mit.edu}

\date{\today}

\begin{abstract}
Completely random measures (CRMs) and their normalizations are a rich source of 
Bayesian nonparametric priors.
Examples include the beta, gamma, and Dirichlet processes. 
In this paper we detail two major classes of sequential CRM representations---\emph{series representations}
and \emph{superposition representations}---within which we organize both novel 
and existing sequential representations that can
be used for simulation and posterior inference.
These two classes and their constituent representations subsume 
existing ones that have previously been developed
in an ad hoc manner for specific processes. 
Since a complete infinite-dimensional CRM cannot be used explicitly for computation,
sequential representations are often truncated for tractability.
We provide truncation error analyses for each type of sequential
representation, as well as their normalized versions, thereby generalizing 
and improving upon existing truncation error bounds in the literature.
We analyze the computational complexity of the sequential representations,
which in conjunction with our error bounds allows us to directly compare representations
and discuss their relative efficiency. 
We include numerous applications of our theoretical results to commonly-used 
(normalized) CRMs, demonstrating that our results 
enable a straightforward representation and analysis of
CRMs that has not previously been available in a Bayesian nonparametric context.

\end{abstract}

\maketitle

\begin{small}
~\\~\\
\renewcommand\contentsname{\!\!\!\!}
\setcounter{tocdepth}{2}
\vspace{-15mm}
\tableofcontents
\end{small}

\section{Introduction}

In many data sets, we can view the data points as exhibiting a collection of underlying traits.
For instance, each document in the \emph{New York Times} might touch on a number of topics or themes,
an individual's genetic data might be a product of the populations to which their ancestors belonged, 
or a user's activity on a social network might be dictated by their varied personal interests. 
When the traits are not directly observed, a common approach is to model each trait as having
some frequency or rate in the broader population \citep{airoldi:2014:handbook}. The inferential goal is to learn
these rates as well as whether---and to what extent---each data point exhibits each trait.
Since the traits are unknown a priori, their cardinality is also typically unknown.

As a data set grows larger, we can reasonably expect the number of traits %
to increase as well. 
In the cases above, for example, we expect to uncover more topics as we read more documents,
more ancestral populations as we examine more individuals' genetic data, and more unique interests
as we observe more individuals on a social network.
\emph{Bayesian nonparametric} (BNP) priors provide a flexible, principled approach to creating 
models in which the number of exhibited traits is random, can grow without bound,
and may be learned as part of the inferential procedure. 
By generating a countable infinity of potential 
traits---where any individual data point exhibits only finitely many---these models
enable growth in the number of observed traits with the size of the data set. 

In practice, however, it is impossible to store a countable infinity of random variables in memory
or learn the distribution over a countable infinity of variables in finite time. 
Conjugate priors and likelihoods have been developed \citep{orbanz:2010:conjugate}
that theoretically circumvent the infinite representation altogether
and perform exact Bayesian posterior inference \citep{Broderick:2014}. However,
these priors and likelihoods are often just a single piece within a more complex generative model,
and ultimately
an approximate posterior inference scheme such as Markov Chain Monte Carlo (MCMC) or 
variational Bayes (VB) is required. These approximation schemes often necessitate a full and explicit representation of the 
latent variables.

One option is to approximate the infinite-dimensional prior with a related finite-dimensional prior:
that is, to replace the infinite collection of random traits by a finite subset of ``likely'' traits.
To do so, first enumerate the countable infinity of traits in the full model and write 
$(\psi_k, \theta_k)$ for each paired trait $\psi_k$ (e.g.\ a topic in a document) and its rate or frequency $\theta_k$. 
Then the discrete measure $\Theta \defined \sum_{k=1}^{\infty} \theta_k \delta_{\psi_k}$
captures the traits/rates in a sequence indexed by $k$.
The $(\psi_k, \theta_k)$ pairs are random in the Bayesian model, so $\Theta$ 
is a random measure. In many cases, the distribution of $\Theta$ can be defined
by specifying a sequence of simple, familiar distributions for the finite-dimensional $\psi_k$ and $\theta_k$, known as a \emph{sequential representation}.
Given a sequential representation of $\Theta$, a natural way to choose a subset of traits is to keep the 
first $K<\infty$ traits and discard the rest, resulting in an approximate measure $\Theta_K$. 
This approach is called \emph{truncation}. Note that it is also possible to truncate 
by removing atoms with weights less than a specified threshold \citep{Argiento:2015,Muliere:1998},
though this approach is not as easily incorporated in posterior inference algorithms.

Sequential representations have been shown to exist for \emph{completely random measures} 
(CRMs)~\citep{Kingman:1967,Ferguson:1972}, a large class of nonparametric priors
that includes such popular models as the
beta process \citep{Hjort:1990,Kim:1999} and the gamma process \citep{Ferguson:1972,Kingman:1975,Brix:1999,Titsias:2008,James:2013}. 
Numerous sequential representations of CRMs have been developed in the literature \citep{Ferguson:1972,Bondesson:1982,Rosinski:1990,Rosinski:2001,James:2014,Broderick:2014}.
CRM priors are often paired with likelihood processes---such as 
the Bernoulli process \citep{Thibaux:2007}, negative binomial process \citep{Zhou:2012,Broderick:2015}, and Poisson likelihood
process \citep{Titsias:2008}. 
The likelihood process determines how much each trait is expressed by each data point. 
Sequential representations also exist
for \emph{normalized completely random measures} (NCRMs) (sometimes referred to as normalized random measures with independent increments) \citep{Perman:1992,Perman:1993,James:2002,Pitman:2003,Regazzini:2003,Lijoi:2010,James:2009ux}, 
which provide random distributions over traits, such as
the Dirichlet process~\citep{Ferguson:1973,Sethuraman:1994}.
NCRMs are typically paired with a likelihood that assigns each data point
to a single trait using the NCRM as a discrete distribution.

Since (N)CRMs have many possible sequential representations, a method is required for determining which 
 to use for the application at hand and, once a representation is selected, for choosing a truncation level. %
Our main contributions enable the principled selection of both representation and truncation 
level using approximation error:
\benum
\item We provide a comprehensive characterization of the
different types of sequential representations for (N)CRMs, %
filling in many gaps in the literature of sequential representations along the way. 
We classify these representations into two major groups: \emph{series representations}, which are constructed by transforming
a homogeneous Poisson point process; and \emph{superposition representations}, which are 
the superposition of infinitely many Poisson point processes with finite rate measures. 
We also introduce two novel sequential representations for (N)CRMs. 
\item We provide theoretical guarantees on the approximation error induced when truncating these sequential representations. 
We give the error as a function of the prior process, the likelihood process, and the level of truncation.
While truncation error bounds for (N)CRMs have been studied previously, 
past work has focused on specific combinations of (N)CRM priors and likelihoods---in particular,
the Dirichlet-multinomial~\citep{Sethuraman:1994,Ishwaran:2001,Ishwaran:2002,Blei:2006a},
beta-Bernoulli~\citep{Paisley:2012,DoshiVelez:2009}, generalized beta-Bernoulli~\citep{Roy:2014},
and gamma-Poisson~\citep{Roychowdhury:2015} processes.
In the current work, we give much more general results for bounding the truncation error.
\eenum

Our results fill in large gaps in the analysis of truncation error, which is often measured in terms of the 
$L^{1}$ (a.k.a.\ total variation) distance between the data distributions induced by the full and truncated priors. 
We provide the first analysis of truncation error
for some sequential representations of the beta process with Bernoulli likelihood \citep{Thibaux:2007},
for the beta process with negative binomial likelihood \citep{Zhou:2012,Broderick:2015}, 
and for the normalization of the generalized gamma process \citep{Brix:1999}, 
the $\sigma$-stable process, and the generalized inverse gamma~\citep{Lijoi:2005,Lijoi:2010} with discrete likelihood.
Moreover, even when truncation results already exist in the literature
\citep{Ishwaran:2001,DoshiVelez:2009,Paisley:2012,Roychowdhury:2015}, we improve
on those error bounds by a factor of two. The reduction arises from
our use of the point process machinery of CRMs, circumventing the total variation bound used 
originally by \citet{Ishwaran:2001,Ishwaran:2002b} upon which most modern truncation analyses are built.
We obtain our truncation error guarantees by bounding the probability that data drawn from the full model
will use a feature that is not available to the truncated model. 
Thinking in terms of this probability provides a more intuitive interpretation of our bounds
that can be communicated to practitioners and used to guide them in their choice of truncation level.

The remainder of this paper is organized as follows. 
In \cref{sec:bkg}, we provide background material on CRMs and establish notation.
In our first main theoretical section, \cref{sec:seqcrms},
we describe seven different sequential CRM representations, including
four series representations %
and three superposition representations, %
two of which are novel.
Next, we provide a general theoretical analysis of the truncation error 
for series and superposition representations in \cref{sec:trunc}. 
We provide analogous theory for the normalized versions of each representation
in \cref{sec:ncrm-trunc} via an infinite extension of the ``Gumbel-max
trick''~\citep{Gumbel:1954,Maddison:2014ww}.
We determine the complexity of simulating each representation in
\cref{sec:sampling}. In \cref{sec:simulation}, we summarize our results
(\cref{tbl:results}) and provide 
advice on how to select sequential representations in practice.
Proofs for all results developed in this paper are provided in the the supplemental 
article \citet{supplement}.

\section{Background}
\label{sec:bkg}
\subsection{CRMs and truncation}

Consider a Poisson point process on $\mathbb{R}_+ \defined \left[0, \infty\right)$ with rate measure $\nu(\dee\theta)$ such that 
\[
  \nu(\mathbb{R}_+)=\infty && \text{and}&&
  \int \min(1, \theta) \nu(\dee\theta) < \infty. \label{eq:nuassump}
\]
Such a process generates a countable infinity of values
$\left(\theta_k\right)_{k=1}^\infty$, $\theta_k \in \mathbb{R}_+$, having an
almost surely finite sum $\sum_{k=1}^{\infty} \theta_k < \infty$.
In a BNP trait model, we interpret each $\theta_k$ as the rate or frequency of the $k$-th trait.
Typically, each $\theta_k$ is paired with a parameter $\psi_k$ associated with the $k$-th trait 
(e.g., a topic in a document or a shared interest on a social network). 
We assume throughout that $\psi_k \in \Psi$ for some space $\Psi$
and $\psi_k \distiid G$ for some
distribution $G$. Constructing
a measure by placing mass $\theta_k$ at atom location $\psi_k$ results in
a \emph{completely random measure} (CRM) \citep{Kingman:1967}. As shorthand, we will write $\distCRM(\nu)$
for the completely random measure generated as just described:
\begin{align}
\Theta &\defined \sum_k \theta_k \delta_{\psi_k} \dist \distCRM(\nu).
\end{align}
The trait distribution $G$ is left implicit in the notation as it has no effect on our results. Further, the possible fixed-location and deterministic components of a 
CRM~\citep{Kingman:1967} are not considered here for brevity;
these components can be added (assuming they are purely atomic)
and the analysis modified without undue effort.
The CRM prior on $\Theta$ is typically combined
with a likelihood that generates trait counts for each data point. 
Let $h(\cdot \given \theta)$ be
a proper probability mass function on $\nats \cup \{0\}$ for all 
$\theta$ in the support of $\nu$ (though the present work may be easily extended to likelihoods with support in $\reals$). 
Then  
a collection of conditionally independent observations
$X_{1:N} \defined \left\{X_n\right\}_{n=1}^N$  given $\Theta$
 are distributed according to the 
\emph{likelihood process} $\distLP(h, \Theta)$, i.e.
\begin{align}
X_n &\defined \sum_k x_{nk}\delta_{\psi_k} \distiid \distLP(h, \Theta),
\end{align}
if 
$x_{nk} \sim h(x \given \theta_k)$ independently across $k$
and \iid\ across $n$.
The desideratum that each $X_n$ expresses a finite 
number of traits is encoded by the assumption that
\begin{align}
\int (1-h(0 \given \theta))\nu(\dee\theta) < \infty.\label{eq:hassump}
\end{align}
Since the trait counts are typically latent in a full generative model
specification, define the observed data $Y_{n} \given X_{n} \distind f(\cdot \given X_{n})$ for
a conditional density $f$ with respect to a measure $\mu$ on some space. 
For instance, if the sequence $\left(\theta_k\right)_{k=1}^\infty$ represents the topic rates
in a document corpus, $X_{n}$ might capture how many words in document
$n$ are generated from each topic and $Y_{n}$ might be the 
observed collection of words for that document.

Since the sequence $(\theta_k)_{k=1}^{\infty}$ is countably infinite, it may be
difficult to simulate or perform posterior
inference in this model. One approximation scheme is to define the \emph{truncation}
$\Theta_{K} \defined \sum_{k=1}^K \theta_k \delta_{\psi_k}$. Since it is finite, the truncation $\Theta_{K}$ can
be used for exact simulation or in posterior inference---but some error arises 
from not using the full CRM $\Theta$. 
To quantify this error, consider its propagation through the above Bayesian model. 
Define $Z_{1:N}$ and $W_{1:N}$ for $\Theta_{K}$ analogous to the definitions of $X_{1:N}$ and $Y_{1:N}$ for $\Theta$:
\[
Z_{n} \given \Theta_{K} &\distiid \distLP(h, \Theta_{K}), & 
W_{n} \given Z_{n} &\distind f(\cdot \given Z_{n}), & n=1,\dots,N.
\]
A standard approach to measuring the distance between $\Theta$ and $\Theta_{K}$
is to use the $L^{1}$ metric between the marginal densities 
$p_{N,\infty}$ and $p_{N,K}$ (with respect to some measure $\mu$)
of the final observations $Y_{1:N}$ and $W_{1:N}$~\citep{Ishwaran:2001,DoshiVelez:2009,Paisley:2012}:
\begin{align}
	& \frac{1}{2}\|p_{N,\infty} - p_{N,K}\|_{1} \defined \frac{1}{2}\int_{y_{1:N}} | p_{N,\infty}(y_{1:N}) - p_{N,K}(y_{1:N}) | \; \mu(\dee y_{1:N}).
\end{align}
All of our bounds on $ \frac{1}{2}\|p_{N,\infty} - p_{N,K}\|_{1}$ are also bounds on the probability that $X_{1:N}$ 
contains a feature that is not in the truncation $\Theta_K$ (cf. \cref{sec:trunc,sec:ncrm-trunc}).
This interpretation may be easier to digest since it does not depend on the observation model $f$ and is instead
framed in terms of the underlying traits the practitioner is trying to estimate. 

\subsection{The gamma-Poisson process}
\label{sec:commonratemeasures}
To illustrate the practical application of the theoretical developments
in this work, we provide a number of examples throughout involving the \emph{gamma process}~\citep{Brix:1999}, denoted $\distGammaP(\gamma, \lambda, d)$,
with discount parameter $d \in [0,1)$, scale parameter $\lambda > 0$, mass
parameter $\gamma > 0$, and rate measure
\begin{align}
\nu(\dee \theta) &= \gamma\frac{\lambda^{1-d}}{\Gamma(1-d)}
\theta^{-d-1}e^{-\lambda \theta}\dee \theta. 
\end{align}
Setting $d = 0$ yields the undiscounted gamma
process~\citep{Ferguson:1972,Kingman:1975,Titsias:2008}. 
The gamma process is often paired with a Poisson likelihood,
\begin{align}
  h(x\given\theta) &= \frac{\theta^x}{x!}e^{-\theta}.
\end{align}
Throughout the present work, we use the rate parametrization of the gamma distribution (to match
the gamma process parametrization), for which the density is given by
\[
\distGam(x; a, b) = \frac{b^{a}}{\Gamma(a)} x^{a-1}e^{-bx}. 
\]
\cref{app:examples} provides additional example applications of our theoretical results for two other CRMs: 
the beta process $\distBP(\gamma, \alpha, d)$~\citep{Teh:2009,Broderick:2012} with Bernoulli or negative binomial likelihood,
and the beta prime process~$\distBPP(\gamma,\alpha,d)$~\citep{Broderick:2014} with odds-Bernoulli likelihood.

\section{Sequential representations}\label{sec:seqcrms}

Sequential representations are at the heart of the study of truncated CRMs.
They provide an iterative method that can be terminated at any point to
yield a finite approximation to the infinite process, where the choice of
termination point determines the accuracy of the approximation. 
Thus, the natural first step in providing a coherent 
treatment of truncation analysis is to do the same for sequential representations.
In past work, two major classes of sequential representation have been used:
\emph{series representations} of the form $\sum_{k=1}^\infty \theta_k\delta_{\psi_k}$,
and \emph{superposition representations} of the form $\sum_{k=1}^\infty \sum_{i=1}^{C_k}\theta_{ki}\delta_{\psi_{ki}}$,
where each inner sum of $C_k$ atoms is itself a CRM.
This section examines four series representations \citep{Ferguson:1972,Bondesson:1982,Rosinski:1990,Rosinski:2001}
and three superposition representations (two of which are novel) \citep{Broderick:2012,Broderick:2014,James:2014}.
We show how previously-developed sequential representations for specific CRMs
fit into these seven general representations.
Finally, we discuss a stochastic mapping procedure that is useful in obtaining new representations 
from the transformation of others. 
Proofs for the results in this section may be found in \cref{app:seqproofs}.

\subsection{Series representations} 
\label{sec:seriesreps}

Series representations arise from the transformation of a homogeneous Poisson point process \citep{Rosinski:2001}.
They tend to be somewhat difficult to analyze due to the dependence 
between the atoms but also tend to produce very simple representations
with small truncation error (cf.~\cref{sec:trunc,sec:simulation}).  
Throughout the paper we let $\Gamma_{k} = \sum_{\ell=1}^{k}E_{\ell}$, 
$E_{\ell} \distiid \distExp(1)$, be the ordered jumps of a 
unit-rate homogeneous Poisson process on $\reals_+$, let $\nu$ 
be a measure on $\reals_+$ satisfying the basic conditions in \cref{eq:nuassump},
and let $\psi_k \distiid G$.

\paragraph{Inverse-L\'evy \citep{Ferguson:1972}}
Define $\nuinv(u) \defined \inf\left\{x : \nu\left([x, \infty)\right) \leq u \right\}$,
the inverse of the tail measure $\nu([x,\infty))$. 
We say $\Theta$ has an \emph{inverse-L\'evy} representation
and write $\Theta \isILrep(\nu)$ if
\[
\Theta &= \sum_{k=1}^{\infty}\theta_k\delta_{\psi_k},  &
& \text{with} & 
\theta_k &= \nuinv(\Gamma_{k}). \label{eq:invlevy}
\]
\citet{Ferguson:1972} showed that $\Theta\isILrep(\nu)$ implies $\Theta\dist\distCRM(\nu)$.
The inverse-L\'evy representation is analogous to the inverse CDF method for generating an arbitrary random
variable from a uniform random variable, with the homogenous Poisson process playing the role of the
uniform random variable.
It is also the optimal sequential representation in the sense that the sequence $\theta_{k}$ that it
generates is non-increasing. 
While an elegant and general approach, simulating the inverse-L\'evy representation 
is difficult, as inverting the function $\nu\left([x, \infty)\right)$ is 
analytically intractable except in a few cases.

\begin{exa}[Gamma process, $\distGammaP(\gamma, \lambda, 0)$]
We have $\nu([x, \infty)) = \gamma \lambda E_{1}(\lambda x)$, where 
$E_{1}(x) \defined \int_{x}^{\infty}u^{-1}e^{-u}\,\dee u$ is the exponential integral
function~\citep{Abramowitz:1964}. 
The inverse-L\'evy representation for $\distGammaP(\gamma, \lambda, 0)$ is thus
\[
\Theta = \sum_{k=1}^{\infty}\lambda^{-1}E_{1}^{-1}(\gamma^{-1}\lambda^{-1}\Gamma_{k})\delta_{\psi_{k}}.
\]
Neither $E_{1}$ nor its inverse can be computed in closed form, so one must
resort to numerical approximations. 
\end{exa}

\paragraph{Bondesson \citep{Bondesson:1982}}
We say $\Theta$ has a 
\emph{Bondesson representation} and write $\Theta \isBDrep(c, g)$ if
for $c > 0$ and $g$ a density on $\reals_{+}$,
\[
\Theta &= \sum_{k=1}^{\infty}\theta_{k}\delta_{\psi_{k}}, & 
&\text{with} &
\theta_{k} &= V_{k} e^{-\Gamma_{k}/c}, &
V_{k} &\distiid g. \label{eq:bondesson}
\]
\cref{prop:bondesson-representation} shows that Bondesson representations can be constructed for a 
large, albeit restricted, class of CRM rate measures. 
We offer a novel proof of \cref{prop:bondesson-representation} in
\cref{app:seqproofs} using the induction strategy introduced by \citet{Banjevic:2002}.
Similar proof ideas are also used to prove truncation error bounds for sequential representations 
in \cref{sec:trunc}.
We use a slight abuse of notation for brevity: if $\nu(\dee\theta)$ is a measure on $\reals_+$
that is absolutely continuous with respect to Lebesgue measure, 
then $\nu(\theta)$ is the density of $\nu(\dee \theta)$ with respect to the Lebesgue measure.

\begin{nthm}[Bondesson representation \citep{Bondesson:1982}] \label{prop:bondesson-representation}
Let $\nu(\dee \theta) = \nu(\theta)\dee \theta$ be a rate measure satisfying \cref{eq:nuassump}. 
If $\theta \nu(\theta)$ is nonincreasing, $\lim_{\theta \to \infty}\theta\nu(\theta) = 0$,
and $c_{\nu} \defined \lim_{\theta \to 0} \theta \nu(\theta) < \infty$, then 
$g_{\nu}(v) \defined -c_{\nu}^{-1}\der{}{v}[v \nu(v)]$ is a density on $\reals_{+}$ and
\[
\Theta \isBDrep(c_{\nu}, g_{\nu}) \quad \text{implies} \quad \Theta \dist \distCRM(\nu). \label{eq:ber}
\]
\end{nthm}

\begin{exa}[Bondesson representation for $\distGammaP(\gamma, \lambda, 0)$] \label{ex:b-rep-gamma-process}
The following representation for the gamma process with $d=0$
was described by \citet{Bondesson:1982} and \citet{Banjevic:2002}.
Since $\theta \nu(\theta) = \gamma \lambda e^{- \lambda \theta}$ is non-increasing
and $c_{\nu} = \lim_{\theta \to 0} \theta \nu(\theta) = \gamma \lambda$, 
we obtain $g_{\nu}(v) = \lambda e^{- \lambda v} = \distExp(v; \lambda)$. 
Thus, it follows from \cref{prop:bondesson-representation} that 
if $\Theta \isBDrep(\gamma\lambda, \distExp(\lambda))$, then $\Theta \dist 
\distGammaP(\gamma, \lambda, 0)$.
The condition that $\theta \nu(\theta)$ is non-increasing fails to hold if $d > 0$, so 
we cannot apply \cref{prop:bondesson-representation} to $\distGammaP(\gamma, \lambda, d)$ when $d > 0$. 
\end{exa}

\paragraph{Thinning \citep{Rosinski:1990}}
Using the nomenclature of \citet{Rosinski:2001}, we say $\Theta$ has a \emph{thinning representation} and write $\Theta \isTrep(\nu, g)$ if
$g$ is a probability measure on $\reals_+$ such that $\nu$ is absolutely continuous with respect to $g$,
i.e.~$\nu \ll g$, and
\[
\Theta &= \sum_{k=1}^{\infty}\theta_{k}\delta_{\psi_{k}}, & 
&\text{with} &
\theta_{k} &= V_k\ind\left(\frac{\dee \nu}{\dee g}(V_k) \geq \Gamma_{k}\right), &
V_{k} &\distiid g.\label{eq:thinning}
\]
\citet{Rosinski:1990} showed that $\Theta\isTrep(\nu,g)$ implies $\Theta\dist\distCRM(\nu)$.
Note that $\Gamma_{k} \convas \infty$ as $k\to\infty$, so the probability that $\frac{\dee \nu}{\dee g}(V_k) \geq \Gamma_{k}$
is decreasing in $k$. Thus, this representation generates atoms with $\theta_k = 0$ (which have no effect and can be removed) 
increasingly frequently and becomes inefficient as $k\to\infty$.

\begin{exa}[Thinning representation for $\distGammaP(\gamma, \lambda, d)$] \label{ex:t-rep-gamma-process}
If we let $g = \distGam(1 - d, \lambda)$, then the thinning representation for $\distGammaP(\gamma, \lambda, d)$ is 
\[
\Theta &= \sum_{k=1}^{\infty} V_{k}\ind(V_{k}\Gamma_{k} \le \gamma)\delta_{\psi_{k}}, &
&\text{with} &
V_{k} &\distiid \distGam(1 - d, \lambda).
\]

\end{exa}

\paragraph{Rejection \citep{Rosinski:2001}}
Using the nomenclature of \citet{Rosinski:2001}, we say $\Theta$ has a \emph{rejection representation} and write $\Theta \isRrep(\nu, \mu)$ if
$\mu$ is a measure on $\reals_+$ satisfying \cref{eq:nuassump} and $\frac{\dee \nu}{\dee \mu} \leq 1$, and
\[
\Theta &= \sum_{k=1}^{\infty}\theta_{k}\delta_{\psi_{k}},  &
 &\text{with} &
 \theta_{k} &= V_k\ind\left(\frac{\dee\nu}{\dee\mu}(V_k) \geq U_k \right), &
 (V_{k})_{k\in\nats} &\dist \distPP(\mu),   \label{eq:rejection} \\
&&&& 
  U_{k} &\distiid \distUnif[0, 1]. 
\]
\citet{Rosinski:2001} showed that $\Theta\isRrep(\nu,\mu)$ implies $\Theta\dist\distCRM(\nu)$.
This representation is very similar to the thinning representation, except that the sequence
$(V_k)_{k\in\nats}$ is generated from a Poisson process on $\reals_+$ rather than \iid 
This allows $V_k \convas 0$ as $k\to\infty$, causing the frequency of generating
ineffective atoms $\theta_k = 0$ to decay as $k\to\infty$,
assuming $\mu$ is appropriately chosen such that $\frac{\dee\nu}{\dee\mu}(\theta) \to 1$ as $\theta\to0$.
This representation can thus be constructed to be more efficient than the thinning representation.
We can calculate the efficiency in terms of the expected number of rejections
(that is, the number of $\theta_{k}$ that are identically zero):
\bnprop \label{prop:r-rep-efficiency}
For $\Rrep(\nu, \mu)$, the expected number of rejections is
\[
\EE\left[\sum_{k=1}^{\infty}\ind(\theta_{k} = 0)\right] = \int \left(1 - \frac{\dee\nu}{\dee\mu}(x)\right) \mu(\dee x). 
\]
\enprop
\brmk
If $\mu$ and $\nu$ can be written as densities with respect to Lebesgue measure, then the
integral in \cref{prop:r-rep-efficiency} can be rewritten as
$\int (\mu(x) - \nu(x))\dee x$. 
\ermk

\begin{exa}[Rejection representation for $\distGammaP(\gamma, \lambda, 0)$] \label{ex:r-rep-gamma-process}
Following \citet{Rosinski:2001}, consider $\mu(\dee \theta) = \gamma \lambda \theta^{-1}(1 + \lambda \theta)^{-1}\dee \theta$.
We call $\distCRM(\mu)$ the Lomax process, $\distLomP(\gamma, \lambda^{-1})$, after the related Lomax distribution. 
We can use the inverse-L\'evy method analytically with $\mu$ since $\muinv(u) = \frac{1}{\lambda(e^{(\gamma\lambda)^{-1}u} - 1)}$.
Thus, the rejection representation of $\distGammaP(\gamma, \lambda, 0)$ is
\[
\Theta &= \sum_{k=1}^{\infty} V_{k} \ind(U_{k} \le (1 + \lambda V_{k})e^{-\lambda V_{k}})\delta_{\psi_{k}}, & 
&\text{with} &
V_{k} &= \frac{1}{\lambda(e^{(\gamma\lambda)^{-1}\Gamma_{k}} - 1)}, &
 U_{k} &\distiid \distUnif[0,1].
\]
Unlike in the thinning construction given in \cref{ex:t-rep-gamma-process},  
only a finite number of rates will be set to zero almost surely. 
In particular, the expected number of rejections is $\gamma \lambda c_{\gamma}$, where $c_{\gamma}$ is the Euler-Mascheroni constant. 
\end{exa}

\begin{exa}[Rejection representation for $\distGammaP(\gamma, \lambda, d)$, $d>0$] \label{ex:r-rep-gamma-process-power-law}
For the case of $d > 0$, we instead use $\mu(\dee \theta) = \gamma\frac{\lambda^{1-d}}{\Gamma(1-d)} \theta^{-1-d}\dee \theta$. 
We can again use the  inverse-L\'evy method analytically with $\mu$ since $\muinv(u) = (\gamma'u^{-1})^{1/d}$, where 
$\gamma' \defined \gamma\frac{\lambda^{1-d}}{d\Gamma(1-d)}$. 
The rejection representation is then
\[
\Theta &= \sum_{k=1}^{\infty} V_{k} \ind(U_{k} \le e^{-\lambda V_{k}})\delta_{\psi_{k}}, & 
&\text{with} &
V_{k} &=  (\gamma'\Gamma_{k}^{-1})^{1/d}, &
 U_{k} &\distiid \distUnif[0,1].
\]
The expected number of rejections is $\gamma\frac{\lambda^{1-d}}{d}$, so the representation is 
efficient for large $d$, but extremely inefficient when $d$ is small. 
\end{exa}

\subsection{Superposition representations} 
\label{sec:superpositionreps}

Superposition representations arise as an infinite sum of CRMs
with finite rate measure. 
These tend to be easier to analyze than series representations 
as they decouple atoms between the summed CRMs, but can produce representations
with larger truncation error (cf.~\cref{sec:trunc,sec:simulation}).  
Throughout, let $\nu$ be a measure on $\reals_+$ satisfying the basic conditions in \cref{eq:nuassump},
and let $\psi_k \distiid G$.

\paragraph{Decoupled Bondesson}
We say $\Theta$ has a \emph{decoupled Bondesson representation}
and write $\Theta \isDBFrep(c, g, \xi)$ if for $c > 0, \xi > 0$, and $g$ a density on $\reals_{+}$,
\vphantom{\cref{eq:f-rep}}%
\[
\Theta &= \sum_{k=1}^{\infty}\sum_{i=1}^{C_{k}}\theta_{ki}\delta_{\psi_{ki}}, &
&\text{with} &
C_{k} &\distiid \distPoiss(c/\xi), &
\theta_{ki} &= V_{ki} e^{-T_{ki}}, \label{eq:f-rep}\\ 
&& && T_{ki} &\distind \distGam(k, \xi), &
V_{ki} &\distiid g. \nonumber
\]
This is a novel superposition representation, though special cases are already known~\citep{Paisley:2010,Roychowdhury:2015}.
\cref{prop:decoupled-bondesson-representation} shows that the decoupled Bondesson
representation applies to the same class of CRMs as the Bondesson representation from \cref{sec:seriesreps}. 
\begin{nthm}[Decoupled Bondesson representation] \label{prop:decoupled-bondesson-representation}
Let $\nu(\dee \theta) = \nu(\theta)\dee \theta$, $c_\nu$, and $g_\nu$ 
be as specified in \cref{prop:bondesson-representation}.
Then for any fixed $\xi > 0$,
\[
&\Theta\isDBFrep(c_{\nu}, g_{\nu}, \xi) \quad \text{implies} \quad \Theta \dist \distCRM(\nu).
\]
\end{nthm}
The proof of \cref{prop:decoupled-bondesson-representation} in \cref{app:seqproofs}
generalizes the arguments from \citet{Paisley:2010} and \citet{Roychowdhury:2015}. 
The free parameter $\xi$ controls the number of atoms generated for each outer sum index $k$; its principled selection
can be made by trading off computational complexity (cf.~\cref{sec:sampling}) and truncation
error (cf.~\cref{sec:trunc}). 

\begin{exa}[Decoupled Bondesson representation for $\distGammaP(\gamma, \lambda, 0)$]
Arguments paralleling those made in \cref{ex:b-rep-gamma-process}
show that the $\distGammaP(\gamma, \lambda, 0)$ representation from \citet{Roychowdhury:2015} follows 
directly from an application of \cref{prop:decoupled-bondesson-representation}:
if $\Theta \isDBFrep(\gamma\lambda, \distExp(\lambda), \xi)$, then $\Theta \dist 
\distGammaP(\gamma, \lambda, 0)$.
As in the Bondesson representation setting, \cref{prop:decoupled-bondesson-representation} does not apply to
$\distGammaP(\gamma, \lambda, d)$ when $d > 0$ because  the condition that $\theta \nu(\theta)$ is non-increasing 
fails to hold.
\end{exa}

\paragraph{Size-biased \citep{Broderick:2014,James:2014}}
Let $\pi(\theta)\defined h(0\given \theta)$.
We say $\Theta$ has a \emph{size-biased representation} and write $\Theta \isVrep(\nu, h)$ if
\[
\Theta &= \sum_{k=1}^\infty\sum_{i=1}^{C_{k}}\theta_{ki}\delta_{\psi_{ki}}, &
& \text{with} &
C_{k} &\distind \distPoiss\left(\eta_{k}\right), \\
&&&& \theta_{ki} &\distind \frac{1}{\eta_{k}}\pi(\theta)^{k-1}\left(1-\pi(\theta)\right)\nu(\dee\theta), \label{eq:crm_theta}\\
&&&& \eta_{k} &\defined \int \pi(\theta)^{k-1}\left(1-\pi(\theta)\right)\nu(\dee\theta).
\]
\citet{Broderick:2014} and \citet{James:2014} showed that $\Theta\isVrep(\nu, h)$ implies $\Theta\dist\distCRM(\nu)$.
If the rate measure $\nu$ and the likelihood $h$ are selected to be a conjugate
exponential family then, noting that 
$\sum_{x=1}^\infty h(x\given\theta) = 1-\pi(\theta)$, 
the rate $\theta_{ki}$ can be sampled
from a mixture of exponential family distributions:
\[
\theta_{ki}\given z_{ki} &\distind \frac{1}{\eta_{kz_{ki}}}h(z_{ki} \given \theta) \pi(\theta)^{k-1}\nu(\dee\theta), &
z_{ki} &\distind \distCat\left( \left(\eta_{kx}/\eta_k\right)_{x=1}^\infty \right),  \\
\eta_{kx} &\defined \int h(x\given \theta)\pi(\theta)^{k-1}\nu(\dee\theta).
\]

\begin{exa}[Size-biased representation for $\distGammaP(\gamma, \lambda, d)$]
For the Gamma process, values for $\eta_{kx}$ and $\eta_k$ can be found using integration by parts and the standard gamma distribution integral,
while $\theta_{ki} \given z_{ki}$ is sampled from a gamma distribution by inspection:
\[
\eta_{kx} &= \frac{\gamma \lambda^{1-d}\Gamma(x - d)}{x!\Gamma(1-d)(\lambda + k)^{x-d}}, &
\eta_k &= \left\{\begin{array}{ll}
  \frac{\gamma\lambda^{1-d}}{d}\left((\lambda+k)^d-(\lambda+k-1)^d\right) & d > 0\\
  \gamma\lambda\left(\log(\lambda+k)-\log(\lambda+k-1)\right) & d = 0
\end{array}\right.,\\
\theta_{ki} &\given z_{ki} \distind \distGam(x - d, \lambda + k).
\] 
\end{exa}

\paragraph{Power-law}
We say $\Theta$ has a \emph{power-law representation} and write $\Theta\isPLFrep(\gamma, \alpha, d, g)$
if for $\gamma > 0, 0 \le d < 1$, $\alpha > -d$, and $g$ a density on $\reals_{+}$, 
\vphantom{\cref{eq:plrep}}
\[
\Theta &= \sum_{k=1}^{\infty} \sum_{i=1}^{C_{k}} \theta_{ki} \delta_{\psi_{ki}}, &
&\text{with}&
C_{k} &\distiid \distPoiss(\gamma), &
\theta_{ki} &= V_{ki}U_{kik}\prod_{j=1}^{k-1}(1 - U_{kij}), \qquad \label{eq:plrep} \\
&& && V_{ki} &\distiid g,  &
U_{kij} &\distind \distBeta(1 - d, \alpha + j d). \nonumber
\]
This is a novel superposition representation, although it was previously developed in 
the special case of the beta process (where $g(v) = \delta_1$)~\citep{Broderick:2012}.
The name of this representation arises from the fact that it exhibits 
Types I and II power-law behavior \citep{Broderick:2012} under mild conditions when $d>0$, 
as we show in \cref{thm:f-rep-power-law-behavior} in the appendix (note, however, that it will not exhibit power-law behavior when $d=0$).
\cref{thm:power-law-rep} below shows the conditions under which $\Theta\isPLFrep(\gamma, \alpha, d, g)$ implies $\Theta \dist\distCRM(\nu)$.
Its proof in \cref{app:seqproofs} relies on the notion of \emph{stochastic mapping} (\cref{lem:stochastic-mapping}), a powerful technique for transforming
one CRM into another.
Note that in \cref{lem:stochastic-mapping}, the case where $u$ is a deterministic function of $\theta$
via the mapping $u = \tau(\theta)$ may be recovered by setting $\kappa(\theta, \dee u) = \delta_{\tau(\theta)}$. 
\begin{nlem}[CRM stochastic mapping] \label{lem:stochastic-mapping}
Let $\Theta = \sum_{k=1}^{\infty} \theta_{k}\delta_{\psi_{k}} \dist \distCRM(\nu)$.
Then for any probability kernel $\kappa(\theta, \dee u)$, we have $\kappa(\Theta)\dist\distCRM(\nu_{\kappa})$, where
\[
\kappa(\Theta) &\defined \sum_{k=1}^{\infty} u_{k}\delta_{\psi_{k}}, &
u_{k}\given \theta_k &\dist \kappa(\theta_{k}, \cdot), &
&\text{and} &
\nu_{\kappa}(\dee u) &\defined \int \kappa(\theta, \dee u)\nu(\dee \theta).
\]
\end{nlem}

\begin{nthm}[Power-law representation] \label{thm:power-law-rep}
Let $\nu(\dee\theta) = \nu(\theta)\dee\theta$ be a rate measure satisfying \cref{eq:nuassump}, and let $g_\nu$ be a density on $\reals_+$ such that
\[
\nu(u) &= \int \theta^{-1} g_\nu\left(u \theta^{-1}\right)\nu_{\distBP}\left(\dee \theta\right),\\
\shortintertext{where}
\nu_{\distBP}(\dee\theta) &= \gamma\frac{\Gamma(\alpha + 1)}{\Gamma(1 - d)\Gamma(\alpha +
d)}\ind\left[\theta \leq 1\right] \theta^{-1-d}(1-\theta)^{\alpha+d-1}\dee\theta
\]
is the rate measure for the beta process $\distBP(\gamma,\alpha,d)$ from \cref{eq:bpmsr}. Then
\[
\Theta \isPLFrep(\gamma,\alpha,d,g_\nu) \quad \text{implies} \quad \Theta\dist\distCRM(\nu).
\]
\end{nthm}

\begin{exa}[Power-law representation for $\distGammaP(\gamma, \lambda, d)$] \label{ex:power-law-f-rep-gamma-process}
If we choose $g_{\nu} = \distGam(\lambda, \lambda)$, then 
using the change of variable $w = u(\theta^{-1} - 1)$,
\[
\lefteqn{\int \theta^{-1} g_\nu\left(u \theta^{-1}\right)\nu_{\distBP}\left(\dee \theta\right)} \\
&=  \gamma\lambda\frac{\lambda^{\lambda}}{\Gamma(1 - d)\Gamma(\lambda + d)}u^{\lambda-1}\int \theta^{-\lambda-d-1}e^{-\lambda u\theta^{-1}}(1-\theta)^{\lambda+d-1}\,\dee\theta\, \dee u \\
&= \gamma\lambda\frac{\lambda^{\lambda}}{\Gamma(1 - d)\Gamma(\lambda + d)}u^{-1-d}e^{-\lambda u}\int w^{\lambda+d-1}e^{-\lambda w}\,\dee w\,\dee u \\
&= \gamma\frac{\lambda^{1-d}}{\Gamma(1 - d)}u^{-1-d}e^{-\lambda u}\,\dee u. %
\]
It follows immediately from \cref{thm:power-law-rep} that if
$\Theta \isPLFrep(\gamma, \lambda, d, \distGam(\lambda, \lambda))$, then
$\Theta \dist \distGammaP(\gamma, \lambda, d)$.
To the best knowledge of the authors, this power-law representation for the gamma process is novel. 
\end{exa}

\section{Truncation analysis}
\label{sec:trunc}

Each of the sequential representations developed in \cref{sec:seqcrms} shares
a common structural element---an outer infinite sum---which is responsible
for generating a countably infinite number of atoms in the CRM.
In this section, we terminate these outer sums at a finite \emph{truncation level} $K\in\nats$,
resulting in a \emph{truncated CRM} $\Theta_K$ possessing a finite number of atoms.
We develop upper bounds on the error induced by this truncation procedure. 
All of the truncated CRM error bounds in this section rely on \cref{lem:protobound-crm}, which
is a tightening (by a factor of two) of 
the bound in~\citet{Ishwaran:2001,Ishwaran:2002b} (for its generalization to arbitrary discrete random measures, see \cref{lem:protobound}). 

\begin{nlem}[CRM protobound]\label{lem:protobound-crm}
Let $\Theta\dist \distCRM(\nu)$. For any truncation $\Theta_K$, if
\begin{align}
X_{n} \given \Theta &\distiid \distLP(h, \Theta), & 
Z_{n} \given \Theta_{K} &\distiid \distLP(h, \Theta_K), \\
Y_{n} \given X_{n} &\distind f(\cdot \given X_{n}), & 
W_{n} \given Z_{n} &\distind f(\cdot \given Z_{n}), 
\end{align}
  then, with $p_{N, \infty}$ and $p_{N, K}$ denoting the marginal densities of $Y_{1:N}$ and
  $W_{1:N}$, respectively,
  \begin{align}
    \frac{1}{2}\|p_{N, \infty} - p_{N, K}\|_{1}
    &\leq 1-\Pr\left(\supp(X_{1:N}) \subseteq \supp(\Theta_K)\right),
  \end{align}
\end{nlem}

The proof of all results in this section (including \cref{lem:protobound-crm}) can be found in \cref{app:truncproofs}.
All of the provided truncation results use the generative model in \cref{lem:protobound-crm},
and are summarized in \cref{tbl:results} in \cref{sec:simulation}.
Throughout this section, for
a given likelihood model $h(x \given \theta)$ we define $\pi(\theta) \defined h(0\given \theta)$ for notational brevity. 
The asymptotic behavior of truncation error bounds is specified with tilde notation:
\begin{align}
 a(K) \sim b(K), \quad K\to\infty  \quad \iff\quad \lim_{K\to\infty} \frac{a(K)}{b(K)} = 1. 
\label{eq:asympdef}
\end{align} 

\subsection{Series representations}
\label{sec:seriestrunc}

Each of the series representations can be viewed a functional of a standard
Poisson point process and a sequence of \iid\ random variables with some 
distribution $g$ on $\reals_+$.
In particular, we may write each in the form 
\[
\Theta &= \sum_{k=1}^{\infty}\theta_{k}\delta_{\psi_{k}}, & 
&\text{with} &
\theta_{k} &= \tau(V_{k}, \Gamma_{k}), &
V_{k} &\distiid g, \label{eq:general-tau-series-rep}
\]
where $\Gamma_k$ are the jumps of a unit-rate homogeneous Poisson point process on $\reals_+$, 
and $\tau : \reals_+\times \reals_+\to\reals_+$ is a non-negative measurable function
such that $\lim_{u\to\infty}\tau(v, u) = 0$ for $g$-almost every $v$. 
The truncated CRM then takes the form 
\[
\Theta_{K} \defined \sum_{k=1}^{K}\theta_{k}\delta_{\psi_{k}}.
\]
\cref{thm:series-rep-trunc} provides a general truncation error bound for series representations of 
the form  \cref{eq:general-tau-series-rep}, specifies its range, and guarantees that the
bound decays to 0 as $K \to \infty$.
\bnthm[Series representation truncation error] \label{thm:series-rep-trunc}
The error in approximating a series representation of $\Theta$ with its 
truncation $\Theta_K$ satisfies
\[
0 &\leq \frac{1}{2}\|p_{N,\infty}-p_{N,K}\|_1 \leq 1 - e^{-B_{N,K}} \leq 1,
\shortintertext{where}
B_{N,K} &\defined \int_{0}^{\infty}\left(1 - \EE\left[\pi\left(\tau\left(V, u + G_{K}\right)\right)^N\right]\right)\,\dee u, \label{eq:series-rep-trunc-inequality}
\]
$G_{0} \defined 0$, $G_{K} \distind \distGam(K, 1)$ for $K \ge 1$, and $V\distind g$. 
Furthermore, $\forall N\in\nats, \, \lim_{K\to\infty}B_{N,K} = 0$.
\enthm
\brmk
An alternate form of $B_{N, K}$ that is sometimes easier to use in practice 
can be found by applying the standard geometric series formula to \cref{eq:series-rep-trunc-inequality},
which yields
\[
 B_{N, K}  = \sum_{n=1}^N\int_{0}^{\infty}\EE\left[\pi\left(\tau\left(V, u + G_{K}\right)\right)^{n-1}\left(1 - \pi\left(\tau\left(V, u + G_{K}\right)\right)\right)\right]\,\dee u.\label{eq:series-geobound}
\]
A simplified upper bound on $B_{N, K}$ can be derived by noting that $\pi(\theta)\leq 1$, so
\[
  B_{N,K} \leq N\int_{0}^{\infty}\left(1 - \EE\left[\pi\left(\tau\left(V, u + G_{K}\right)\right)\right]\right)\,\dee u.\label{eq:series-ezbound}
\]
This bound usually gives the same asymptotics in $K$ as \cref{eq:series-rep-trunc-inequality}.
\ermk
The main task in using \cref{thm:series-rep-trunc} to develop a truncation error bound for a series representation
is evaluating the integrand in the definition of $B_{N, K}$. 
    Thus, we next evaluate the integrand and provide expressions 
of the truncation error bound for the four series representations outlined in \cref{sec:seriesreps}.
Throughout the remainder of this section, $G_K$ is defined as in \cref{thm:series-rep-trunc}, $F_0 \equiv 1$, and $F_K$ is the CDF of $G_K$. 

\paragraph{Inverse-L\'evy representation}
For this representation we have
\[
\tau(v, u) = \nuinv(u) \defined \inf\left\{y : \nu\left([y, \infty)\right) \leq u\right\}.
\]
To evaluate the bound in \cref{eq:series-ezbound}, we use  the transformation of variables 
$x = \nuinv(u+G_K)$ and the fact that for $a, b \geq 0$, $\nuinv(a) \ge b \iff a \le \nu\left(\left[b, \infty\right)\right)$
to conclude that
\[
B_{N,K} \leq N\int_{0}^{\infty}F_{K}(\nu[x, \infty))(1 - \pi(x))\,\nu(\dee x). \label{eq:il-trunc-bound}
\]
Recent work on the inverse-L\'evy representation has developed Monte Carlo estimates of the error 
of the truncated random measure moments for those $\nu\left([x, \infty)\right)$ with
known inverse $\nuinv$~\citep{Arbel:2015}. In contrast, the result above
provides an explicit bound on the $L^1$ truncation error. Our bound does not require
knowing $\nuinv$, which is often the most challenging 
aspect of applying the inverse-L\'evy representation. 

\begin{exa}[\ILrep~truncation for $\distLomP(\gamma, \lambda^{-1})$ with Poisson likelihood] \label{ex:il-lomax-trunc}
Recall from \cref{ex:r-rep-gamma-process} that the Lomax process $\distLomP(\gamma, \lambda^{-1})$
is the CRM with rate measure  $\nu(\dee \theta) = \gamma \lambda \theta^{-1}(1 + \lambda \theta)^{-1}\dee \theta$,
so $\nu[x, \infty) = \gamma\lambda \log\{1 + (\lambda x)^{-1}\}$. 
Using \cref{eq:il-trunc-bound}, we have
\[
B_{N,K} 
&\le N \gamma\lambda \int_{0}^{\infty} F_{K}(\gamma\lambda \log\{1 + (\lambda x)^{-1}\})(1 - e^{-x})x^{-1}(1 + \lambda x)^{-1}\dee x.
\]
Since $F_{K}(t) \le t^{K}/K! \le (3 t / K)^{K}$, for any $a > 0$ the integral is upper bounded by
\[
\int_{0}^{a}(1 + \lambda x)^{-1} \dee x &+ F_{K}\left(\gamma\lambda \log\left\{1 + \frac{1}{\lambda a}\right\}\right)\int_{a}^{\infty}x^{-1}(1 + \lambda x)^{-1}\dee x \\
&\le a + F_{K}\left(\gamma\lambda \log\left\{1 + \frac{1}{\lambda a}\right\}\right) \log\left\{1 + \frac{1}{\lambda a}\right\} \\
&\le \lambda^{-1}(e^{b} - 1)^{-1} + b (3 \gamma\lambda b / K)^{K} \quad \text{where} \quad b \defined \log\{1 + (\lambda a)^{-1}\}. \label{eq:il-rep-lomp-intermediate-bound}
\]
Replacing $(e^{b} - 1)^{-1}$ with the approximation $e^{-b}$ and then setting the two terms in \cref{eq:il-rep-lomp-intermediate-bound} equal,
we obtain $b = K W_{0}\left(\left\{3\gamma\lambda^{\frac{K+2}{K+1}}(K+1)^{\frac{1}{K+1}}\right\}^{-1}\right)$, where $W_{0}$ 
is the product logarithm function, i.e.~
\[
W_0(y) = x \iff x e^x = y.\label{eq:lambertdefn}
\] 
Thus, using the fact that $e^{-t} \le (e^{t} - 1)^{-1}$ and $\lambda^{\frac{K+2}{K+1}}(K+1)^{\frac{1}{K+1}}$ reaches its maximum at $K = \max(0, e\lambda^{-1}-1)$, we conclude that
\[
B_{N,K} 
&\le \frac{2 N \gamma[1 + (3\gamma \lambda)^{-1}]}{\exp\left(K W_{0}\left(\left\{3\gamma\lambda\max(\lambda, e)\right\}^{-1}\right)\right) - 1} \\
&\sim 2 N \gamma[1 + (3\gamma \lambda)^{-1}]e^{-K W_{0}(\{3\gamma\lambda\max(\lambda, e)\}^{-1})} \qquad K \to \infty.
\]
\end{exa}

\paragraph{Bondesson representation}
For this representation we have
\[
\tau(v, u) &= v e^{-u/c}, &
g(\dee v) &= -c^{-1}\frac{\dee}{\dee v}\left(v\nu(v)\right)\,\dee v.
\]
Writing the expectation over $V$ explicitly as an integral with measure $g(v)\dee v$, 
using the transformation of variables $u = -c \log x/v$ (so $x = ve^{-u/c}$), and 
given the definition of $g(v) = -c^{-1}\frac{\dee}{\dee v}\left(v\nu(v)\right)$ for the Bondesson representation, we have 
\[
B_{N,K} \leq N \int_{0}^{\infty} \left(1 - \EE\left[\pi\left(ve^{-G_{K}/c}\right)\right]\right)\nu(\dee v).
\]

\begin{exa}[Truncation of the Bondesson representation for $\distGammaP(\gamma, \lambda, 0)$] \label{ex:d-rep-gpp-trunc}
Let $\tG_K \eqD G_K/(\gamma\lambda)$. 
Since $\pi(\theta) = e^{-\theta}$ and $c = \gamma\lambda$, we have
\[
\int_{0}^{\infty} \left(1 - \EE\left[\pi(ve^{-G_{K}/c})\right]\right)\nu(\dee v)
&= \gamma\lambda\,\EE\left[\int_{0}^{\infty}(1 - e^{-ve^{-\tG_{K}}})v^{-1}e^{-\lambda v} \dee v\right] \\
&= \gamma\lambda\,\EE\left[\log(1 + e^{-\tG_{K}}/\lambda)\right]  \\
&\le \gamma\,\EE\left[e^{-\tG_{K}}\right]  
= \gamma \left(\frac{\gamma\lambda}{1 + \gamma\lambda}\right)^{K}.
\]
The second equality follows by using the power series for the exponential 
integral~\citep[Chapter 5]{Abramowitz:1964}. 
Thus,
\[
B_{N,K} 
&\le N \gamma\left(\frac{\gamma\lambda}{1 + \gamma\lambda}\right)^{K}.
\]
\end{exa}

\paragraph{Thinning representation}
For this representation we have
\[
\tau(v, u) &= v \ind\left[\frac{\dee\nu}{\dee g}(v) \ge u\right], &
& g\text{ any distribution on $\reals_+$ s.t. } \nu \ll g.
\]
Since $\pi(0) = 1$ by \cref{lem:behavior-at-zero}, we have that 
$1-\pi\left(v\ind(A)\right) = \left(1-\pi\left(v\right)\right)\ind(A)$ for any event $A$.
Using this fact, we have
\[
B_{N,K} \leq N\int_{0}^{\infty}(1 - \pi(v))\int_{0}^{\frac{\dee \nu}{\dee g}(v)}F_{K}\left(u\right)\dee u\,g(\dee v). \label{eq:thinningrepbd}
\]
Analytic bounds for the thinning representation of specific processes tend to be opaque and notationally cumbersome, so
we simply compare its truncation error in \cref{sec:simulation} to the
other representations by numerical approximation of \cref{eq:thinningrepbd}.

\paragraph{Rejection representation}
Assume that we can use the inverse-L\'evy representation to simulate $\distPP(\mu)$. 
Then for the rejection representation we have
\[
\tau(v, u) &= \muinv(u)\ind\left[\frac{\dee\nu}{\dee \mu}(\muinv(u)) \ge v \right], &
g(\dee v) &= \ind[0 \le v \le 1] \dee v,
\]
where $\mu$ satisfies $\nu \ll \mu$,
$\frac{\dee\nu}{\dee \mu} \le 1$, and $\muinv(u) \defined \inf \left\{x : \mu\left([x, \infty)\right) \leq u\right\}$.
Using the same techniques as for the thinning and inverse-L\'evy representations, we have that
\[
B_{N,K} \leq N\int_{0}^{\infty}F_{K}(\mu[x, \infty))(1 - \pi(x))\,\nu(\dee x). \label{eq:rej-trunc-bound}
\]

\begin{exa}[\Rrep~truncation for $\distGammaP(\gamma, \lambda, 0)$ with Poisson likelihood]\label{ex:rrepgamp}
Using \cref{eq:rej-trunc-bound} and the fact that $1 - e^{-x} \le x$, we have
\[
B_{N,K} 
&\le N \gamma\lambda \int_{0}^{\infty}F_{K}(\gamma\lambda \log\{1 + (\lambda x)^{-1}\})e^{-\lambda x}\dee x. \label{eq:r-rep-gamp-formal-bound}
\]
Arguing as in \cref{ex:il-lomax-trunc}, we see that the integral in \cref{eq:r-rep-gamp-formal-bound} is upper bounded by
\[
\int_{0}^{a} e^{-\lambda x}\dee x &+ F_{K}\left(\gamma\lambda \log\left\{1 + \frac{1}{\lambda a}\right\}\right)\int_{a}^{\infty}e^{-\lambda x}\dee x \\
&\le a + \lambda^{-1}F_{K}\left(\gamma\lambda \log\left\{1 + \frac{1}{\lambda a}\right\}\right) \\
&= \lambda^{-1}\left((e^{b} - 1)^{-1} + (3 \gamma\lambda b / K)^{K}\right), \label{eq:r-rep-gamp-intermediate-bound}
\]
where $b \defined \log\{1 + (\lambda a)^{-1}\}$. 
Replacing $(e^{b} - 1)^{-1}$ with the approximation $e^{-b}$ and then setting the two terms in \cref{eq:r-rep-gamp-intermediate-bound} equal
to each other, we obtain $b = K W_{0}(\{3\gamma\lambda\}^{-1})$ (where $W_0$ is defined in \cref{eq:lambertdefn})
and conclude that
\[
B_{N,K} \le \frac{2 N \gamma}{e^{K W_{0}(\{3\gamma\lambda\}^{-1})} - 1}
\sim 2 N \gamma e^{-K W_{0}(\{3\gamma\lambda\}^{-1})} \qquad K \to \infty.
\]
\end{exa}

\begin{exa}[\Rrep~truncation for $\distGammaP(\gamma, \lambda, d)$ with Poisson likelihood, $d>0$]\label{ex:rrepgamp-power-law}
We have
\[
B_{N,K} 
&\le N \frac{\gamma \lambda^{1-d}}{\Gamma(1-d)} \int_{0}^{\infty} F_{K}(\gamma'x^{-d})(1 - e^{-x})x^{-1-d}e^{-\lambda x}\dee x.
\]
The integral can be upper bounded as
\[
\lefteqn{\int_{0}^{a}x^{-d}\dee x +  F_{K}(\gamma'a^{-d})\int_{a}^{\infty}(1+e^{-x}) x^{-1-d}e^{-\lambda x}\dee x} \\
&\le (1-d)^{-1}a^{1-d} +  \Gamma(-d)(\lambda^{d} - (1+\lambda)^{d}) (3\gamma' K^{-1}a^{-d})^{K}.
\]
Setting the two terms equal and solving for $a$, we obtain
\[
B_{N,K} 
&\le 2 N \frac{\gamma\lambda^{1-d}}{\Gamma(2-d)}[(1-d)\Gamma(-d)]^{\frac{d(1-d)}{d(1-d)+K}}\left[\frac{3\gamma\lambda^{1-d}}{d\Gamma(1-d) K}\right]^{\frac{Kd(1-d)}{d(1-d)+K}} \\
&\sim 2 N \frac{\gamma\lambda^{1-d}}{\Gamma(2-d)}\left[\frac{3\gamma\lambda^{1-d}}{d\Gamma(1-d)}\right]^{d(1-d)} K^{-d(1-d)} \qquad K \to \infty.
\]
\end{exa}

\subsection{Superposition representations} \label{sec:superpositiontrunc}

For superposition representations, the truncated CRM takes the form
\[
\Theta_{K} &\defined \sum_{k=1}^{K}\sum_{i=1}^{C_{k}}\theta_{ki}\delta_{\psi_{ki}}.
\]
Let $\Theta_{K}^{+} \defined \Theta - \Theta_{K}$
denote the \emph{tail measure}. 
By the superposition property of Poisson point processes \citep{Kingman:1993}, the tail measure
is itself a CRM with some rate measure $\nu_{K}^{+}$ and is independent of $\Theta_{K}$:
\[
\Theta_K^+ &= \sum_{k=K+1}^\infty\sum_{i=1}^{C_{k}}\theta_{ki}\delta_{\psi_{ki}} \dist\distCRM\left(\nu_K^+\right), &
\Theta_K^+ &\indep \Theta_K, &
\Theta = \Theta_K+\Theta_K^+. \label{eq:superpositionsplitproperty}
\]
The following result provides a general truncation error bound for superposition representations, specifies its range, and guarantees
that the bound decays to 0 as $K \to \infty$.
\begin{nthm}[Superposition representation truncation error]\label{thm:superpositiontrunc}
The error in approximating a superposition representation of $\Theta\dist\distCRM(\nu)$ with its 
truncation $\Theta_K$ satisfies
\[
0 &\leq \frac{1}{2}\|p_{N,\infty}-p_{N,K}\|_1\leq 1 - e^{-B_{N,K}}\leq 1,
\shortintertext{where}
B_{N,K} &\defined \int \left(1-\pi(\theta)^N\right)\nu_K^+(\dee\theta).\label{eq:superposition-rep-trunc-inequality}
\]
Furthermore, $\forall N\in\nats, \, \lim_{K\to\infty}B_{N,K} = 0$.
\end{nthm}
\brmk
As for series representations, an alternate form of $B_{N, K}$ that is sometimes easier to use 
can be found by applying the standard geometric series formula to \cref{eq:superposition-rep-trunc-inequality}:
\[
 B_{N, K} = \sum_{n=1}^N\int\pi(\theta)^{n-1}\left(1-\pi(\theta)\right)\nu_K^+(\dee\theta).\label{eq:superposition-geobound}
\]
A simplified upper bound on $B_{N, K}$ can be derived by noting that $\pi(\theta)\leq 1$, so
\[
  B_{N,K} \leq N\int_{0}^{\infty}\left(1 - \pi\left(\theta\right)\right)\,\nu_K^+(\dee \theta).\label{eq:superposition-ezbound}
\]
This bound usually gives the same asymptotics in $K$ as \cref{eq:superposition-rep-trunc-inequality}.
\ermk
The main task in using \cref{thm:superpositiontrunc} to develop a truncation error bound for a superposition representation
is determining its tail measure $\nu_K^+$. In the following, we provide the tail measure for the three superposition representations
outlined in \cref{sec:superpositionreps}.

\paragraph{Decoupled Bondesson representation} 
For each point process in the superposition, an average of $c/\xi$ atoms are generated
with independent weights of the form $Ve^{-T_k}$ where $V\distind g$ and $T_k \distind \distGam(k, \xi)$.
Therefore, the tail measure is
\[
\nu_K^+(\dee\theta) &= 
\frac{c}{\xi}\sum_{k=K+1}^\infty \tg_{k,\xi}(\theta)\dee\theta,
\]
where $\tg_{k,\xi}$ is the density of $Ve^{-T_k}$.
The bound for the decoupled Bondesson representation can therefore be expressed as
\[
 B_{N, K} 
&\leq N\frac{c}{\xi}\sum_{k=K+1}^\infty\EE\left[1-\pi\left(Ve^{-T_k}\right)\right].
\]

\begin{exa}[Decoupled Bondesson representation truncation for $\distGammaP(\gamma,\lambda,0)$] \label{ex:db-rep-gpp-trunc}
Using the fact that $1 - e^{-\theta} \le \theta$, we have
\[
B_{N,K}
&= \frac{N\gamma\lambda}{\xi}\sum_{k=K+1}^{\infty}\EE[1 - \pi(V_{k1}e^{-T_{k1}})] 
\le \frac{N\gamma\lambda}{\xi}\sum_{k=K+1}^{\infty}\EE[V_{k1}e^{-T_{k1}}] \\
&=\frac{N\gamma\lambda}{\xi}\sum_{k=K+1}^{\infty}\frac{1}{\lambda}\left(\frac{\xi}{1+\xi}\right)^{k} 
= N\gamma\left(\frac{\xi}{1+\xi}\right)^{K}, 
\]
which is equivalent (up to a factor of 2) to the bound in \citet{Roychowdhury:2015}. 
\end{exa}

\paragraph{Size-biased representation} 
The constructive derivation of the size-biased representation~\citep[proof of Theorem 5.1]{Broderick:2014} immediately yields
\[
\nu_K^+(\dee\theta) = \pi(\theta)^K\nu(\dee\theta).
\]
Therefore, the size-biased representation truncation error bound can be expressed 
using the formula for $\eta_k$ from \cref{eq:crm_theta} as
\[
B_{N,K} &= \sum_{n=1}^N\int \pi(\theta)^{K+n-1}(1-\pi(\theta))\nu(\dee\theta) = \sum_{n=1}^N\eta_{K+n}. \label{eq:sb-trunc}
\]

\begin{exa}[Size-biased representation truncation for $\distGammaP(\gamma, \lambda, d)$]
For $d>0$, the standard gamma integral yields
\[
\eta_k = \int \pi(\theta)^{k-1}(1-\pi(\theta))\nu(\dee\theta)
&= \frac{\gamma \lambda^{1-d}}{d}\left((\lambda + k)^{d} - (\lambda + k - 1)^{d}\right). \label{eq:gprateint}
\]
The sum from \cref{eq:sb-trunc} is telescoping, so canceling terms,
\begin{align}
B_{N,K} 
&\le \frac{\gamma \lambda^{1-d}}{d}\left( (\lambda+K+N)^d - (\lambda+K)^d\right) 
\sim \gamma N \lambda^{1-d} K^{d-1} \qquad K \to \infty,
\end{align}
where the asymptotic result follows from \cref{lem:taylorlimitapprox}.
To analyze the $d=0$ case, we use L'Hospital's rule to take the 
limit of the integral:
\begin{align}
\lim_{d\to0}
\int \pi(\theta)^{k-1}(1-\pi(\theta))\nu(\dee\theta)
&=\gamma\lambda\left(\log(\lambda+k)-\log(\lambda+k-1)\right).
\end{align}
Canceling terms in the telescopic sum yields 
\begin{align}
B_{N,K} 
&\le \gamma \lambda\left( \log(\lambda+K+N) - \log(\lambda+K) \right)
\sim \gamma \lambda N K^{-1} \qquad K \to \infty,
\end{align}
where the asymptotic result follows from an application of \cref{lem:taylorlimitapprox}.
\end{exa}

\paragraph{Power-law representation} 
For each point process in the superposition, an average of $\gamma$ atoms are generated
with independent weights of the form $VU_k\prod_{\ell=1}^{k-1}(1-U_\ell)$, where $V\distind g$ and $U_\ell \distind \distBeta(1-d, \alpha+\ell d)$.
Therefore, the tail measure is
\[
\nu_K^+(\dee\theta) = \gamma\sum_{k=K+1}^\infty \tg_{k}(\theta)\dee\theta,
\]
where $\tg_k$ is the density of the random variable $VU_k\prod_{\ell=1}^{k-1}(1-U_\ell)$. 
The truncation error bound may be expressed as
\[
B_{N,K} \leq N\gamma \sum_{k=K+1}^\infty\EE\left[1-\pi\left(VU_k\prod_{\ell=1}^{k-1}(1-U_\ell)\right)\right].
\]

\begin{exa}[Power-law representation truncation for $\distGammaP(\gamma,\lambda,d)$] \label{ex:gpp-pl-rep-error}
Let $\beta_k$ be a random variable with density $\tg_k$ (with $\lambda$ in the place of $\alpha$). 
Using  $1 - e^{-\theta} \le \theta$, we have
\[
\sum_{k=K+1}^{\infty}\EE[1  - \pi(\beta_{k})]
&\le \sum_{k=K+1}^{\infty}\EE[\beta_{k}] 
= \EE\left[\sum_{k=K+1}^{\infty}\beta_{k}\right] 
= \prod_{k=1}^{K}\frac{\lambda + kd}{\lambda + kd - d + 1}, \label{eq:pyp-prod}
\]
where the final equality follows from \citet[Theorem 1]{Ishwaran:2001}. 
Thus, 
\[
\hspace{-.2cm}B_{N,K}
&\le \gamma N \prod_{k=1}^{K}\frac{\lambda + kd}{\lambda + kd - d + 1} 
\sim \gamma N \left\{\begin{array}{ll}
                      \left(\frac{\lambda}{\lambda+1}\right)^K & d = 0\\
                      \frac{\Gamma\left(\frac{\lambda+1}{d}\right)}{\Gamma\left(\frac{\lambda+d}{d}\right)}K^{1-d^{-1}} & 0<d<1
                      \end{array}\right. \quad K\to\infty,
\label{eq:gpp-pl-rep-error}
\]
where the $0<d<1$ case in \cref{eq:gpp-pl-rep-error} follows by \cref{lem:gautschi} applied to 
\[
\prod_{k=1}^{K}\frac{\lambda + kd}{\lambda + kd - d + 1}
&= \frac{\Gamma((\lambda + 1)/d)}{\Gamma((\lambda + d)/d)}\frac{\Gamma(\lambda/d + K+1)}{\Gamma(\lambda/d + K + d^{-1})}. %
\label{eq:pyp-prod-asymp}
\] 
\end{exa}

\subsection{Stochastic mapping
}\label{sec:stochmaptrunc}

We now show how truncation bounds developed elsewhere in this paper can be applied to CRM 
representations that have been transformed using \cref{lem:stochastic-mapping}.
For $\Theta \dist \distCRM(\nu)$, we denote its transformation by 
$\tTheta = \kappa(\Theta)$. 
For any object defined with respect to $\Theta$, the corresponding object for $\tTheta$ is denoted 
with a tilde.
For example, in place of $N$ and $X_{1:N}$ (for $\Theta$), 
we use $\tN$ and $\tX_{1:\tN}$ (for $\tTheta$). 
We make $B_{N,K}$ a function of $\pi(\theta)$ in the notation of
\cref{prop:stochastic-mapping-trunc}; when one applies stochastic mapping
to a CRM, one usually also wants to change the likelihood $h(x\given \theta)$,
and thus also changes $\pi(\theta) = h(0\given\theta)$.
The proof of \cref{prop:stochastic-mapping-trunc} may be found in \cref{app:truncproofs}.

\bnprop[Truncation error under a stochastic mapping] \label{prop:stochastic-mapping-trunc}
Consider a representation for $\Theta \dist \distCRM(\nu)$ with truncation
error bound $B_{N,K}(\pi)$. Then for any likelihood $\tilde{h}(x \given u)$,
if $\tTheta$ is a stochastic mapping of $\Theta$ under the probability kernel $\kappa(\theta, \dee u)$,
 its truncation error bound is $B_{1,K}(\pi_{\kappa,\tN})$, where
$\pi_{\kappa,\tN}(\theta) \defined \int \tilde{h}(0 \given u)^{\tN} \kappa(\theta, \dee u)$.
\enprop

\subsection{Hyperpriors}

In practice, prior distributions are often placed on the hyperparameters of the CRM rate measure (i.e.~$\gamma$, $\alpha$, $\lambda$, $d$, etc.).
We conclude our investigation of CRM truncation error by showing how bounds developed in this section can be modified
to account for the use of hyperpriors. 
Note that we make the dependence of $B_{N,K}$ on the hyperparameters $\Phi$ 
explicit in the notation of \cref{prop:trunccrm-hyperprior}.
\bnprop[CRM truncation error with a hyperprior] \label{prop:trunccrm-hyperprior}
Given hyperparameters $\Phi$, consider a representation for $\Theta \given \Phi \dist \distCRM(\nu)$,
and let $B_{N,K}(\Phi)$ be given by \cref{eq:series-rep-trunc-inequality} (for a series representation) or \cref{eq:superposition-rep-trunc-inequality} (for a superposition representation).
The error of approximating $\Theta$ with its truncation $\Theta_K$ satisfies
\begin{align}
0\leq \frac{1}{2}\|p_{N,\infty}&-p_{N,K}\|_1 \leq  1- e^{-\EE\left[B_{N,K}(\Phi)\right]}\leq 1.
\end{align}
\enprop

\begin{exa}[Decoupled Bondesson representation truncation for $\distGammaP(\gamma,\lambda,0)$]
A standard choice of hyperprior for the mass $\gamma$ is a gamma distribution, i.e.~$\gamma \dist \distGam(a, b)$.
Combining \cref{prop:trunccrm-hyperprior} and \cref{ex:db-rep-gpp-trunc}, we have that
\[
\EE\left[B_{N,K}\left(\Phi\right) \right]
&\le N\frac{a}{b}\left(\frac{\xi}{\xi+1}\right)^K.
\]
\end{exa}

\section{Normalized truncation analysis} \label{sec:ncrm-trunc}

In this section, we provide truncation error bounds for \emph{normalized CRMs} (NCRMs). 
Examples include the Dirichlet process~\citep{Ferguson:1973},
the normalized gamma process~\citep{Brix:1999,James:2002,Lijoi:2003,Pitman:2003,Lijoi:2007,Lijoi:2010},
and the normalized $\sigma$-stable process~\citep{Kingman:1975,Lijoi:2010}.
Given a CRM $\Theta$ on $\Psi$, we define 
the corresponding NCRM $\Xi$ via $\Xi(S) \defined
\Theta(S)/\Theta(\Psi)$ for each measurable subset $S\subseteq \Psi$. 
Likewise, given a truncated CRM $\Theta_K$, we define 
its normalization $\Xi_K$ via $\Xi_K(S) \defined \Theta_K(S)/\Theta_K(\Psi)$. 
Note that any simulation algorithm
for $\Theta_K$ can be used for $\Xi_K$ by simply normalizing the result. 
This does not depend on the particular representation 
of the CRM, and thus applies equally to all the representations in \cref{sec:seqcrms}.

The first step in the analysis of NCRM truncations is to define their approximation error
in a manner similar to that of CRM truncations.  
Since $\Xi$ and $\Xi_K$ are both normalized, they are distributions on $\Psi$; thus, 
observations $X_{1:N}$ are generated \iid from $\Xi$, and $Z_{1:N}$ are generated $\iid$ from $\Xi_K$. 
$Y_{1:N}$ and $W_{1:N}$ have the same definition as for CRMs.
As in the developments of \cref{sec:trunc}, the theoretical
results of this section rely on a general upper bound, provided by \cref{lem:protobound-ncrm}. %
\begin{nlem}[NCRM protobound]\label{lem:protobound-ncrm}
Let $\Theta\dist\distCRM(\nu)$, and let its truncation be $\Theta_K$. 
Let their normalizations be $\Xi$ and $\Xi_K$ respectively.
If
\[
X_{n} \given \Xi &\distiid \Xi, & 
Z_{n} \given \Xi_{K} &\distiid \Xi_K, \\
Y_{n} \given X_{n} & \distind f(\cdot \given X_{n}), &
W_{n} \given Z_{n} &\distind f(\cdot \given Z_{n}), 
\]
then
\[
\frac{1}{2}\|p_{N, \infty} - p_{N, K}\|_{1}
&\leq 1- \Pr\left(X_{1:N} \subseteq \supp(\Xi_K)\right),
\]
where $p_{N, \infty}$, $p_{N, K}$ are the marginal densities of $Y_{1:N}$ and
$W_{1:N}$, respectively. 
\end{nlem}

The analysis of CRMs in \cref{sec:trunc} relied heavily on the Poisson process stucture of the
rates in $\Theta$ and $X_{1:N}$; unfortunately, the rates in $\Xi$ do not possess
the same structure and thus lack many useful independence properties (the rates must sum to one). Likewise,
sampling $X_n$ for each $n$ does not depend on the atoms of $\Xi$ independently
($X_n$ randomly selects a single atom based on their rates).
Rather than using the basic definitions of the above random quantities to derive an error
bound, we decouple the atoms of $\Xi$ and $X_{1:N}$ using a technique from
extreme value theory. 
A \emph{Gumbel} random variable $T$ with location $\mu\in \reals$ and 
scale $\sigma > 0 $, denoted $T\dist\distGumbel(\mu, \sigma)$,
is defined by the cumulative distribution function and corresponding density
\[
\Pr(T \leq t) &= e^{-e^{-\frac{t-\mu}{\sigma}}} &&\text{and}&  \frac{1}{\sigma}e^{-\left(\frac{t-\mu}{\sigma}\right)-e^{-\left(\frac{t-\mu}{\sigma}\right)}}.
\]
An interesting property of the Gumbel distribution is that if one perturbs
the log-probabilities of a finite discrete distribution by \iid $\distGumbel(0, 1)$
random variables, the $\argmax$ of the resulting set is a sample from
the discrete distribution~\citep{Gumbel:1954,Maddison:2014ww}. This
technique is invariant to normalization, as the $\argmax$ is invariant
to the corresponding constant shift in the log-transformed space.
For present purposes, we develop the infinite extension of this result:

\begin{nlem}[Infinite Gumbel-max sampling] \label{lem:gumbelmax}
Let $(p_i)_{i=1}^\infty$ be a collection of positive numbers
such that $\sum_i p_i < \infty$ and let $\bp_{j} \defined \frac{p_j}{\sum_i p_i}$. 
If $(T_i)_{i=1}^\infty$ are \iid $\distGumbel(0,1)$ random variables,
then $\argmax_{i\in\nats} \, \, T_i+\log p_i$ exists, is unique \as, and
has distribution
\[
\argmax_{i\in\nats} \, \, T_i+\log p_i \dist \distCat\left(\left(\bp_{j}\right)_{j=1}^\infty \right).
\]
\end{nlem}

The proof of this result, along with the others in this section, 
may be found in \cref{app:normproofs}. The utility of \cref{lem:gumbelmax} is that it allows the construction of $\Xi$ and $X_{1:N}$
without the problematic coupling of the underlying CRM atoms due to normalization;
rather than dealing directly with $\Xi$, we 
log-transform the rates of $\Theta$, perturb
them by \iid $\distGumbel(0, 1)$ random variables, and characterize the distribution of the maximum rate
in this process. The combination of this distribution with \cref{lem:gumbelmax} yields 
the key proof technique used to develop the truncation bounds in \cref{thm:ndbound,thm:nfvbound}.
The results presented in this section are summarized in \cref{tbl:results} in \cref{sec:simulation}.
\subsection{Series representations}
The following result provides a general truncation error bound for normalized series representations,
specifies its range, and guarantees that it decays to 0 as $K\to\infty$. 
We again use the general series representation notation from \cref{eq:general-tau-series-rep},
where $g$ is a distribution on $\reals_+$, and $\tau : \reals_+\times\reals_+\to\reals_+$ is a measurable function
such that $\lim_{u\to\infty} \tau(v, u) = 0$ for $g$-almost every $v$.
\begin{nthm}[Normalized series representation truncation error bound]\label{thm:ndbound}
The error of approximating a series representation of $\Xi\dist\distNCRM(\nu)$ with its
truncation $\Xi_K$ satisfies
\[
0&\leq \frac{1}{2}\|p_{N,\infty}-p_{N,K}\|_1 \leq  1- \left(1- B_K\right)^N\leq 1,
\]
where
\[
B_K &\defined \EE\bigg[\int_{0}^{\infty}\!\!\! J\left(\Gamma_K, t\right)\left(\int_0^1 J\left(\Gamma_K u, t\right)\,\dee u\right)^{K-1} 
\left(-\frac{\dee}{\dee t}e^{\int_0^\infty\left( J\left(u+\Gamma_K, t\right) - 1 \right)\dee u}\right)\dee t\bigg], \label{eq:seriesBK} 
\]
\[
J(u, t) &= \EE\left[e^{-t\cdot \tau(V, u)}\right], &
V &\dist g, &
& \text{and} &
\Gamma_K &\dist \distGam(K, 1).
\]
Furthermore, $\lim_{K\to\infty}B_{K} = 0$.
\end{nthm}
\begin{exa}[Dirichlet process, $\distDP(\gamma)$, $\BDrep$]\label{eg:drepdpbound}
The Dirichlet process with concentration $\gamma>0$ is a normalized gamma
process $\distNGammaP(\gamma, 1, 0)$. %
From \cref{ex:b-rep-gamma-process} we have $c_\nu = \gamma$ and $g_\nu = \distExp(1)$, and from \cref{sec:seriestrunc}
we have $\tau(v, u) = v e^{-u/c_\nu}$. Therefore $J$ and its antiderivative are
\[
J(u, t) &= \EE\left[e^{-t V e^{-u/\gamma}}\right] = \left(1+te^{-u/\gamma}\right)^{-1} 
\quad \text{and} \quad 
\int J(u, t)\dee u = \gamma \log\left(e^{u/\gamma}+t\right).
\]
Using the antiderivative to evaluate the integrals in the formula for $B_K$, writing the expectation over $\Gamma_K\dist\distGam(K, 1)$ explicitly,
and making a change of variables we have 
\[
B_K &= \frac{\gamma^{K+1}}{\Gamma(K)} \int_0^\infty  \int_1^\infty
\left(\log\left(\frac{s+t}{1+t}\right)\right)^{K-1}
\left(s+t\right)^{-(\gamma+2)}
 \,\dee s\,\dee t = \left(\frac{\gamma}{1+\gamma}\right)^K,
\]
where the last equality is found by multiplying and dividing the integrand by $(1+t)^{-(\gamma+2)}$, and making the change of variables from $s$ to $x = \log\frac{s+t}{1+t}$.
Therefore, the truncation error can be bounded by
\[
\frac{1}{2}\left\|p_{N, \infty}-p_{N, K}\right\|_1 
&\leq 1-\left(1-\left(\frac{\gamma}{\gamma+1}\right)^K\right)^N
\sim N \left(\frac{\gamma}{\gamma+1}\right)^{K} \qquad K\to\infty.
\]
\end{exa}
The bound in \cref{eg:dpbound} has exponential decay, and reproduces earlier
DP truncation error bound rates due to \citet{Ishwaran:2001} and \citet{Ishwaran:2002}.
However, the techniques used in past work do not generalize beyond the
Dirichlet process, while those developed here apply to any NCRM. 

\subsection{Superposition representations}
The following result provides a general truncation error bound for normalized superposition representations,
specifies its range, and guarantees that it decays to 0 as $K\to\infty$. We once again
rely on the property that the truncation $\Theta_K$ and tail $\Theta_K^+$
are mutually independent CRMs, as expressed in \cref{eq:superpositionsplitproperty}, with the tail measure denoted $\nu_K^+$.
\begin{nthm}[Truncation error bound for normalized superposition representations]\label{thm:nfvbound}
The error of approximating a superposition representation of $\Xi\dist\distNCRM(\nu)$ with its
truncation $\Xi_K$ satisfies
\[
0 &\leq \frac{1}{2}\|p_{N,\infty}-p_{N,K}\|_1 \leq  1- \left(1- B_K\right)^N\leq 1,
\shortintertext{where}
B_K &\defined \int_{0}^\infty \left(\int\theta e^{-\theta t}\nu_K^+(\dee\theta)\right)e^{\int \left(e^{-\theta t}-1\right)\nu(\dee\theta)}\,\dee t. \label{eq:superBK}
\] 
Furthermore, $\lim_{K\to\infty}B_{K} = 0$.
\end{nthm}
This bound can be applied by using the tail measures derived earlier in \cref{sec:superpositiontrunc}.

\begin{exa}[Dirichlet process, $\distDP(\gamma)$, $\DBFrep$]\label{eg:dpbound}
As in \cref{eg:drepdpbound}, we view the Dirichlet process with concentration $\gamma>0$ as a normalized gamma
process $\distNGammaP(\gamma, 1, 0)$. %
First, by \cref{lem:em1nint}, the integral in the exponential is
\[
\exp\left(\int(e^{-t\theta}-1)\nu(\dee\theta)\right)
&= \exp\left(\gamma \int_0^\infty (e^{-t\theta}-1)\theta^{-1}e^{-\theta}\dee\theta \right)
= (t+1)^{-\gamma}.
\]
\cref{ex:b-rep-gamma-process} shows
$c_{\nu}=\gamma$ and  
$g_{\nu}(v) = e^{-v}$, 
  and \cref{eq:rate-measure-for-theta-r} provides the tail measure $\nu_{K}^+$ for the decoupled Bondesson representation,
\[
\nu_K^+(\dee\theta) &= \frac{\gamma}{\xi}\sum_{k=K+1}^\infty \frac{\xi^k}{\Gamma(k)}\left(\int_0^1(-\log x)^{k-1}x^{\xi-2} e^{-\theta x^{-1}}\dee x\right) \dee\theta.
\]
Substituting this result, using Fubini's theorem to swap the order of integration and summation, evaluating the integral over $\theta$, and making the substitution
$x=e^{-s}$ yields
\[
B_K &=\frac{\gamma}{\xi}\sum_{k=K+1}^\infty \frac{\xi^k}{\Gamma(k)}\iint_{s, t\geq 0}\hspace{-.3cm}\frac{s^{k-1}e^{-(\xi-1)s}(t+1)^{-\gamma}}{(e^s+t)^2}\dee s\,\dee t.
\]
Noting that $\forall s \geq 0, \, e^{s} \geq 1$, we have for any $a\in(0, 1]\cap (0, \gamma)$,
\[
 \hspace{-.4cm}B_K &\leq \frac{\gamma}{\xi}\sum_{k=K+1}^\infty \frac{\xi^k}{\Gamma(k)}\int_0^\infty s^{k-1}e^{-(\xi+a)s}\dee s\int_0^\infty(t+1)^{-(\gamma+1-a)}\dee t\\
  &= \frac{\gamma}{(\gamma-a)\xi}\sum_{k=K+1}^\infty \left(\frac{\xi}{\xi+a}\right)^k
  = \frac{\gamma}{a(\gamma-a)}\left(\frac{\xi}{\xi+a}\right)^K.
\]
Therefore, for any $a\in(0, 1]\cap(0, \gamma)$,
\[
\frac{1}{2}\left\|p_{N, \infty}-p_{N, K}\right\|_1 
&\leq 1-\left(1-\frac{\gamma}{a(\gamma-a)}\left(\frac{\xi}{\xi+a}\right)^K\right)^N 
\!\!\!\sim \frac{N \gamma}{a(\gamma-a)} \left(\frac{\xi}{\xi+a}\right)^K \quad K\to\infty.\label{eq:optimizeddpbound}
\]
To find the tightest bound, one can minimize with respect to $a$ given $\gamma, \xi, K$.
\end{exa}

\begin{exa}[Normalized gamma process, $\distNGammaP(\gamma, \lambda, d)$, $\Vrep$]\label{eg:gendpvrep}
By \cref{lem:em1nint}, the integral in the exponential is
\[
\exp\left(\int (e^{-\theta t}-1)\nu(\dee\theta)\right)
&= \left\{\begin{array}{ll}
\exp\left(-\gamma\lambda^{1-d} d^{-1}((t+\lambda)^d - \lambda^d)\right) & d>0\\
\left(\frac{\lambda}{t+\lambda}\right)^{\gamma\lambda} & d=0,
\end{array}\right. \label{eq:gendpvrep-part1}
\]
and the standard gamma integral yields
\[
\int\theta e^{-\theta t}\nu_K^+(\dee\theta) 
&= \gamma\frac{\lambda^{1-d}}{\Gamma(1-d)}\int \theta^{-d}e^{-(K+t+\lambda)\theta}\,\dee\theta 
= \gamma\lambda^{1-d}(K+t+\lambda)^{d-1}. \label{eq:gendpvrep-part2}
\]

When $d > 0$, multiplying the previous two displays and integrating over $t\geq 0$ yields
\[
B_K
&= \gamma \lambda^{1-d}e^{\gamma \lambda/d} \int_\lambda^\infty (K+t)^{d-1}e^{-\gamma\lambda^{1-d} t^d/d}\,\dee t 
\leq C_{\gamma,\lambda,d}\left(K+\lambda\right)^{d-1},
\]
where we have used $\left(K+t\right)^{d-1} \leq \left(K+\lambda\right)^{d-1}$ for $t\geq \lambda$
and the change of variables $u = \gamma\lambda^{1-d}d^{-1}t^d$ to find that 
$C_{\gamma, \lambda, d} = e^{\sigma}\sigma^{1-d^{-1}}\lambda^{1-d}\Gamma\left(d^{-1}, \sigma\right)$, where $\sigma=\gamma\lambda d^{-1}$ and
$\Gamma(a, x) \defined \int_x^\infty \theta^{a-1}e^{-\theta}\dee\theta$ is the upper incomplete gamma function.
Therefore,
\[
\frac{1}{2}\left\|p_{N, \infty}-p_{N, K}\right\|_1 \leq 1-\left(1-C_{\gamma, \lambda, d}(K+\lambda)^{d-1}\right)^N \sim NC_{\gamma,\lambda,d}K^{d-1}\quad  K\to\infty.
\]

When $d = 0$, multiplying \cref{eq:gendpvrep-part1,eq:gendpvrep-part2} and integrating over $t\geq 0$ yields
\[
\frac{B_K}{\gamma\lambda^{1+\gamma\lambda}} &= \int_\lambda^\infty (K+t)^{-1}t^{-\gamma\lambda}\dee t 
\leq \left\{\begin{array}{ll}
(K+\lambda)^{-1}\left(\frac{\frac{1}{\gamma\lambda}(K+\lambda)^{1-\gamma\lambda}-\lambda^{1-\gamma\lambda}}{1-\gamma\lambda}\right) & \gamma\lambda \neq 1\\
K^{-1}\log\left(\frac{K+\lambda}{\lambda}\right) & \gamma\lambda = 1 ,
\end{array}\right.
\]
where we obtain the bound for $\gamma\lambda \neq 1$ by splitting the integral into the intervals $[\lambda, K+\lambda]$ and $[K+\lambda, \infty)$ and bounding each section separately,
and we obtain the bound for $\gamma\lambda = 1$ via the transformation $u = t/(K+t)$.
Therefore, asymptotically
\[
\frac{1}{2}\left\|p_{N, \infty}-p_{N, K}\right\|_1 \lesssim N
\left\{\begin{array}{ll}
C_{\gamma,\lambda}K^{-\min(1, \gamma\lambda)} & \gamma\lambda \neq 1\\
\lambda K^{-1}\log K & \gamma\lambda = 1
\end{array}\right. \qquad K\to\infty,
\]
where $C_{\gamma,\lambda} \defined \max\left(\frac{\lambda^{\gamma\lambda}}{1-\gamma\lambda}, \frac{\gamma\lambda^2}{\gamma\lambda-1}\right)$.
\end{exa}
Truncation of the $\distNGammaP(\gamma, \lambda, d)$ has been studied previously: \citet{Argiento:2015}
threshold the weights of the unnormalized CRM to be beyond a fixed level $\epsilon > 0$ prior to normalization, and develop 
error bounds for that method of truncation. These results are not directly comparable to those of the present work due to the different
methods of truncation (i.e.~sequential representation termination versus weight thresholding).

\subsection{Hyperpriors}
As in the CRM case, we can place priors on the hyperparameters of the NCRM rate measure (i.e.~$\gamma$, $\alpha$, $\lambda$, $d$, etc.). 
We conclude our investigation of NCRM truncation error by
showing how bounds developed in this section can be modified
to account for hyperpriors. 
Note that we make the dependence of $B_K$ on the hyperparameters $\Phi$ 
explicit in the notation of \cref{prop:truncncrm-hyperprior}. 
\bnprop[NCRM truncation error with a hyperprior] \label{prop:truncncrm-hyperprior}
Given hyperparameters $\Phi$, consider a representation for $\Theta \given \Phi \dist \distCRM(\nu)$, let $\Xi \given \Phi$ be its normalization,
and let $B_K(\Phi)$ be given by \cref{eq:seriesBK} (for a series representation) or \cref{eq:superBK} (for a superposition representation).
The error of approximating $\Xi$ with its truncation $\Xi_K$ satisfies
\begin{align}
0\leq \frac{1}{2}\|p_{N,\infty}&-p_{N,K}\|_1 \leq  1- \left(1- \EE\left[B_K(\Phi)\right]\right)^N\leq 1.
\end{align}
\enprop

\begin{exa}[Dirichlet process, $\distDP(\gamma)$, $\BDrep$]
If we place a Lomax prior on $\gamma$, i.e.~$\gamma\dist\distLomP(a, 1)$,
then 
combining \cref{prop:truncncrm-hyperprior} and \cref{eg:drepdpbound} yields
\[
\frac{1}{2}\left\|p_{N, \infty}-p_{N, K}\right\|_1 &\leq  
1-\left(1-\frac{\Gamma(a+1)\Gamma(K+1)}{\Gamma(a+K+1)}\right)^N 
\sim N\Gamma(a+1)(K+1)^{-a}  \quad K\to\infty.
\]
\end{exa}

\section{Simulation and computational complexity}
\label{sec:sampling}

The sequential representations in \cref{sec:seqcrms} are each
generated from a different finite sequence of distributions, resulting in a
different expected computational cost for the same truncation level. 
Thus, the truncation level itself is not an appropriate parameter with which to compare
the error bounds for different representations and we require a characterization of
the computational cost. 
We investigate the mean complexity $\EE[R]$ of each representation,
where $R$ is the number of random variables sampled, as a 
function of the truncation level for each of the representations in \cref{sec:seqcrms}. 

We begin with the series representations. For each value of $k=1, \dots, K$, 
each series representation generates a single trait $\psi_{k} \dist G$ and a rate $\theta_k$
composed of some transformation of random variables. Thus, all of the series representations in this 
work satisfy $\EE[R] = rK$ for some constant $r$:  by inspection,
the inverse-L\'evy representation has $r=2$, and all the remaining series representations have $r=3$.

The superposition representations, on the other hand, generate a Poisson
random variable to determine the number of atoms at each value of $k=1, \dots, K$,
and then generate those atoms. Therefore, the mean simulation complexity
takes the form $\EE[R] = \sum_{k=1}^K1+r_k\EE[C_k]$ for some 
constants $r_k$ that might depend on the value of $k$.
For the decoupled Bondesson representation, $r_k = 3$ since each atom requires
generating three values ($\psi_{ki}$, $V_{ki}$, and $T_{ki}$), and $\EE[C_k] = c/\xi$, so
$\EE[R] = \left(\frac{3c}{\xi}+1\right)K$. 
For the size-biased representation, $r_k = 3$ since each atom requires
generating three values ($\psi_{ki}$, $z_{ki}$, and $\theta_{ki}$), and $\EE[C_k] = \eta_k$, so
$\EE[R] = K+3\sum_{k=1}^K\eta_k$. Note that here $\EE[R] \sim K,$ for $K\to\infty$ since $\eta_k$ is a decreasing sequence.
For the power-law representation,
$r_k = k+2$, since each atom requires generating $\psi_{ki}$, $V_{ki}$,
and $k$ beta random variables, and therefore $\EE[R] = \left(1+\frac{5\gamma}{2}\right)K+\frac{\gamma}{2}K^{2}$.

\section{Summary of results}\label{sec:simulation}
\cref{tbl:results} summarizes our truncation and simulation cost results as
applied to the beta, (normalized) gamma, and beta prime processes. 
Results for the Bondesson representation of $\distBP(\gamma, 1, 0)$
as well as the decoupled Bondesson representations of $\distBP(\gamma, \lambda, 0)$ and 
$\distGammaP(\gamma, \lambda, 0)$ were previously known, and are reproduced by our results. 
All other results in the table are novel to the best of the authors' knowledge. 
It is interesting to note that the bounds and expected costs within each of the 
representation classes often have the same form, aside from some constants. 
Across classes, however, they vary significantly, indicating that the chosen sequential
representation of a process has more of an influence on the truncation error than the process itself.

\begin{table}[t!]
\footnotesize
\caption{Asymptotic error bounds and simulation cost summary. 
Error bounds are presented up to a constant that varies between models.
Be = Bernoulli, OBe = odds Bernoulli, Poi = Poisson. %
}
\setlength{\tabcolsep}{3px}
\vspace{.5cm}
\begin{tabular}{l|l|l|c|c}
Rep. & Random Measure & $h$ & Asymptotic Error Bound & Complexity  \\ \hline \hline
IL &
$\distLomP(\gamma, \lambda^{-1})$ & Poi &  
$\vphantom{\prod\limits_{i=1}^\infty} N e^{-KW_0\left(\{3\gamma\lambda\max(\lambda, e)\}^{-1}\right)}$ &
$2K$ 
\\ \hline
\multirow{4}{*}{B} &
$\distBP(\gamma, \lambda \ge 1, 0)$ & Be &  
\multirow{3}{*}{$N\gamma\Big(\frac{\gamma\lambda}{\gamma\lambda + 1}\Big)^{K}$} &
\multirow{4}{*}{$3K$} 
\\ 
& $\distGammaP(\gamma, \lambda, 0)$ & Poi &
& \\
& $\distBPP(\gamma, \lambda, 0)$ & OBe & 
& \\ \cline{4-4}
& $\distDP(\gamma)$ & --- & $\vphantom{\prod\limits_{i=1}^\infty}N\big(\frac{\gamma}{\gamma+1}\big)^K$ & 
\\ \hline
T &--- &--- &  
See \cref{eq:thinningrepbd} &
$3K$  \\  \hline
\multirow{4}{*}{R} &
$\distBP(\gamma, \lambda, 0)$ & Be &  
\multirow{4}{*}{$ N 
\left\{\begin{array}{ll} 
   e^{-K W_0\left(\{3\gamma\lambda\}^{-1}\right)} & d = 0\\
   K^{-d(1-d)} & d > 0 \, (\distGammaP)\\
   K^{-1/d} & d > 0 \, (\distBP, \distBPP)
\end{array}\right.
$} &
\multirow{4}{*}{$3K$} 
\\ 
& $\distGammaP(\gamma, \lambda, 0)$ & Poi &
& \\
& $\distBPP(\gamma, \lambda, 0)$ & OBe & 
& 
\\ 
& & & & \\ \hline
\multirow{4}{*}{DB} &
$\distBP(\gamma, \lambda \ge 1, 0)$ & Be &  
\multirow{3}{*}{$ N \Big(\frac{\xi}{\xi + 1}\Big)^{K}$} &
\\ 
& $\distGammaP(\gamma, \lambda, 0)$ & Poi &
&\multirow{4}{*}{$\left(\frac{3c}{\xi}+1\right)K$}  \\ 
& $\distBPP(\gamma, \lambda>1, 0)$ & OBe & & \\ \cline{4-4}
& $\distDP(\gamma)$ & --- & $\vphantom{\prod\limits_{i=1}^\infty}\frac{N\gamma}{a(\gamma-a)}\big(\frac{\xi}{\xi+a}\big)^K, \, a\in(0, 1]\cap(0, \gamma)$ & 
\\ \hline
\multirow{4}{*}{SB} &
$\distBP(\gamma, \lambda, d)$ & Be &  
\multirow{3}{*}{$ N K^{d-1}$}&
\multirow{4}{*}{$K$} 
\\ 
& $\distGammaP(\gamma, \lambda, d)$ & Poi &
& \\
& $\distBPP(\gamma, \lambda, d)$ & OBe & 
& \\ \cline{4-4}
& $\distNGammaP(\gamma, \lambda, d)$ & --- 
& $\vphantom{\prod\limits_{i=1}^\infty}
 N\left\{\begin{array}{ll} 
K^{-1}\log K & d = 0, \gamma\lambda = 1\\
K^{-\min(1, \gamma\lambda)} & d = 0, \gamma\lambda \neq 1\\
K^{d-1} &  d > 0
\end{array}\right.$
\\ \hline
\multirow{3}{*}{PL} &
$\distBP(\gamma, \lambda, d)$ & Be &  
\multirow{4}{*}{$ N
\left\{\begin{array}{ll} 
\Big(\frac{\lambda}{\lambda+1}\Big)^K & d = 0 \, (\distBP, \distGammaP)\\ 
2^{-K} & d = 0 \, (\distBPP)\\ 
K^{1-1/d} & d > 0 \end{array}\right.$}&
\multirow{3}{*}{$\frac{\gamma}{2}K^2$} 
\\ 
& $\distGammaP(\gamma, \lambda, d)$ & Poi &
& \\
& $\distBPP(\gamma, \lambda>1, d)$ & OBe & 
& 
\\ 
&&&&\\\hline

\end{tabular}
\label{tbl:results}
\end{table}

\begin{figure}
  
\begin{subfigure}[b]{.48\textwidth}
    \raisebox{1.35cm}{\rotatebox{90}{$d=0.0$}} \includegraphics[width=.94\columnwidth, trim=0 30 0 0, clip]{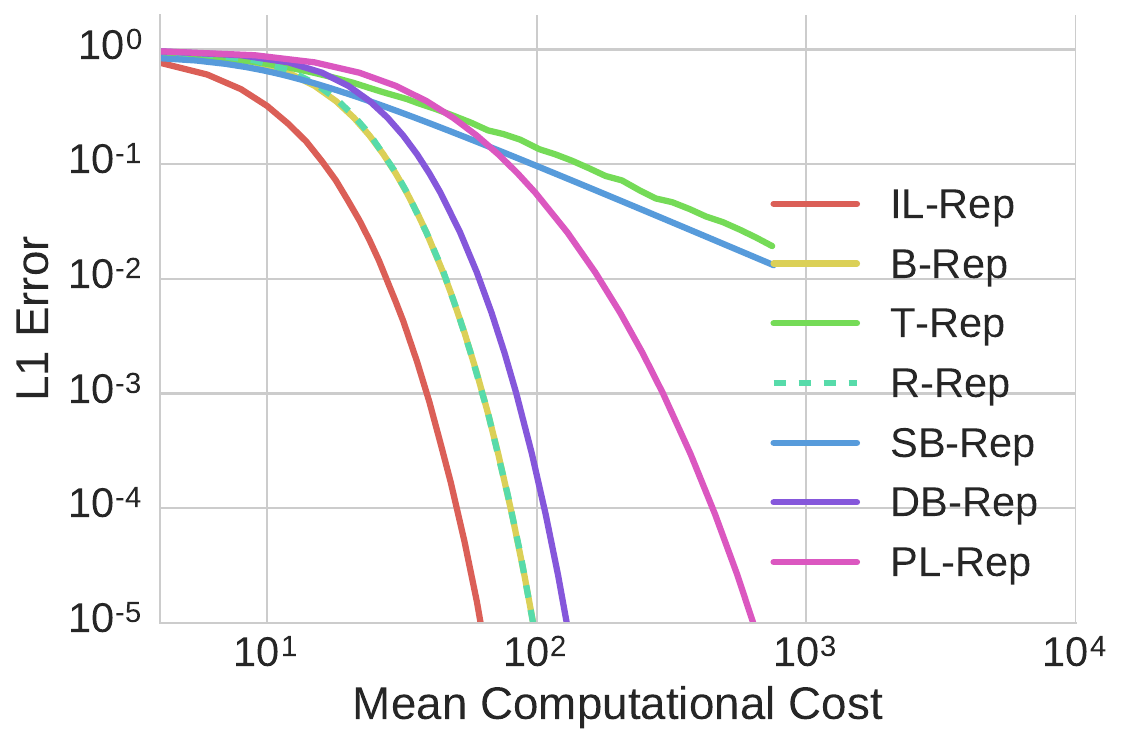}
  \end{subfigure}
 \begin{subfigure}[b]{.48\textwidth}
    \includegraphics[width=.88\columnwidth, trim=30 30 0 0, clip]{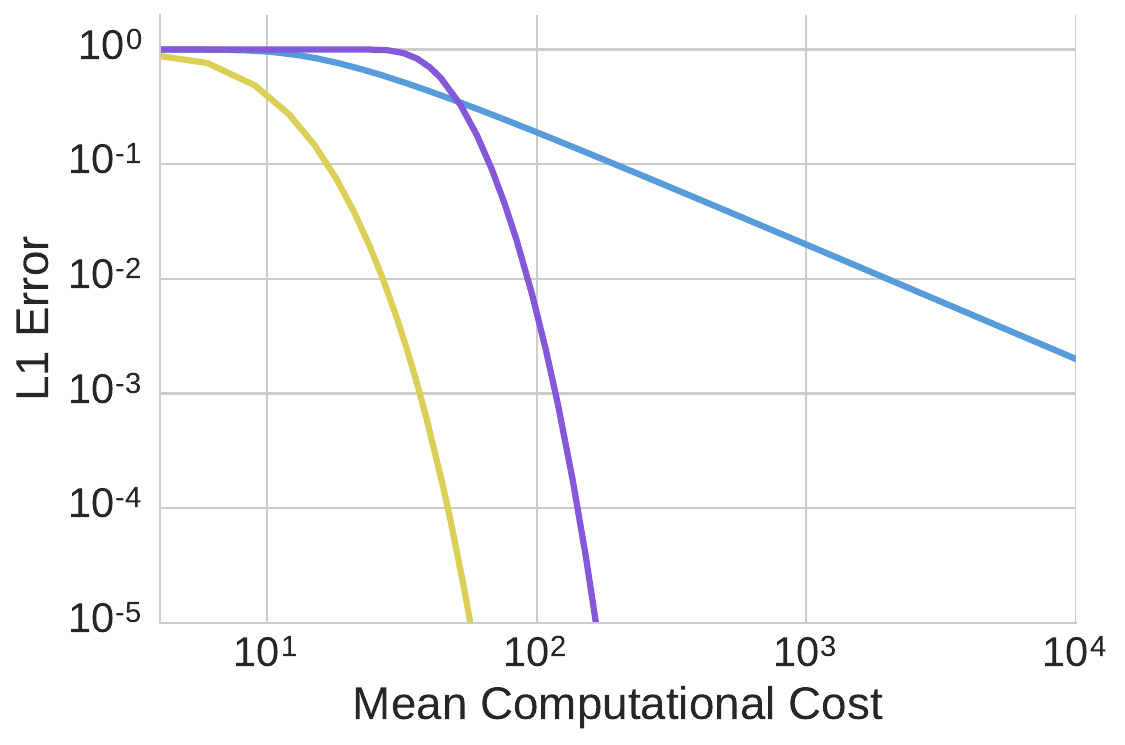}
  \end{subfigure}
\begin{subfigure}[b]{.48\textwidth}
    \raisebox{1.35cm}{\rotatebox{90}{$d=0.1$}} \includegraphics[width=.94\columnwidth, trim=0 30 0 0, clip]{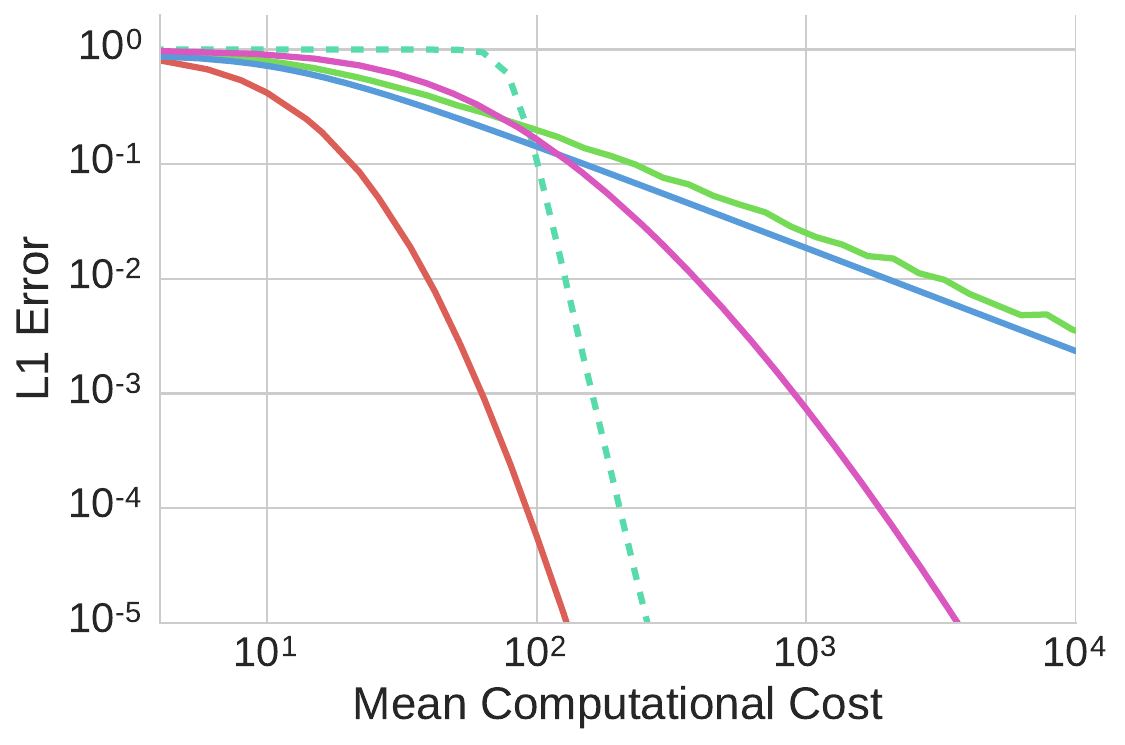}
  \end{subfigure}
 \begin{subfigure}[b]{.48\textwidth}
    \includegraphics[width=.88\columnwidth, trim=30 30 0 0, clip]{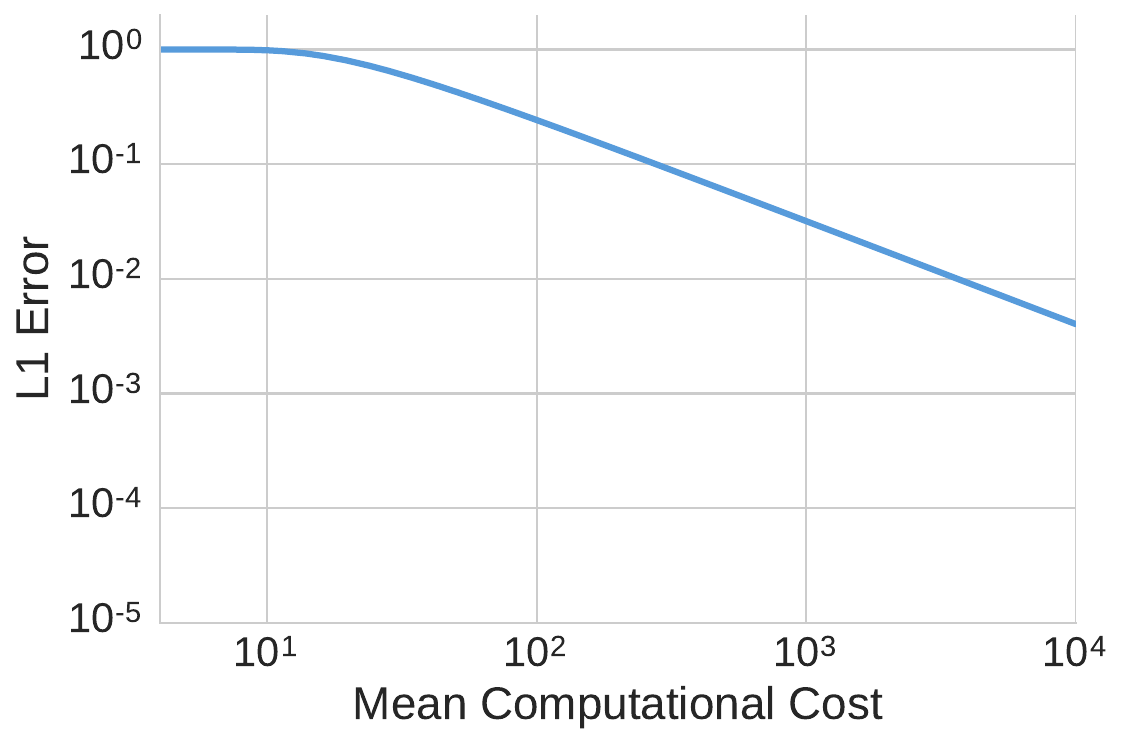}
  \end{subfigure}
\begin{subfigure}[b]{.48\textwidth}
    \raisebox{1.6cm}{\rotatebox{90}{$d=0.5$}} \includegraphics[width=.94\columnwidth, trim=0 0 0 0, clip]{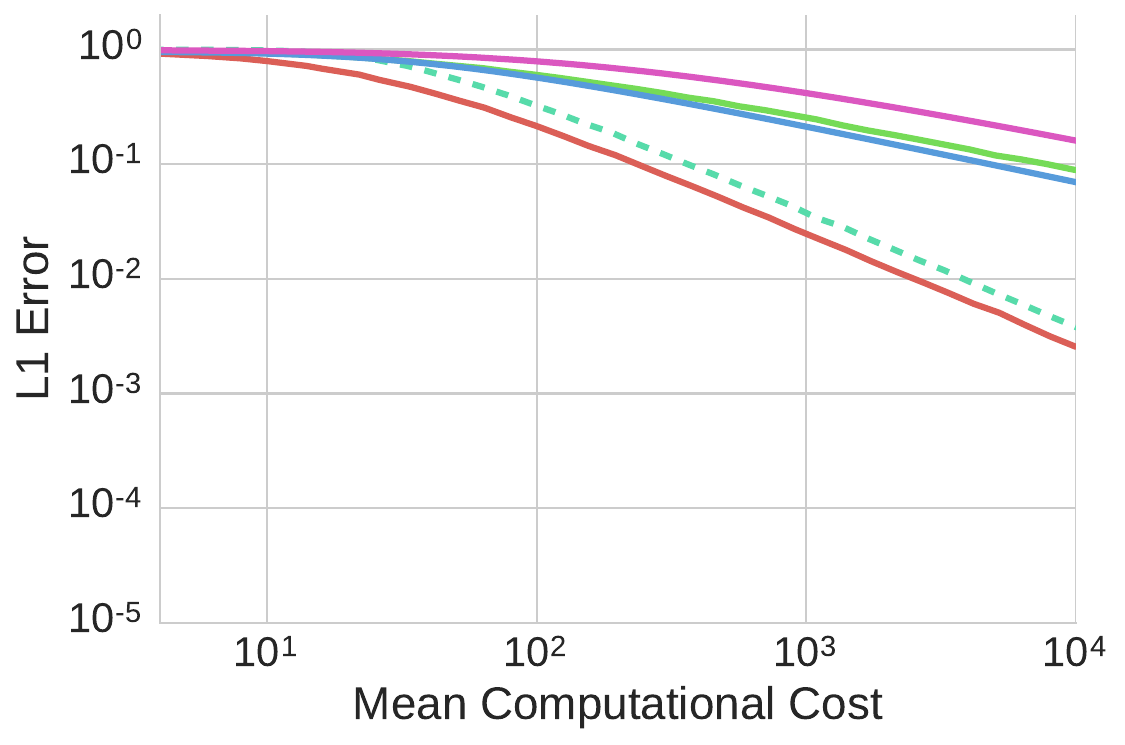}
    \caption*{$\distGammaP$}
  \end{subfigure}
 \begin{subfigure}[b]{.48\textwidth}
    \includegraphics[width=.88\columnwidth, trim=30 0 0 0, clip]{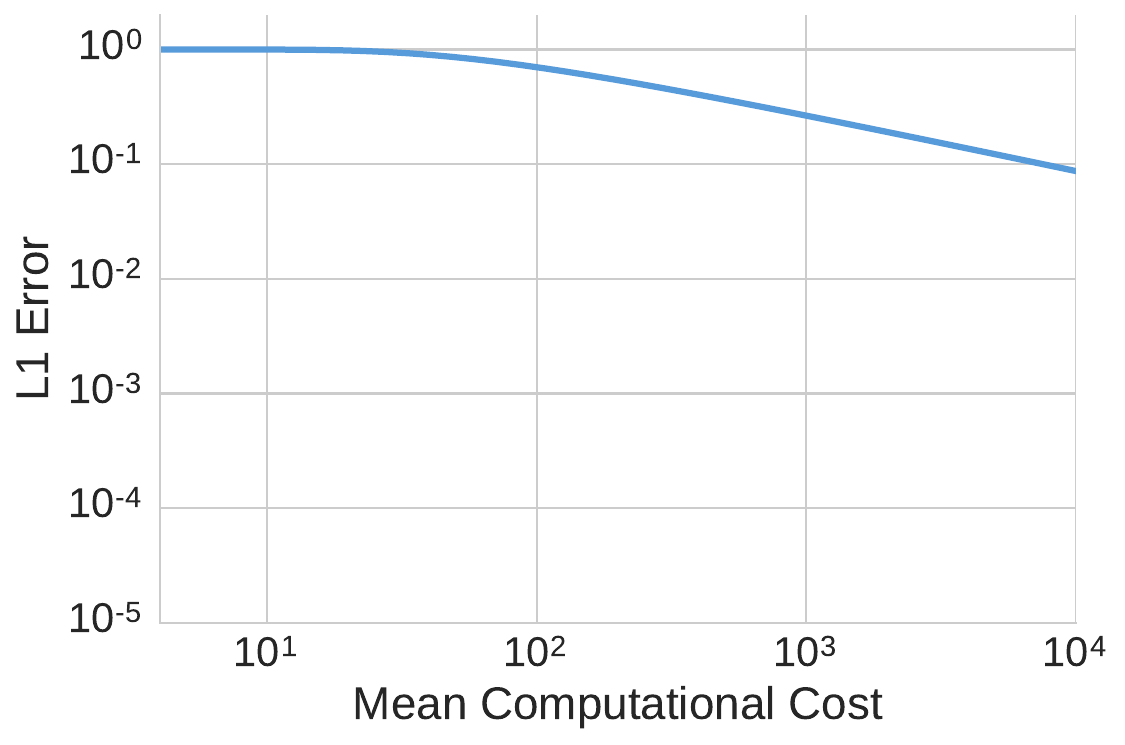}
    \caption*{$\distNGammaP$}
  \end{subfigure}

  \caption{Truncation error bounds for representations of the (normalized) gamma-Poisson process, with $\gamma = 1$, $\lambda = 2$,
  and $\xi = \gamma\lambda$. The left column 
is for the unnormalized process, while the right column is for the normalized process. Each row displays results for 
a different value of the discount parameter $d \in \{0, 0.1, 0.5\}$.}\label{fig:simstudy}
\end{figure}

\cref{fig:simstudy} shows a comparison of how the truncation error bounds vary with the expected computational cost $\EE[R]$ of simulation
for the (normalized) gamma process and Poisson likelihood with $N=5$ observations. Results shown for the thinning, rejection, and inverse-L\'evy representations
are computed by Monte-Carlo approximation of the formula for $B_{N,K}$ in \cref{eq:series-ezbound}, while all others use closed-form expressions from the examples
in \cref{sec:trunc,sec:ncrm-trunc}.
Note that the Bondesson and decoupled Bondesson representations
do not exist when $d > 0$. Further, only those representations for which we provide closed-form bounds in the examples are shown for the normalized gamma process;
we leave the numerical approximation of the results from \cref{thm:ndbound,thm:nfvbound} as an open problem. 
Similar figures for other processes (in particular, the beta-Bernoulli and beta prime-odds Bernoulli) are provided in \cref{app:examples}. 
Note that all bounds presented are improved
by a factor of two versus comparable past results in the literature, due to the reliance on
\cref{lem:protobound-crm,lem:protobound-ncrm} rather than the earlier bound found in~\citet{Ishwaran:2001}.

In \cref{fig:simstudy}, the top row shows results for the light-tailed process ($\gamma=1$, $\lambda=1$, $d=0$, and $\xi = c=\gamma\lambda$).
All representations except for thinning and size-biased capture its exponential truncation error decay.
This is due to the fact that the thinning representation generates increasingly many atoms of weight 0 as $K\to\infty$,
and the expected number of atoms at each outer index for the size-biased representation decays as $K\to\infty$. 
The inverse-L\'evy representation has the lowest truncation error as expected, as it is the only representation that generates a 
nonincreasing sequence of weights (and so must be the most efficient \citep{Arbel:2015}).  
Based on this figure and those in \cref{app:examples} for other processes, it appears that the Bondesson representation typically provides the best tradeoff between simplicity and efficiency,
and should be used whenever its conditions in \cref{prop:bondesson-representation} are satisfied. 
When the technical conditions are not satisfied, the rejection representation is a good alternative.
If ease of theoretical analysis is a concern,
the decoupled Bondesson representation provides comparable efficiency with the analytical simplicity of a superposition representation. 

The bottom two rows of \cref{fig:simstudy} show results for the heavy-tailed process ($\gamma=1$, $\lambda=2$, and $d\in\{0.1,0.5\}$).
The representation options are more limited, as the technical conditions of the Bondesson and decoupled Bondesson representations are not satisfied.
Here the rejection representation is often the best choice due to its simplicity and competitive performance with the inverse-L\'evy representation.
However, one must take care to check its efficiency beforehand using \cref{prop:r-rep-efficiency} given a particular choice of $\mu(\dee \theta)$. For example, 
the choice of $\mu(\dee \theta) \propto \theta^{-1-d}\dee\theta$ in the present work makes the rejection representation very inefficient when $d \ll 1$ for both
the gamma-Poisson (\cref{fig:simstudy}) and beta prime-odds Bernoulli (\cref{fig:bppfigs}) processes, but efficient for the beta-Bernoulli process (\cref{fig:bpfigs}).
If no $\mu(\dee\theta)$ yields reasonable results,
 the power-law representation is a good choice for $d\ll1$ as its truncation bound approaches the exponential decay of the light-tailed process.
For larger $d>0$ the size-biased representation is a good alternative.

Based on the results in \cref{fig:simstudy}, it appears that there is no single dominant representation
for all situations (provided the inverse-L\'evy representation is intractable, as it most often is).
However, as a guideline, the rejection and Bondesson representations tend to be good choices
for light-tailed processes, while the rejection, size-biased, and power-law representations are good choices for
heavy-tailed processes.

\section{Discussion}\label{sec:discussion}

We have investigated sequential representations, truncation error bounds, and simulation algorithms
for (normalized) completely random measures. In past work, the development and analysis of
these tools has occurred only on an ad hoc basis.
The results in the present paper, in contrast, provide a comprehensive characterization 
and analysis of the different types of sequential (N)CRM representations available to the practitioner. 
However, there are a number of remaining open questions and limitations.

First, this work does not consider the influence of observed data:
all analyses assume an \emph{a priori} perspective, as
truncation is typically performed before data are incorporated via posterior inference
(e.g.~in variational inference for the DP mixture~\citep{Blei:2006a} and 
BP latent feature model~\citep{DoshiVelez:2009}).
However, analysis of \emph{a posteriori} truncation has been studied in past
work as well~\citep{Ishwaran:2001,Gelfand:2002,Ishwaran:2002}.  
In the language of CRMs, observations introduce a fixed-location component in the posterior process,
while the unobserved traits are drawn from the (possibly normalized) ordinary component of a CRM~\citep{Ishwaran:2002,Broderick:2014}. 
We anticipate that this property makes observations
reasonably simple to include: the truncation tools provided
in the present paper can be used directly on the unobserved ordinary component, while the fixed-location
component may be treated exactly.

In addition, there are important open questions regarding the sequential representations developed in this work. 
It is unknown whether generalized versions of the Bondesson and decoupled Bondesson  representations can be developed for larger classes of rate measures.
The power-law representation does provide a partial answer in the decoupled Bondesson case.
Regarding size-biased representations, one might expect that the use of conjugate exponential family CRMs~\citep{Broderick:2014}
would yield a closed-form expression for the truncation 
bound. In all of the cases provided in this paper, this was indeed the case; the integrals were evaluated exactly and a closed-form 
expression was found. However, we were unable to identify a general expression applicable to all conjugate exponential family CRMs. 
Based on the examples provided, we conjecture that such an expression exists.
Finally, fundamental connections between some of the representations were left largely unexplored in this work.
This is an open area of research, although progress has been made by connecting decoupled Bondesson and size-biased representations
for (hierarchies of) generalized beta processes~\citep[Sec. 6.4]{Roy:2014}.

A final remark is that one of the primary uses of sequential representations in past work has been in the development of 
posterior inference procedures~\citep{Paisley:2010,Blei:2006a,DoshiVelez:2009}.
The present work provides no guidance on which truncated representations are best paired with which inference
methods. We leave this as an open direction for future research, which will require both theoretical 
and empirical investigation.

\section*{Acknowledgments}
The authors thank the anonymous reviewers for their thoughtful comments, which led to substantial
improvements in both the presentation and content of the paper.
The authors also thank Tin Nguyen for pointing out a flaw in the proof of a result appearing in a preprint draft of the paper.
All authors are supported by the Office of Naval Research under MURI grant N000141110688. 
J.~Huggins is supported by the U.S. Government under FA9550-11-C-0028 and awarded by the DoD, Air Force
Office of Scientific Research, National Defense Science and Engineering Graduate
(NDSEG) Fellowship, 32 CFR 168a. 
T.~Campbell and T.~Broderick are supported by DARPA award FA8750-17-2-0019.

\appendix

\section{Additional Examples} 
\label{app:examples}

In this section, we provide example applications of our theory to the beta process
and to the beta prime process.

\subsection{Beta process}\label{app:bpegs}

The beta process $\distBP(\gamma, \alpha, d)$~\citep{Teh:2009,Broderick:2012}
with discount parameter $d \in [0,1)$, concentration parameter $\alpha > -d$, and mass
parameter $\gamma > 0$, is a CRM with rate measure\phantom{\ref{eq:bpmsr}}
\begin{align}
\nu(\dee\theta) = \gamma\frac{\Gamma(\alpha + 1)}{\Gamma(1 - d)\Gamma(\alpha +
d)}\ind\left[\theta \leq 1\right] \theta^{-1-d}(1-\theta)^{\alpha+d-1}\dee\theta. \label{eq:bpmsr}
\end{align}
Setting $d = 0$ yields the standard beta process~\citep{Hjort:1990,Thibaux:2007}.
The beta process is often paired with a Bernoulli likelihood or negative binomial likelihood with $s\in\nats$ failures:
\begin{align}
&  \text{Bern: } &  h(x \given \theta) &= \ind\left[x\leq 1\right]\theta^x(1-\theta)^{1-x}, \\
&  \text{NegBinom: } & h(x\given \theta) &= {x + s - 1 \choose x}(1 - \theta)^{s}\theta^{x}. 
\end{align}
Note that for the Bernoulli likelihood $\pi(\theta) = 1 - \theta$ and for
the negative binomial likelihood $\pi(\theta) = (1-\theta)^{s}$.

\begin{figure}[t!]
\begin{subfigure}[b]{.32\textwidth}
    \includegraphics[width=\columnwidth]{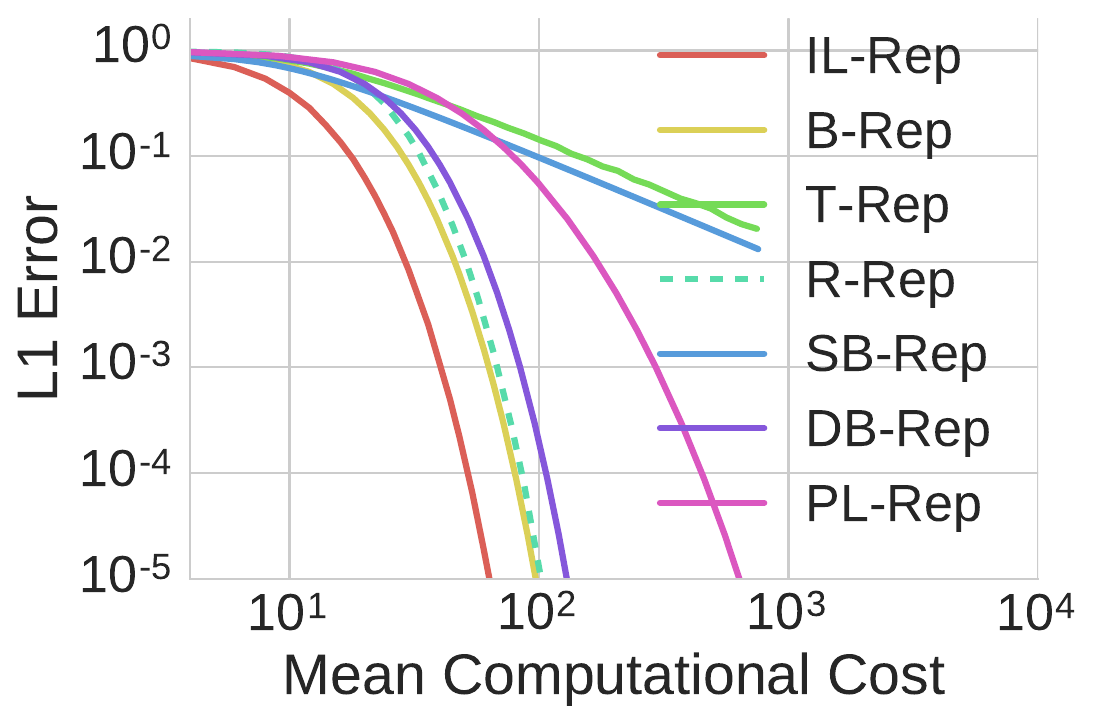}
    \caption*{$d=0.0$}
  \end{subfigure}
 \begin{subfigure}[b]{.32\textwidth}
    \includegraphics[width=.95\columnwidth, trim=33 0 0 0, clip]{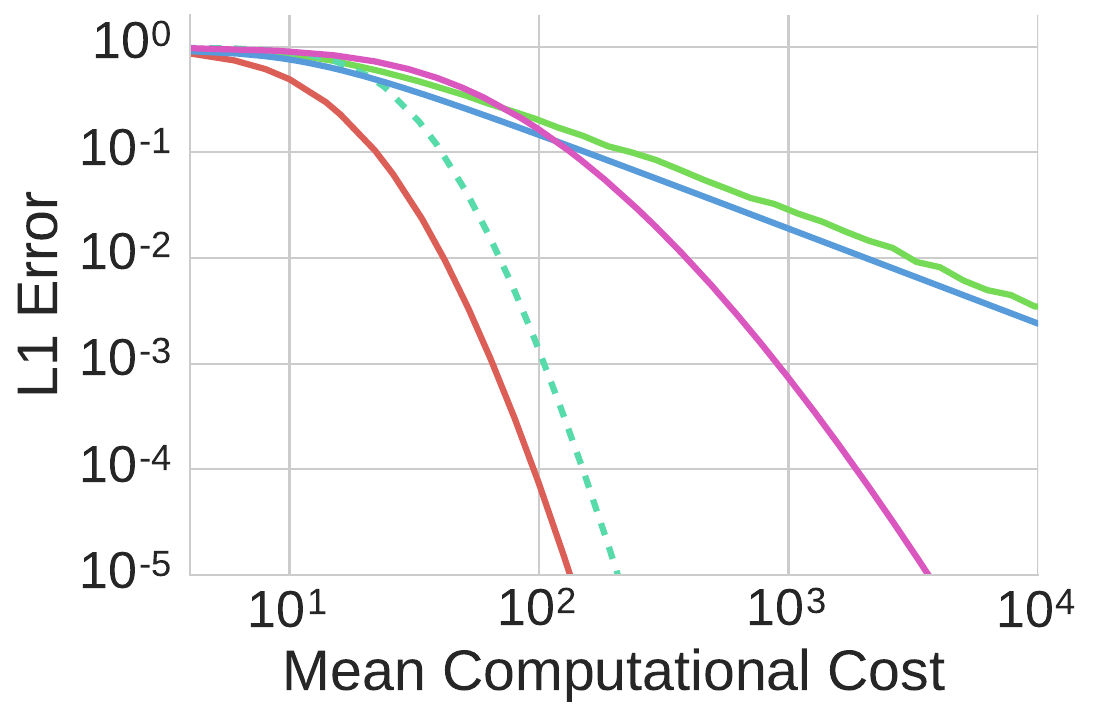}
    \caption*{$d=0.1$}
  \end{subfigure}
\begin{subfigure}[b]{.32\textwidth}
    \vspace{.7cm}\includegraphics[width=.95\columnwidth, trim=33 0 0 0, clip]{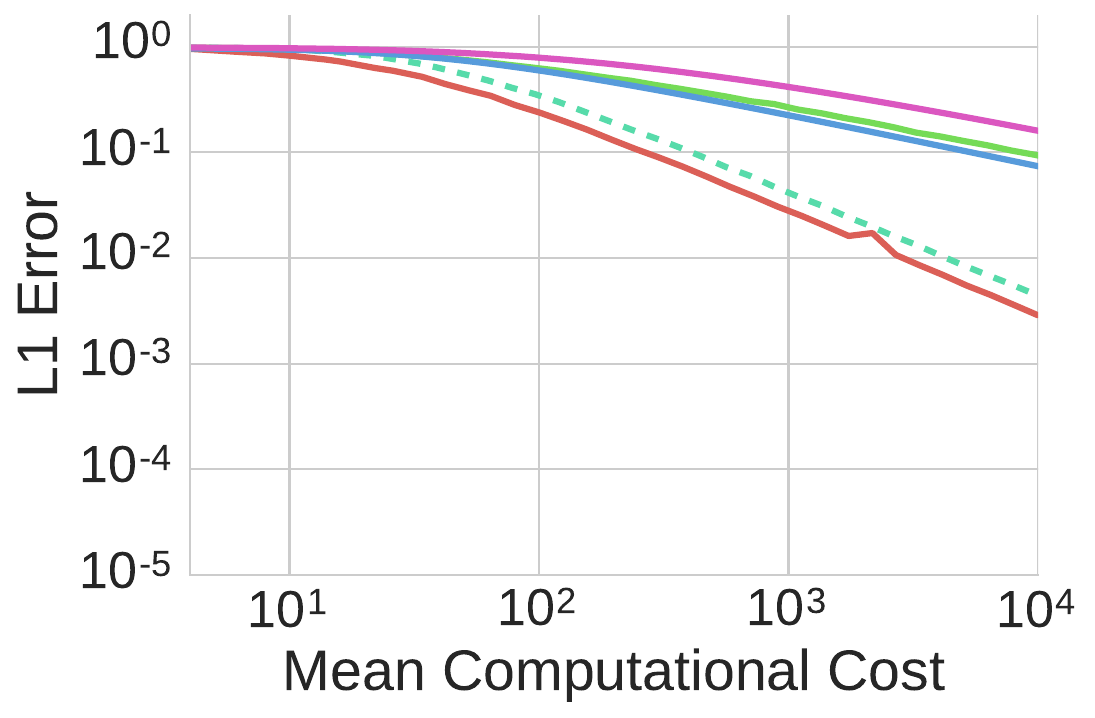}
    \caption*{$d=0.5$}
  \end{subfigure}
  \caption{Truncation error bounds for the beta-Bernoulli process. }\label{fig:bpfigs}
\end{figure}

\paragraph{Bondesson representation}
If $\alpha > 1$ and $d = 0$, then $\theta \nu(\theta) = \gamma\alpha(1 - \theta)^{\alpha-1}\ind[\theta \le 1]$
is non-increasing, $c_{\nu} = \lim_{\theta\to 0}\theta \nu(\theta) = \gamma\alpha$, and
$g_{\nu}(v) = (\alpha-1)(1 - v)^{\alpha-2} = \distBeta(v; 1, \alpha - 1)$. 
Thus, it follows from \cref{prop:bondesson-representation} that 
if $\Theta \isBDrep(\gamma\alpha, \distBeta(1, \alpha-1))$, then $\Theta \dist \distBP(\gamma, \alpha, 0)$. 
In the case of $\alpha = 1$, $g_{\nu}(v) = \delta_{1}$, so $V_{k} \equiv 1$ and the Bondesson representation
is equivalent to the inverse-L\'evy representation.
Since $\exp(-E_{k}/c) \dist \distBeta(1, c)$, the representation used 
in \citet{Teh:2007} is equivalent to the Bondesson representation for $\distBP(\gamma, 1, 0)$.

To obtain a truncation bound in the Bernoulli likelihood case, we 
can argue as in \cref{ex:d-rep-gpp-trunc}:
\[
B_{N,K}
&\le N \int_{0}^{\infty} \left(1 - \EE\left[\pi(ve^{-G_{K}/(\gamma\alpha)})\right]\right)\nu(\dee v) \\
&= N\gamma\alpha\,\EE[e^{-G_{K}/(\gamma\alpha)}]\int_{0}^{1}(1-v)^{\alpha - 1} \dee v \label{eq:b-rep-bbp-trunc} \\
&= N\gamma \left(\frac{\gamma\alpha}{1 + \gamma\alpha}\right)^{K}. 
\]
This result generalizes that in \citet{DoshiVelez:2009}, which applies only when $\alpha = 1$. 
\cref{prop:bondesson-representation} does not apply directly when $\alpha < 1$, since 
$\lim_{\theta \to 0} \theta \nu(\theta) = \infty$.
However, a representation can be obtained by using a trick from \citet{Paisley:2012}.
For $\alpha > 0$, let 
\[
\Theta' = \sum_{k=1}^{C} \theta'_{k} \delta_{\psi'_{k}},
\]
where $C \dist \distPoiss(\gamma)$, $\theta'_{k} \distiid \distBeta(1, \alpha)$,
and $\psi'_{k} \distiid G$.
Thus, $\Theta'$ is a CRM with rate measure 
$\gamma\alpha(1 - \theta)^{\alpha-1}\ind[\theta \le 1]\dee \theta$. 
If $\Theta \dist \distBP(\gamma\alpha/(\alpha+1), \alpha + 1, 0)$, which
can be generated according to \cref{prop:bondesson-representation}, then 
$\Theta'' = \Theta + \Theta'$ is a CRM with rate density on $[0,1]$ given by
\[
\gamma\alpha\theta^{-1}(1 - \theta)^{\alpha} + \gamma\alpha(1 - \theta)^{\alpha-1}
&= \gamma\alpha\theta^{-1}(1 - \theta)^{\alpha-1},
\]
hence $\Theta'' \dist \distBP(\gamma, \alpha, 0)$. 

\paragraph{Thinning representation}
If we let $g = \distBeta(1 - d, \alpha + d)$, then the thinning representation for $\distBP(\gamma, \lambda, d)$ is
\[
\Theta &= \sum_{k=1}^{\infty}V_{k}\ind(V_{k}\Gamma_{k} \le \gamma)\delta_{\psi_{k}}, &
& \text{with} &
V_{k} \distiid \distBeta(1 - d, \alpha + d). 
\]

\paragraph{Rejection representation}
To obtain a rejection representation for any $d$ when $\alpha \ge 1 - d$, let $\mu$ be 
the rate measure for $\distBP( \gamma \frac{\Gamma(\alpha + 1)}{\Gamma(2 - d)\Gamma(\alpha + d)}, 1-d, d)$.
We then have that
\[
\mu[x,\infty) &= 
\begin{cases}
\gamma' d^{-1}(x^{-d} - 1) & d > 0 \\
-\gamma'\log x & d = 0
\end{cases} &
& \text{and} &
\muinv(u) &= 
\begin{cases}
(1 + du/\gamma')^{-1/d} & d > 0 \\
e^{-u/\gamma'} & d = 0,
\end{cases}
\]
where $\gamma' \defined  \gamma \frac{\Gamma(\alpha + 1)}{\Gamma(1 - d)\Gamma(\alpha + d)}$. 
Thus, we can apply the inverse-L\'evy method analytically for $\mu$. 
Since we have constructed $\mu$ such that $\dee \nu/\dee \mu \le 1$, 
we can use $\mu$ to construct the rejection representation
\[
\Theta = \sum_{k=1}^{\infty} V_{k} \ind(U_{k} \le (1 - V_{k})^{\alpha + d - 1}) \dist \distBP(\gamma, \alpha, d), \qquad \alpha \ge 1 - d, \\
\text{with} \qquad V_{k} = \begin{cases}
(1 + d \Gamma_{k}/\gamma')^{-1/d} & d > 0 \\
e^{-\Gamma_{k}/\gamma'} & d = 0
\end{cases},
\qquad 
U_{k} \distiid \distUnif[0,1].
\]
The expected number of rejections is
\[
-d^{-1}\left[1 + 2 d\,\leftidx{_2^{}}{F}{_1^{(0,0,1,0)}}(-\alpha-d,-d,-d; 1) + d \,\leftidx{_2^{}}{F}{_1^{(0,1,0,0)}}(-\alpha-d,-d,-d; 1)\right],
\]
where ${_2^{}}{F}{_1}$ is the ordinary hypergeometric function and the parenthetical superscripts
indicate partial derivatives.
This quantity monotonically diverges to $\infty$ as $d \to 1$. 

To obtain a truncation bound in the Bernoulli likelihood case, we consider the $d > 0$ and $d = 0$
settings separately.
If $d > 0$, we have 
\[
\frac{B_{N,K}}{N\gamma'}
&\le  \int_{0}^{1}F_{K}(\gamma' d^{-1}(x^{-d} - 1)) x^{-d}(1 - x)^{\alpha + d - 1}\dee x \\
&\le  \int_{0}^{a} x^{-d}(1 - x)^{a + d - 1}\dee x +  F_{K}(\gamma' d^{-1}(a^{-d} - 1)) \int_{a}^{1}x^{-d}(1 - x)^{\alpha + d - 1}\dee x \\
&\le  \int_{0}^{a} x^{-d}\dee x + a^{-d} F_{K}(\gamma' d^{-1}(a^{-d} - 1)) \int_{a}^{1}(1 - x)^{\alpha + d - 1}\dee x \\
&\le (1-d)^{-1}a^{1-d} + a^{-d} \left(3\gamma' d^{-1}(a^{-d} - 1)/K\right)^{K} \\
&\le  (1-d)^{-1}a^{1-d} +  \left(3\gamma' d^{-1}/K\right)^{K}a^{-(K+1)d}.
\]
Setting the two terms equal and solving for $a$ we obtain $a^{1+d(K-1)} = \left(3\gamma' (d^{-1}-1)/K\right)^{K}$ and conclude that
\[
B_{N,K} 
&\le 2N\gamma' \left(\frac{3\gamma' (d^{-1}-1)}{K}\right)^{\frac{K}{1+d(K-1)}}  
\sim 2N\gamma' \left(\frac{3\gamma' (d^{-1}-1)}{K}\right)^{1/d} \qquad K \to \infty.
\]
If $d = 0$, we have 
\[
\frac{B_{N,K}}{N\gamma'} 
&\le \int_{0}^{1}F_{K}(-\gamma' \log x) (1 - x)^{\alpha - 1}\dee x \\
&\le  \int_{0}^{a} (1 - x)^{\alpha - 1}\dee x + F_{K}(-\gamma' \log a) \int_{a}^{1}(1 - x)^{\alpha - 1}\dee x \\
&\le a + \left(-3\gamma' \log a /K\right)^{K}.
\]
Setting the two terms equal and solving for $a$ we conclude that
\[
B_{N,K}
&\le 2N\gamma\alpha e^{-K W_{0}(\{3\gamma\alpha\}^{-1})}, 
\]
where $W_0$ is as defined in \cref{eq:lambertdefn}.

\paragraph{Decoupled Bondesson and power-law representations}
The decoupled Bondesson representation for $\distBP(\gamma, \alpha, 0)$ from \citet{Paisley:2010} 
was extended by \citet{Broderick:2012} to the $\distBP(\gamma, \alpha, d)$ setting. 
The \citet{Broderick:2012} construction for the $\distBP(\gamma, \alpha, d)$ is in fact the ``trivial''
power-law representation $\PLFrep(\gamma, \alpha, d, \delta_{1})$ (the decoupled Bondesson representation
is the special case when $d=0$). 

In the Bernoulli likelihood case, the truncation bound
for the decoupled Bondesson representation is
\[
B_{N,K}
&\le \frac{N\gamma\alpha}{\xi}\sum_{k=K+1}^{\infty}\EE[V e^{-T_{k}}] 
= \frac{N\gamma\alpha}{\xi}\sum_{k=K+1}^{\infty}\frac{1}{\alpha}\left(\frac{\xi}{1+\xi}\right)^{K} 
= N\gamma\left(\frac{\xi}{1+\xi}\right)^{K}.
\]
For the power law representation, by the same arguments as \cref{ex:gpp-pl-rep-error},
\[
B_{N,K} \le  N\gamma\prod_{k=1}^{K}\frac{\alpha + kd}{\alpha + kd - d + 1}.
\]
This result generalizes that in \citet{Paisley:2012}, which applies only when $d = 0$. 
\paragraph{Size-biased representation} 
The size-biased representation of the beta process is well-established
and we refer the reader to \citet{Broderick:2012,Broderick:2014} for details.
We note that the standard beta integral yields
\[
\eta_{k} = \eta_{k1}
= \int \pi(\theta)^{k-1}(1-\pi(\theta))\nu(\dee\theta)
= \gamma \frac{\Gamma(\alpha + 1)}{\Gamma(\alpha + d)}\frac{\Gamma(\alpha + d + k - 1)}{\Gamma(\alpha + k)}, \label{eq:bprateint}
\]
and for $i > 1$, $\eta_{ki} = 0$. 
Hence $z_{ki} = 1$ almost surely and $\theta_{ki} \dist \distBeta(1-d, \alpha+d+k-1)$, demonstrating that the 
construction due to~\citet{Thibaux:2007} is a special case of the size-biased representation for $\distBP(\gamma, \alpha, d)$. 

To obtain a truncation bound for the Bernoulli likelihood case, first consider the $d > 0$ setting.
Using \cref{lem:gamma-ratio-sum} to simplify the sum in \cref{eq:sb-trunc}, we have
\[
B_{N,K}
&\le \frac{\gamma}{d}\frac{\Gamma(\alpha + 1)}{\Gamma(\alpha + d)}\left(\frac{\Gamma(\alpha+d+K)}{\Gamma(\alpha+K)}-\frac{\Gamma(\alpha+d+K+N)}{\Gamma(\alpha+K+N)}\right) \\
&\sim \gamma N \frac{\Gamma(\alpha + 1)}{\Gamma(\alpha + d)} K^{d-1} \qquad K \to \infty,
\]
where the asymptotic result follows from \cref{lem:taylorlimitapprox,lem:gamma-ratio-deriv,lem:asymp-properties}.
When $d=0$, we can again use \cref{lem:gamma-ratio-sum} to arrive at
\[
B_{N,K} 
&\le \gamma\alpha\left(\psi(\alpha+K)-\psi(\alpha+N+K)\right) 
\sim \gamma\alpha NK^{-1} \qquad K \to \infty,
\]
where $\psi(\cdot)$ is the digamma function, and the asymptotic result follows from \cref{lem:taylorlimitapprox}. 
We can also bound the truncation error in the case of the negative binomial likelihood. 
For a fixed number of failures $s > 0$, and assuming $\alpha + d+ (k-1)s > 1$, integration
by parts yields 
\[
\lefteqn{\int \pi(\theta)^{k-1}(1-\pi(\theta))\nu(\dee\theta)} \\
&= \frac{\gamma}{d}\frac{\Gamma(\alpha + 1)}{\Gamma(\alpha + d)}\left(\frac{\Gamma(\alpha+d+ks)}{\Gamma(\alpha+ks)} - \frac{\Gamma(\alpha+d+(k-1)s)}{\Gamma(\alpha+(k-1)s)} \right), \label{eq:bpnbrateint}
\]
When $d>0$, the sum from \cref{eq:sb-trunc} is telescoping, so canceling terms,
\[
B_{N,K}
&\leq  \frac{\gamma}{d}\frac{\Gamma(\alpha + 1)}{\Gamma(\alpha+d)}\left(\frac{\Gamma(\alpha+d+Ks)}{\Gamma(\alpha+Ks)} - \frac{\Gamma(\alpha+d+(K+N)s)}{\Gamma(\alpha+(K+N)s)}\right)
 \\
&\sim \gamma N s^d \frac{\Gamma(\alpha + 1)}{\Gamma(\alpha + d)} K^{d-1} \qquad K \to \infty,
\]
where the asymptotic result follows from \cref{lem:taylorlimitapprox,lem:gamma-ratio-deriv,lem:asymp-properties}.
To analyze the case where $d=0$, we can use L'Hospital's rule to take the limit of \cref{eq:bpnbrateint} as $d\to 0$,
yielding
\[
\lim_{d\to0} \int \pi(\theta)^{k-1}(1-\pi(\theta))\nu(\dee\theta)
&=\gamma\alpha(\psi(\alpha+ks)-\psi(\alpha+(k-1)s)). 
\]
Again computing the error bound by canceling terms in the telescoping sum,
\[
B_{N,K}
&\leq  \gamma\alpha\left(\psi(\alpha+Ks)-\psi(\alpha+(K+N)s)\right) 
\sim \gamma \alpha NK^{-1} \qquad K \to \infty,
\]
where the asymptotic result follows from an application of \cref{lem:taylorlimitapprox}.
\paragraph{Stochastic mapping}
We can transform the gamma process $\distGammaP(\gamma, \lambda, 0)$ into the beta process
$\distBP(\gamma, \alpha, 0)$ by applying the stochastic mapping
\[
\theta &\mapsto \theta/(\theta + G), & G &\dist \distGam(\alpha,\alpha), & \kappa(\theta, \dee u) &= \distGam(\theta/(\theta+u); \alpha, \alpha)\frac{\theta}{u^2}\dee u.
\]
Using \cref{lem:stochastic-mapping} yields $\kappa(\Theta)\dist\distBP(\gamma,\alpha,0)$.
Applying this result to the Bondesson representation for $\distGammaP(\gamma, \alpha, 0)$ yields
\[
\sum_{k=1}^{\infty}\theta_{k}\delta_{\psi_{k}} &\dist \distBP(\gamma, \alpha, 0), 
\quad \text{with} &
\theta_{k} &\defined (1 + G_{k}V_{k}^{-1}e^{\Gamma_{k,\alpha\gamma}})^{-1}, &
\psi_{k} &\distiid G, \\
&& G_{k} &\distiid \distGam(\alpha, \alpha), &
V_{k} &\distiid \distExp(\alpha). \nonumber
\]
which, unlike the Bondesson representation, applies for all $\alpha > 0$. 

\paragraph{Hyperpriors}
Consider truncating the Bondesson representation of the beta process, but with a hyperprior on the mass parameter $\gamma$.
A standard choice of hyperprior for $\gamma$ is a gamma distribution, i.e.~$\gamma \dist \distGam(a, b)$.
Combining \cref{prop:trunccrm-hyperprior} and the beta-Bernoulli truncation bound in \cref{eq:b-rep-bbp-trunc}, we have that
\[
B_{N,K} \le N\frac{a}{b}\left(\frac{\xi}{\xi+1}\right)^K. 
\]

\subsection{Beta prime process}\label{app:bppegs}

The beta prime process $\distBPP(\gamma,\alpha,d)$~\citep{Broderick:2014} 
with discount parameter $d \in [0,1)$, concentration parameter $\alpha > -d$, 
and mass parameter $\gamma > 0$, is a CRM with rate measure
\begin{align}
  \nu(\dee \theta) &= \gamma  \frac{\Gamma(\alpha+1)}{\Gamma(1-d)\Gamma(\alpha+d)} \theta^{-1-d}(1+\theta)^{-\alpha} \dee \theta.
\end{align}
The beta prime process is often paired with an odds Bernoulli likelihood,
\begin{align}
  h(x \given \theta) &= \ind\left[x\leq 1\right]\theta^x(1+\theta)^{-1},
\end{align}
in which case $\pi(\theta) = (1+\theta)^{-1}$. 
All truncation results are for the odds Bernoulli likelihood.

\begin{figure}[t!]
\begin{subfigure}[b]{.32\textwidth}
    \includegraphics[width=\columnwidth]{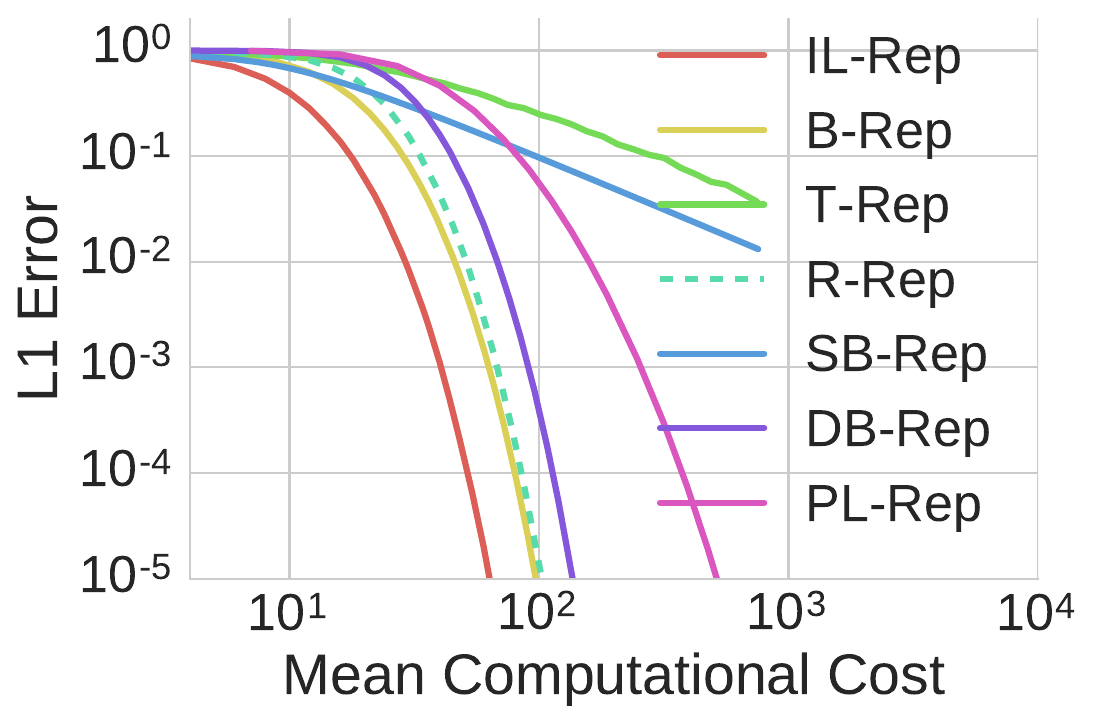}
    \caption*{$d=0.0$}
  \end{subfigure}
 \begin{subfigure}[b]{.32\textwidth}
    \includegraphics[width=.95\columnwidth, trim=33 0 0 0, clip]{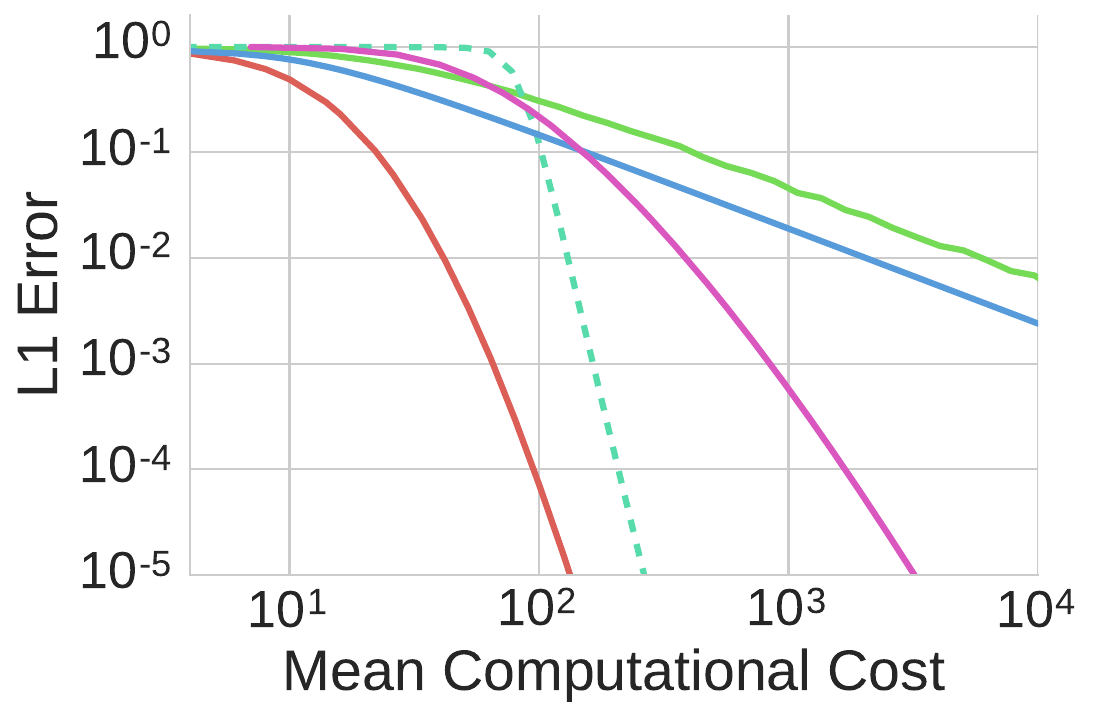}
    \caption*{$d=0.1$}
  \end{subfigure}
\begin{subfigure}[b]{.32\textwidth}
    \vspace{.7cm}\includegraphics[width=.95\columnwidth, trim=33 0 0 0, clip]{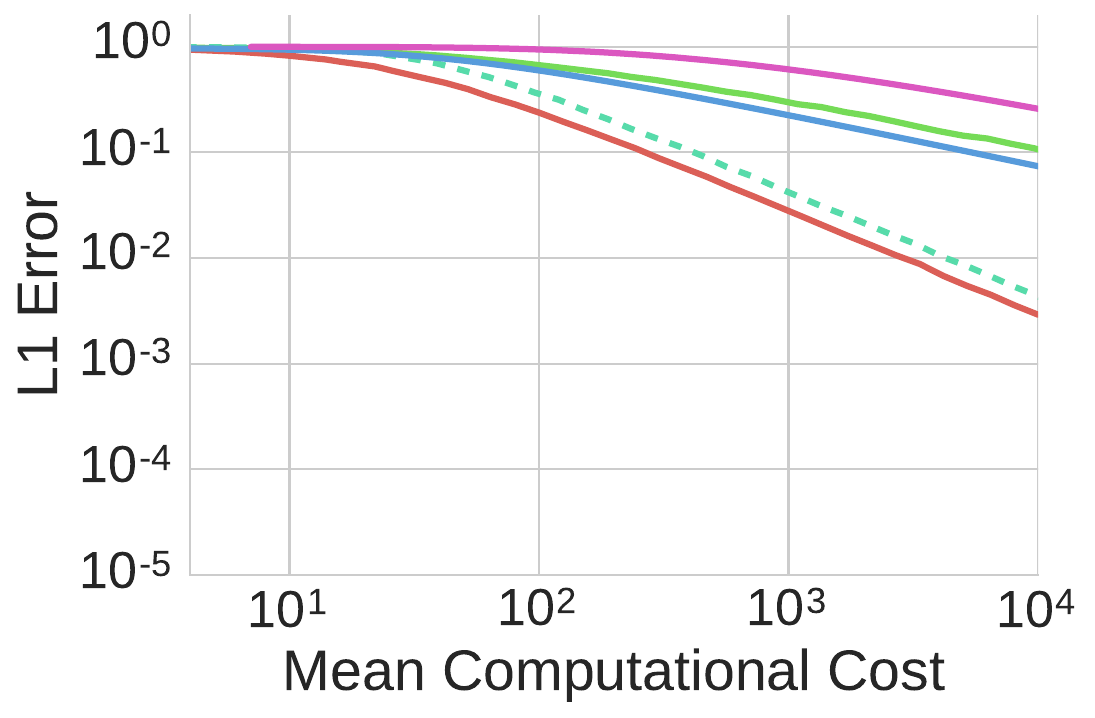}
    \caption*{$d=0.5$}
  \end{subfigure}
  \caption{Truncation error bounds for the beta prime-odds Bernoulli process. }\label{fig:bppfigs}
\end{figure}

\paragraph{Bondesson representation}
If $d = 0$, then $\theta \nu(\theta) = \gamma\alpha(1+\theta)^{-\alpha}$ is non-increasing
and $c_{\nu} = \lim_{\theta \to 0} \theta \nu(\theta) = \gamma\alpha$, so 
$g_\nu(v) = \alpha(1 + v)^{-\alpha-1} = \distBetaPrime(v; 1, \alpha)$.
Thus, it follows from \cref{prop:bondesson-representation} that  if
$\Theta\isBDrep(\gamma\alpha, \distBetaPrime(1, \alpha))$, then $\Theta \dist \distBPP(\gamma, \alpha, 0)$.
For the truncation bound we have
\[
B_{N,K}
&= N\int_{0}^{\infty} \left(1 - \EE\left[\pi(ve^{-G_{K}/(\gamma\alpha)})\right]\right)\nu(\dee v) \\
&= N\gamma\alpha\,\EE\left[e^{-G_{K}/(\gamma\alpha)}\int_{0}^{\infty}(1+v)^{-\alpha}(1+ve^{-G_{K}/(\gamma\alpha)})^{-1} \dee v\right]  \\
&\le N\gamma\alpha\,\EE[e^{-G_{K}/(\gamma\alpha)}]\int_{0}^{\infty}(1+v)^{-\alpha}(1+v\EE[e^{-G_{K}/(\gamma\alpha)}])^{-1} \dee v \\
&\le N\gamma\alpha \left(\frac{\gamma\alpha}{1 + \gamma\alpha}\right)^{K}\int_{0}^{\infty}(1+v)^{-\alpha-1} \dee v \\
&= N\gamma \left(\frac{\gamma\alpha}{1 + \gamma\alpha}\right)^{K},
\]
where the first upper bound follows from Jensen's inequality. 
Thus, the error bound is the same as for the beta-Bernoulli process. 

\paragraph{Thinning representation}
If we let $g = \distBetaPrime(1-d, \alpha + d)$, then the thinning representation for $\distBPP(\gamma, \alpha, d) $ is
\[
\Theta &= \sum_{k=1}^{\infty}V_{k}\ind(V_{k}(\Gamma_{k} - \gamma) \le \gamma)\delta_{\psi_{k}}, &
& \text{with} &
V_{k} \distiid \distBetaPrime(1 - d, \alpha + d). 
\]

\paragraph{Rejection representation}
For $d = 0$ and $\alpha \ge 1$, we take $\mu$ to be the rate measure for $\distLomP(\gamma, 1)$, so the
rejection representation is 
\[
\Theta = \sum_{k=1}^{\infty} V_{k} \ind(U_{k} \le (1 + V_{k})^{\alpha - 1}) \dist \distBPP(\gamma, \alpha, 0), \qquad \alpha \ge 1, \\
\text{with} \qquad V_{k} = (e^{\gamma^{-1}\alpha^{-1}\Gamma_{k}} - 1)^{-1},
\qquad 
U_{k} \distiid \distUnif[0,1].
\]
The expected number of rejections is $c_{\gamma} + \psi(\alpha)$, where $c_{\gamma}$ is the
Euler-Mascheroni constant and $\psi$ is the digamma function. 
Since $\psi(\alpha) \sim \log(\alpha)$ for $\alpha \to \infty$, the representation remains
efficient even for fairly large values of $\alpha$. 
To obtain a truncation bound, we use the same approach as in \cref{ex:rrepgamp}:
\[
\frac{B_{N,K}}{N\gamma\alpha}
&= \int_{0}^{\infty}F_{K}(\gamma\alpha\log(1+x^{-1}))(1+x)^{-\alpha-1}\dee x \\
&\le \int_{0}^{a}(1+x)^{-\alpha-1}\dee x  + F_{K}(\gamma\alpha\log(1+a^{-1}))\int_{a}^{\infty}(1+x)^{-\alpha-1}\dee x \\
&\le a + \alpha^{-1}(3\gamma\alpha\log(1+a^{-1})/K)^{K}(1+a)^{-\alpha} \\
&\sim e^{-b} + \alpha^{-1}(3\gamma\alpha/K)^{K}b^{K},
\]
where $b \defined \log(1+a^{-1})$. 
Setting the two terms equal and solving for $b$, we conclude that
\[
B_{N,K} 
&\le 2N\gamma\alpha e^{-K W_{0}(\{3\gamma\alpha\}^{-1})},
\]
where $W_0$ is the product log function, as defined in \cref{eq:lambertdefn}.

Similarly to \cref{ex:r-rep-gamma-process-power-law}, for the case of $d > 0$ and $\alpha \ge 0$, we instead use 
$\mu(\dee \theta) = \gamma\frac{\Gamma(\alpha+1)}{\Gamma(1-d)\Gamma(\alpha+d)} \theta^{-1-d}\dee \theta$. 
Since $\muinv(u) = (\gamma'u^{-1})^{1/d}$, where 
$\gamma' \defined \gamma\frac{\Gamma(\alpha+1)}{d\Gamma(1-d)\Gamma(\alpha+d)}$,
the rejection representation is
\[
\Theta &= \sum_{k=1}^{\infty} V_{k} \ind(U_{k} \le (1+V_{k})^{-\alpha})\delta_{\psi_{k}}, & 
&\text{with} &
V_{k} &=  (\gamma'\Gamma_{k}^{-1})^{1/d}, \\
&&&& U_{k} &\distiid \distUnif[0,1].
\]
The expected number of rejections is $\frac{\gamma\alpha}{d}$, so the representation is 
efficient for large $d$, but extremely inefficient when $d$ is small. 
We have
\[
B_{N,K} 
&= N\gamma \frac{\Gamma(\alpha+1)}{\Gamma(1-d)\Gamma(\alpha+d)} \int_{0}^{\infty} F_{K}(\gamma'x^{-d})x^{-d}(1+x)^{-1-\alpha}\dee x.
\]
Following the approach of \cref{ex:rrepgamp-power-law}, the integral can be upper bounded as
\[
\lefteqn{\int_{0}^{a}x^{-d}\dee x +  F_{K}(\gamma'a^{-d})\int_{a}^{\infty}x^{-d}(1+x)^{-1-\alpha}\dee x} \\
&\le (1-d)^{-1}a^{1-d} + \alpha^{-1}a^{-d} (3\gamma' K^{-1}a^{-d})^{K}.
\]
Setting the two terms equal and solving for $a$, we obtain
\[
B_{N,K} 
&\le 2 N \left( \frac{\gamma\Gamma(\alpha+1)}{\Gamma(1-d)\Gamma(\alpha+d)}\right)^{\frac{(d+1)K + 1}{dK + 1}}  (\alpha(1-d)^{dK})^{-\frac{1}{dK+1}}\left(\frac{dK}{3}\right)^{-\frac{K}{dK+1}} \\
&\sim 2 N \left( \frac{\gamma\Gamma(\alpha+1)}{\Gamma(2-d)\Gamma(\alpha+d)}\right)^{1+1/d}\left(\frac{dK}{3(1-d)^{d}}\right)^{-1/d} \qquad K \to \infty.
\]

\paragraph{Decoupled Bondesson representation}
It follows from \cref{prop:decoupled-bondesson-representation} and the same arguments as those in the 
Bondesson case that if $\Theta\isDBFrep(\gamma\alpha, \distBetaPrime(1, \alpha), \xi)$, then $\Theta \dist \distBPP(\gamma, \alpha, 0)$.
Using the trivial bound $\theta/(1+\theta) \le \theta$ and calculations
analogous to those in the beta-Bernoulli case, for $\alpha > 1$ we
obtain the upper bound 
\[
B_{N,K} \le 1 - \exp\left\{- N \frac{\gamma\alpha}{\alpha - 1}\left(\frac{\xi}{1+\xi}\right)^{K} \right\}
\sim N \frac{\gamma\alpha}{\alpha - 1}\left(\frac{\xi}{1+\xi}\right)^{K} \qquad K\to\infty.
\]

\paragraph{Power-law representation}
We can transform the gamma process $\distGammaP(\gamma, 1, d)$ into the beta prime process $\distBPP(\gamma, \alpha, d)$
by applying the stochastic mapping 
\[
\theta &\mapsto u=\theta/G, & G &\dist \distGam(\alpha+d, 1), & \kappa(\theta, \dee u) &= \distGam(\theta/u; \alpha+d, 1)\frac{\theta}{u^2}\dee u.
\]
Using \cref{lem:stochastic-mapping} yields $\kappa(\Theta) \dist \distBPP(\gamma, \alpha, d)$. Applying
this result to the power-law representation $\distGammaP(\gamma, 1, d)$ from \cref{ex:power-law-f-rep-gamma-process} 
yields the novel power-law representation
\[
\Theta\isPLFrep(\gamma\alpha, 1, d, \distBetaPrime(1, \alpha + d)) \quad\text{implies}\quad\Theta \dist \distBPP(\gamma, \alpha, d).
\]
Using the trivial bound $\theta/(1+\theta) \le \theta$ and calculations
analogous to those in the beta-Bernouli case, for $\alpha > 1$ we obtain the upper bound and asymptotic simplification
\[
B_{N,K} 
&\le N \frac{\gamma\alpha}{\alpha - 1} \prod_{k=1}^{K}\frac{1 + kd}{2 + kd - d} \\
&\sim  \frac{\gamma N\alpha}{\alpha - 1}\left\{\begin{array}{ll}
                          2^{-K} & d=0\\
                         \frac{\Gamma\left(\frac{2}{d}\right)}{\Gamma\left(\frac{1+d}{d}\right)}K^{1-d^{-1}}  & 0<d<1
                          \end{array}\right.  \qquad K\to\infty.                                
\]

\paragraph{Size-biased representation}
We have
\[
\eta_{k1} = \eta_{k} 
&= \int \pi(\theta)^{k-1}(1-\pi(\theta))\nu(\dee\theta)
=\gamma  \frac{\Gamma(\alpha+1)\Gamma(d+k+\alpha-1)}{\Gamma(\alpha+d)\Gamma(k+\alpha)},
\]
which is the same as for the beta-Bernoulli process.
Thus, the error bound is also the same as the beta-Bernoulli case.

\section{Technical lemmas} \label{app:technical-lemmas}

\bnlem \label{lem:convergence-to-zero}
If $Y_{k}$ is a uniformly bounded, non-negative sequence of random variables such that $\lim_{k \to \infty}\EE[Y_{k}] = 0$,
then $Y_{k} \convP 0$. 
\enlem
\bprf
Without loss of generality we assume that $Y_{k} \in [0,1]$ \as 
Then for all $\eps, \delta > 0$, by hypothesis there exists $k'$ such that for all $k \ge k'$, $\EE[Y_{k}] \le \eps\delta$.
It then follows from Markov's inequality that for all $k \ge k'$, $\Pr(Y_{k} > \eps) \le \delta$. 
\eprf

\bnlem \label{lem:small-measure-sets}
If $\mu$ is a non-atomic measure on $\reals^{d}$, then
for any $x \in \reals^{d}$ and $\delta > 0$, there exists $\eps_{x,\delta} > 0$ such that 
$\mu(\{ y \in \reals^{d} \given \|x - y\|_{2} \le \eps_{x,\delta}\}) \le \delta$.
\enlem
\bprf
Without loss of generality let $x = 0$. 
Suppose the implication does not hold. 
Then there exists $\delta > 0$ such that for all $\eps > 0$, $\mu(B_{\eps}) > \delta$,
where $B_{\eps} \defined \{ y \in \reals^{d} \given \|y\|_{2} \le \eps\}$. 
Let $\eps_{n}$ be a sequence such that $\eps_{n} \to 0$ as $n \to \infty$. 
Then by continuity $\mu(\{0\}) = \lim_{n\to\infty}\mu(B_{\eps_{n}}) > \delta$,
hence $\mu$ is a atomic, which is a contradiction. 
\eprf

\bnlem \label{lem:behavior-at-zero}
If $\nu(\dee \theta)$, an absolutely continuous $\sigma$-finite measure on $\reals_{+}$, 
and continuous $\phi : \reals_{+} \to [0,1]$ satisfy
\[
\nu(\mathbb{R}_+) &= \infty, &
\int \min(1, \theta) \nu(\dee\theta) &< \infty, &
    &\text{and}&
\int \phi(\theta)\nu(\dee\theta) &< \infty,
\]
then
\[
\lim_{\theta \to 0} \phi(\theta) = \phi(0) = 0.
\]
\enlem
\bprf
$\int \min(1, \theta) \nu(\dee\theta) < \infty$ implies that 
$\int_{0}^{1}\theta\nu(\dee \theta) < \infty$ and $\int_{1}^{\infty} \nu(\dee\theta) < \infty$.
But $\nu(\mathbb{R}_+) = \infty$, so $\int_{0}^{1}\nu(\dee \theta) = \infty$. 
In fact, for all $\eps > 0$, $\nu([0,\eps]) = \infty$ since otherwise 
$\int_{\eps}^{1}\nu(\dee\theta) = \infty \implies \int_{\eps}^{1}\theta\nu(\dee\theta) = \infty$,
a contradiction.
Since $\phi$ is continuous and has bounded range, $\lim_{\theta \to 0} \phi(\theta) = c$ exists and is finite. 
Assume $c > 0$, so $\exists \eps > 0$ such that $\forall \eps' < \eps$,
$|\phi(\eps') - c| < c/2$, and in particular $\phi(\eps') > c/2$. 
Thus, for any $\theta' < \eps$, $\int_{0}^{\theta'} \phi(\theta)\nu(\dee\theta) \ge c/2\int_{0}^{\theta'} \nu(\dee\theta) = \infty$,
a contradiction. 
Thus, $c = 0$. 
\eprf

\begin{nlem}\label{lem:taylorlimitapprox}
Assume $\phi(x)$ is a twice continuously differentiable function with the following properties:
\benum
\item $\phi''(x)/\phi'(x) \to 0$ as $x \to \infty$
\item for all $\delta > 0$ there exists $B_{\delta} > 0$ such that
for any increasing sequence $(x_{n})_{n=1}^{\infty}$,
\[
\limsup_{n \to \infty} \sup_{y \in [x_{n}, x_{n} + \delta]}\frac{\phi''(y)}{\phi''(x_{n})} = B_{\delta} < \infty.
\]
\eenum
Then for any constant $c > 0$ and any increasing sequence $(x_{n})_{n=1}^{\infty}$, 
\[
\phi(x_{n}+c) - \phi(x_{n}) &\sim c \phi'(x_{n}) & \text{for $n \to \infty$}. 
\]
\end{nlem}
\bprf
A second-order Taylor expansion of $\phi(x_{n}+c)$ about $x_{n}$ yields
\[
\phi(x_{n} + c) - \phi(x_{n}) 
&= c \phi'(x_{n}) + \frac{c^{2}}{2}\phi''(x_{n}^{*}),
\]
where $x_{n}^{*} \in [x_{n}, x_{n}+c]$. 
Our assumptions on $\phi$ ensure that,
\[
\lim_{n \to \infty} \frac{\phi''(x_{n}^{*})}{\phi'(x_{n})}
&\le  \lim_{n \to \infty} \frac{B_{c}\phi''(x_{n})}{\phi'(x_{n})} = 0
\]
and hence
\[
\lim_{n \to \infty} \frac{\phi(x_{n} + c) - \phi(x_{n})}{c\phi'(x_{n})} 
&= \lim_{n \to \infty} 1 + \frac{c}{2}\frac{\phi''(x_{n}^{*})}{\phi'(x_{n})} = 1.
\]
\eprf

\begin{nlem}[\citet{Gautschi:1959}]\label{lem:gautschi}
\[
 (1+x)^{d-1} \leq \frac{\Gamma(x+d)}{\Gamma(x+1)} \leq x^{d-1} \qquad 0 \leq d \leq 1, x \geq 1,
\]
and thus for $0\leq d \leq 1$,
\[
 \frac{\Gamma(x+d)}{\Gamma(x+1)} &\sim x^{d-1} &  x \to\infty.
\]
\end{nlem}

\begin{nlem}\label{lem:gamma-ratio-sum}
For $\alpha > 0$ and $x\geq -1$,
\[
\sum_{m=1}^M\frac{\Gamma(\alpha+m+x)}{\Gamma(\alpha+m)} =\left\{\begin{array}{ll} 
\frac{1}{1+x}\left(\frac{\Gamma(\alpha+M+x+1)}{\Gamma(\alpha+M)}-\frac{\Gamma(\alpha+x+1)}{\Gamma(\alpha)} \right) & x > -1\\
\psi(\alpha+M)-\psi(\alpha) & x=-1
\end{array}\right.
\]
where $\psi(\cdot)$ is the digamma function.
\end{nlem}
\bprf
When $M=1$ and $x > -1$, analyzing the right hand side yields
\[
\frac{\Gamma(\alpha+M+x+1)}{\Gamma(\alpha+M)}-\frac{\Gamma(\alpha+x+1)}{\Gamma(\alpha)} &= 
\frac{\Gamma(\alpha+x+2)}{\Gamma(\alpha+1)}-\frac{\Gamma(\alpha+x+1)}{\Gamma(\alpha)}\\
&= \frac{\Gamma(\alpha+x+1)}{\Gamma(\alpha)}\left(\frac{\alpha+x+1}{\alpha}-1\right)\\
&=(x+1)\frac{\Gamma(\alpha+x+1)}{\Gamma(\alpha+1)}.
\]
By induction, supposing that the result is true for $M-1 \geq 1$ and $x > -1$,
\[
\sum_{m=1}^M\frac{\Gamma(\alpha+m+x)}{\Gamma(\alpha+m)} &= \sum_{m=1}^{M-1}\frac{\Gamma(\alpha+m+x)}{\Gamma(\alpha+m)} + \frac{\Gamma(\alpha+M+x)}{\Gamma(\alpha+M)}\\
&= \frac{1}{1+x}\left(\frac{\Gamma(\alpha+M+x)}{\Gamma(\alpha+M-1)}-\frac{\Gamma(\alpha+x+1)}{\Gamma(\alpha)} \right) + \frac{\Gamma(\alpha+M+x)}{\Gamma(\alpha+M)}\\
&=\frac{\Gamma(\alpha+M+x)}{\Gamma(\alpha+M-1)}\frac{\alpha+M+x}{(1+x)(\alpha+M-1)} 
-\frac{\Gamma(\alpha+x+1)}{(1+x)\Gamma(\alpha)}\\ 
&=\frac{1}{1+x}\left(\frac{\Gamma(\alpha+M+x+1)}{\Gamma(\alpha+M)}-\frac{\Gamma(\alpha+x+1)}{\Gamma(\alpha)}\right).
\]
This demonstrates the desired result for $x > -1$. Next, when $x = -1$, we have that
\[
\sum_{m=1}^M\frac{\Gamma(\alpha+m-1)}{\Gamma(\alpha+m)} &= 
\sum_{m=1}^M\frac{1}{\alpha+m-1}.
\] 
We proceed by induction once again. 
For $M=1$, using the recurrence relation $\psi(x+1)=\psi(x)+x^{-1}$~\citep[Chapter 6]{Abramowitz:1964}, the right hand side evaluates to
\[
\psi(\alpha+1)-\psi(\alpha) = \psi(\alpha)+\alpha^{-1}-\psi(\alpha) = \alpha^{-1}. 
\]
Supposing that the result is true for $M-1 \geq 1$ and $x = -1$,
\[
\sum_{m=1}^M\frac{1}{\alpha+m-1} &= \sum_{m=1}^{M-1}\frac{1}{\alpha+m-1} + \frac{1}{\alpha+M-1}\\
&= \psi(\alpha+M-1)-\psi(\alpha) + \frac{1}{\alpha+M-1}\\
&= \psi(\alpha+M)-\psi(\alpha),
\]
demonstrating the result for $x=-1$.
\eprf

\begin{nlem}\label{lem:gamma-ratio-deriv}
For $a > 0$, $d \in \reals$, and $x_{n} \to \infty$,
\[
\der{}{x_{n}}\frac{\Gamma(a+x_{n}+d)}{\Gamma(a+x_{n})} \sim dx_{n}^{d-1}. 
\]
\end{nlem}
\bprf
We have
\[
\der{}{x}\frac{\Gamma(a+x+d)}{\Gamma(a+x)}
&= \frac{\Gamma(a+x+d)}{\Gamma(a+x)}(\psi(a+x+d) - \psi(a+x)),
\]
where $\psi$ is the digamma function. 
Using \cref{lem:taylorlimitapprox} and the asymptotic expansion of 
$\psi'$~\citep[Chapter 6]{Abramowitz:1964}, we obtain
\[
\psi(a+x_{n}+d) - \psi(a+x_{n}) \sim  d \psi'(a + x_{n}) \sim \frac{d}{x_{n} + a}\sim dx_n^{-1}.
\]
Since
\[
\frac{\Gamma(a+x_{n}+d)}{\Gamma(a+x_{n})} \sim (a + x_{n})^{d}\sim x_n^d,
\]
using \cref{lem:asymp-properties}(2) with the previous two displays yields the result. 
\eprf

\begin{nlem}\label{lem:em1nint}
For $0\leq d < 1$,
\[
\int_0^\infty \left(e^{-t\theta}-1\right)\theta^{-1-d}e^{-\lambda\theta}\dee\theta &= 
\left\{\begin{array}{ll}
\Gamma(-d)\left( (\lambda+t)^d-\lambda^d\right) & 0<d<1\\
\log \left(\frac{\lambda}{t+\lambda}\right) & d=0
\end{array}\right. .
\]
\end{nlem}
\bprf
By integration by parts and the standard gamma integral,
\[
\int_0^\infty \left(e^{-t\theta}-1\right)\theta^{-1-d}e^{-\lambda\theta}\dee\theta
&=\int_0^\infty \left[(\lambda+t)e^{-(\lambda+t)\theta}-\lambda e^{-\lambda \theta}\right] \frac{\theta^{-d}}{-d}\dee\theta\\
&=\Gamma(-d)\left( (\lambda+t)^d - \lambda^d\right).
\]
Taking the limit as $d\to 0$ via L'Hospital's rule yields
\[
\lim_{d\to 0}\Gamma(-d)\left( (\lambda+t)^d - \lambda^d\right) &= \log\left(\frac{\lambda}{\lambda+t}\right).
\]
\eprf

\begin{nlem}[Standard asymptotic equivalence properties]\label{lem:asymp-properties}
~\benum
\item If $a_{n} \sim b_{n}$ and $b_{n} \sim c_{n}$, then $a_{n} \sim c_{n}$.
\item If $a_{n} \sim b_{n}$ and $c_{n} \sim d_{n}$, then $a_{n}c_{n} \sim b_{n}d_{n}$. 
\item If $a_{n} \sim b_{n}$, $c_{n} \sim d_{n}$ and $a_{n}c_{n} > 0$, then $a_{n} + c_{n} \sim b_{n} + d_{n}$. 
\eenum
\end{nlem}

\section{Proofs of sequential representation results}
\label{app:seqproofs}

\subsection{Correctness of $\BDrep$, $\DBFrep$, and $\PLFrep$}

\begin{proofof}{\cref{prop:bondesson-representation}}
First, we show that $g_{\nu}(v)$ is a density. 
Since $v\nu(v)$ is nondecreasing, $\der{}{v}[v \nu(v)]$ exists almost everywhere, $\der{}{v}[v \nu(v)]\le 0$ and hence $g_{\nu}(v) \ge 0$.
Furthermore, 
\[
\int_{0}^{\infty} g_{\nu}(v) \dee v
&= -  c_{\nu}^{-1} \int_{0}^{\infty}\der{}{v}[v \nu(v)] \dee v 
= - c_{\nu}^{-1} v \nu(v) \Big|_{v=0}^{\infty} 
= 1,
\]
where the final equality follows from the assumed behavior of $v \nu(v)$ at
0 and $\infty$. 
Since for a partition $A_{1},\dots, A_{n}$, the random variables 
$\Theta(A_{1}), \dots, \Theta(A_{n})$ are independent, it suffices to 
show that for any measurable set $A$ (with complement $\bA$), the random variable $\Theta(A)$ has 
the correct characteristic function. 
Define the family of random measures 
\[
\Theta_{t} &= \sum_{k=1}^{\infty} V_{k} e^{-(\Gamma_{k}+t)/c_{\nu}}\delta_{\psi_{k}}, & t \ge 0,
\]
so $\Theta_{0} = \Theta$. 
Conditioning on $\Gamma_{11}$, 
\[
(\Theta_{t}(A) \given \Gamma_{1} = u) \eqD V_{1} e^{-(u+t)/c_{\nu}} \ind[\psi_{1} \in A] + \Theta_{t+u}(A),
\]
and note that the two terms on the left hand side are independent. 
We can thus write the characteristic function of $\Theta_{t}(A)$ as
\[
\varphi(\xi, t, A) &\defined \EE[e^{i\xi \Theta_{t}(A)}] \\
&= \EE[\EE[e^{i\xi \Theta_{t}(A)} \given \Gamma_{1} = u]] \\
&= \EE[\EE[e^{i\xi V_{1} e^{-(u+t)/c_{\nu}}\ind[\psi_{1} \in A]}e^{i\xi \Theta_{t+u}(A)} \given \Gamma_{1} = u]]\\
&= \EE[\EE[(G(A)e^{i\xi V_{1} e^{-(u+t)/c_{\nu}}} + G(\bA))\varphi(\xi, t+u, A) \given \Gamma_{1} = u]] \\
\begin{split}
&= G(A) \int_{0}^{\infty}\int_{0}^{\infty} e^{i\xi v e^{-(u+t)/c_{\nu}}}\varphi(\xi, t+u, A) g_{\nu}(v)e^{-u}\,\dee u\,\dee v \\
&\phantom{=~} + G(\bA) \int_{0}^{\infty} \varphi(\xi, t+u, A)e^{-u}\,\dee u,
\end{split}
\]
where $\bA$ is the complement of $A$. 
Multiplying both sides by $e^{-t}$ and making the change of variable $w = u + t$ yields
\[
\begin{split}
e^{-t}\varphi(\xi, t, A) 
&= G(A) \int_{t}^{\infty}\int_{0}^{\infty} e^{i\xi v e^{-w/c_{\nu}}}\varphi(\xi, w, A)g_{\nu}(v)e^{-w}\,\dee v\,\dee w \\
&\phantom{=~} + G(\bA) \int_{t}^{\infty} \varphi(\xi, w, A)e^{-w}\,\dee w 
\end{split} \\
\begin{split}
&= G(A) \int_{t}^{\infty} \varphi_{g_{\nu}}(\xi e^{-w/c_{\nu}})  \varphi(\xi, w, A) e^{-w}\dee w \\
&\phantom{=~} + G(\bA) \int_{t}^{\infty} \varphi(\xi, w, A)e^{-w}\dee w,
\end{split}
\]
where $\varphi_{g_{\nu}}(a) \defined \int_{0}^{\infty} e^{i a v}g_{\nu}(v)\,\dee v$ is the characteristic
function of a random variable with density $g_{\nu}$. 
Differentiating both sides with respect to $t$ and rearranging yields
\[
\D{\varphi(\xi, t, A)}{t}
&= \varphi(\xi, t, A) - G(A)\varphi_{g_{\nu}}(\xi e^{-t/c_{\nu}})\varphi(\xi, t, A) - (1-G(A))\varphi(\xi, t, A) \\
&= \varphi(\xi, t, A)G(A)(1 - \varphi_{g_{\nu}}(\xi e^{-t/c_{\nu}})),
\]
so we conclude that 
\[
\varphi(\xi, t, A) &= \exp\left(-G(A)\int_{t}^{\infty}(1 - \varphi_{g_{\nu}}(\xi e^{-u/c_{\nu}}))\,\dee u\right).
\]
Using integration by parts and the definition of $g_{\nu}$, rewrite 
\[
\varphi_{g_{\nu}}(a) 
&= -c_{\nu}^{-1}\int_{0}^{\infty} \der{}{v}[v\nu(v)]e^{ia v} \,\dee v \\
&= -c_{\nu}^{-1} v \nu(v) e^{i a v}\Big|_{v=0}^{\infty} + c_{\nu}^{-1}\int_{0}^{\infty} i a v \nu(v) e^{i a v}\,\dee v \\
&= 1 + \int_{0}^{\infty} \frac{i a v}{c_{\nu}} \nu(v) e^{i a v}\,\dee v,
\]
where the final equality follows from the assumed behavior of $v \nu(v)$ at $0$ and $\infty$. 
Combining the previous two displays and setting $t=0$ concludes the proof:
\[
\varphi(\xi, 0, A)
&= \exp\left(-G(A)\int_{0}^{\infty}\int_{0}^{\infty} \frac{i \xi v }{c_{\nu}} e^{-u/c_{\nu}} e^{i \xi v  e^{-u/c_{\nu}}} \nu(v) \,\dee v\,\dee u\right) \\
&= \exp\left(-G(A)\int_{0}^{\infty}\int_{0}^{\infty} \D{}{u}\left[-e^{i \xi v  e^{-u/c_{\nu}}}\right]\nu(v) \,\dee u\,\dee v\right) \\
&= \exp\left(-G(A)\int_{0}^{\infty}(e^{i\xi v} - 1) \nu(v)\,\dee v\right).
\]
\end{proofof}

\begin{proofof}{\cref{prop:decoupled-bondesson-representation}}
It was already shown that $g_{\nu}(v)$ is a density.
Let $\Theta'_{k} = \sum_{i=1}^{C_{k}}\theta_{ki}\delta_{\psi_{ki}}$, so 
$\Theta = \sum_{k=1}^{\infty} \Theta'_{k}$. 
Each $\Theta'_{k}$ is a CRM with rate measure $\frac{c_{\nu}}{\xi}\nu'_{k}(\dee \theta)$,
where $\nu'_{k}(\dee \theta)$ is the law of $\theta_{ki}$. 
Using the product distribution formula we have 
\vphantom{\cref{eq:rate-measure-for-theta-r}} %
\[
\nu'_{k}(\dee \theta) 
&= \int_{0}^{1} \frac{\xi^{k}}{\Gamma(k)}(-\log w)^{k-1}w^{\xi-2} g_{\nu}(\theta/w)\,\dee w\,\dee \theta. \label{eq:rate-measure-for-theta-r}
\]
Let $G_{\nu}(v) = \int_{0}^{v} g_{\nu}(x) \,\dee x$ be the cdf derived from $g_{\nu}$. 
From the preceding arguments, conclude that the rate measure of $\Theta$ is 
\[
\nu'(\dee \theta) 
&\defined \frac{c_{\nu}}{\xi} \sum_{k=1}^{\infty}\nu'_{k}(\dee \theta) \\
&= c_{\nu} \int_{0}^{1} \sum_{k=1}^{\infty} \frac{\xi^{k-1}}{\Gamma(k)}(-\log w)^{k-1}w^{\xi-2} g_{\nu}(\theta/w)\,\dee w\,\dee \theta \\
&= c_{\nu} \int_{0}^{1} \xi w^{-2} g_{\nu}(\theta/w)\,\dee w\,\dee \theta \\
&= c_{\nu} \int_{0}^{1} \xi \theta^{-1} \D{}{w}[-G_{\nu}(\theta/w)]\,\dee w\,\dee \theta\\
&= -c_{\nu} \theta^{-1} G_{\nu}(\theta/w) \Big|_{w=0}^{1}\,\dee \theta \\
&= c_{\nu} \theta^{-1}(1 - G_{\nu}(\theta))\,\dee \theta. 
\]
The cdf can be rewritten as
\[
1 - G_{\nu}(\theta) 
&= 1 + c_{\nu}^{-1}\int_{0}^{\theta} \der{}{x}[x \nu(x)] \,\dee x
= 1 + c_{\nu}^{-1} x \nu(x) \Big|_{x=0}^{\theta} 
=  c_{\nu}^{-1} \theta \nu(\theta).
\]
Combining the previous two displays, conclude that 
$\nu'(\dee \theta) = \nu(\theta)\,\dee \theta$.
\end{proofof}

\begin{proofof}{\cref{thm:power-law-rep}}
Since the power-law representation in \cref{eq:plrep} for the case when $V_{ki} = 1$ almost surely
was previously shown to be $\distBP(\gamma, \alpha, d)$~\citep{Broderick:2012}, we simply apply the stochastic mapping
result in \cref{lem:stochastic-mapping} with $\kappa(\theta, \dee u) = \theta^{-1}g_\nu(u \theta^{-1})\dee u$,
where $\theta^{-1}g_\nu(u \theta^{-1})$ is the density of $U = V\theta \given \theta$ under $V\dist g_\nu$.
\end{proofof}

\subsection{Proof of the expected number rejections of the $\Rrep$}

\begin{proofof}{\cref{prop:r-rep-efficiency}}
We have
\[
\EE\left[\sum_{k=1}^{\infty}\ind(\theta_{k} = 0)\right] 
&= \EE\left[\sum_{k=1}^{\infty}\ind\left(\frac{\dee\nu}{\dee\mu}(V_k) \geq U_k \right)\right]  \\
&= \EE\left[\sum_{k=1}^{\infty}\left(1 - \frac{\dee\nu}{\dee\mu}(V_k)\right)\right] \\
&= \int \left(1 - \frac{\dee\nu}{\dee\mu}(x)\right) \mu(\dee x),
\]
where the equalities follow from the definition of $\theta_{k}$, integrating out $U_{k}$, 
and applying Campbell's theorem.
\end{proofof}

\subsection{Power-law behavior of the $\PLFrep$}

We now formalize the sense in which power-law representations do in fact produce power-law behavior. 
Let $Z_{n} \given \Theta \distiid \distLP(\distPoiss, \Theta)$ and
$y_{k} \defined \sum_{n=1}^{N}\ind[z_{nk} \ge 1]$.
We analyze the number of non-zero features after $N$ observations,
\[
K_{N} \defined \sum_{k=1}^{\infty} \ind[y_{k} \ge 1],
\]
and the number of features appearing $j > 1$ times after $N$ observations,
\[
K_{N,j} \defined \sum_{k=1}^{\infty} \ind[y_{k} = j].
\]

In their power law analysis of the beta process, \citet{Broderick:2012} use a Bernoulli 
likelihood process.
However, the Bernoulli process is only applicable if $\theta_{k} \in [0,1]$, whereas in general
$\theta_{k} \in \reals_{+}$. 
Replacing the Bernoulli process with a Poisson likelihood process is a natural choice
since  $\ind[z_{nk} \ge 1] \dist \distBern(1 - e^{-\theta_{k}})$, and asymptotically
$1 - e^{-\theta_{k}} \sim \theta_{k}~\as$ for $k \to \infty$ since $\lim_{k \to \infty} \theta_{k} = 0~\as$ 
Thus, the Bernoulli and Poisson likelihood processes behave the same asymptotically, which is what is 
relevant to our asymptotic analysis. 
We are therefore able to show that all CRMs with power-law representations, not just the beta process, 
have what \citet{Broderick:2012} call Types I and II power law behavior.
Our only condition is that the tails of $g$ are not too heavy.

\bnthm \label{thm:f-rep-power-law-behavior}
Assume that $g$ is a continuous density such that for some $\eps > 0$,
\[
g(x) = O(x^{-1-d-\eps}). \label{eq:g-limiting-behavior}
\]
Then for $\Theta \isPLFrep(\gamma, \alpha, d, g)$ with $d > 0$, there exists a constant $C$
depending on $\gamma, \alpha, d$, and $g$ such that, almost surely, 
\[
K_{N} &\sim \Gamma(1-d) C N^{d},  & N &\to \infty \\
K_{N,j} & \sim \frac{d\,\Gamma(j-d)}{j!} C N^{d}, & N &\to \infty \qquad (j > 1). 
\]
\enthm

In order to prove \cref{thm:f-rep-power-law-behavior}, we require a number of additional 
definitions and lemmas. 
Our approach follows that in \citet{Broderick:2012}, which the reader is encouraged to 
consult for more details and further discussion of power law behavior of CRMs. 
Throughout this section, $\Theta \isPLFrep(\gamma, \alpha, d, g)$ with $d > 0$.
By \cref{lem:stochastic-mapping}, $\Theta \dist \distCRM(\nu)$, where
\[
\nu(\dee \theta) \defined \int g(\theta/u)u^{-1} \nu_{\distBP}(\dee u)\,\dee\theta
\]
and $\nu_{\distBP}(\dee \theta)$ is the rate measure for $\distBP(\gamma, \alpha, d)$. 
Let $\Pi_{k}$ be a homogeneous Poisson point process on $\reals_{+}$ with rate 
$\theta_{k}$ and define
\[
K(t) &\defined \sum_{k=1}^{\infty} \ind[|\Pi_{k} \cap [0,t]| > 0] \\
K_{j}(t) &\defined \sum_{k=1}^{\infty} \ind[|\Pi_{k} \cap [0,t]| = j].
\]
Furthermore, for $N \in \nats$, let 
\[
\Phi_{N} \defined \EE[K_{N}] 
\qquad \text{and} \qquad
\Phi_{N,j} \defined \EE[K_{N,j}] \quad (j > 1)
\]
and for $t > 0$, let 
\[
\Phi(t) \defined \EE[K(t)] 
\qquad \text{and} \qquad
\Phi_{j}(t) \defined \EE[K_{j}(t)] \quad (j > 1).
\]
If follows from Campbell's Theorem~\citep{Kingman:1993} that
\[
\Phi(t) &= \EE\left[\sum_{k}(1-e^{-t\theta_{k}})\right] = \int (1 - e^{-t\theta}) \nu(\dee \theta) \\
\Phi_{N} &= \EE\left[\sum_{k}(1-e^{-N\theta_{k}})\right] = \Phi(N) \\
\Phi_{j}(t) &= {N \choose j}\EE\left[\sum_{k}\frac{t^{j}(1-e^{-\theta_{k}})^{j}}{j!}e^{-t\theta_{k}}\right] = \frac{t^{j}}{j!}\int (1-e^{-\theta})^{j}e^{-t\theta} \nu(\dee \theta) \\
\Phi_{N,j} &= {N \choose j}\EE\left[\sum_{k}(1-e^{-\theta_{k}})^{j}e^{-(N-j)\theta_{k}}\right] = {N \choose j}\int (1-e^{-\theta})^{j}e^{-(N-j)\theta} \nu(\dee \theta).
\]

The first lemma characterizes the power law behavior of $\Phi(t)$ and $\Phi_{j}(t)$. 
A \emph{slowly varying function} $\ell$ satisfies $\ell(ax)/\ell(x) \to 1$ as $x \to \infty$ for all $a > 0$. 
\bnlem[{\citet{Broderick:2012}}, Proposition 6.1]
If for some $d \in (0,1)$, $C > 0$, and slowly varying function $\ell$,
\[
\bar\nu[0,x] \defined \int_{0}^{x} \theta \nu(\dee \theta) \sim \frac{d}{1-d} C \ell(1/x) x^{1-d}, \qquad x \to 0,
\label{eq:var-nu-asymptotics}
\]
then 
\[
\Phi(t) &\sim \Gamma(1-d) C t^{d}, && t \to \infty \\
\Phi_{j}(t) &\sim \frac{d\,\Gamma(j-d)}{j!} C t^{d}, && t \to \infty \qquad (j > 1).
\]
\enlem
Transferring the power law behavior from $\Phi(t)$ to $\Phi_{N}$ is trivial since $\Phi(N) = \Phi_{N}$. 
The next lemma justifies transferring the power law behavior from $\Phi_{j}(t)$ to $\Phi_{N,j}$. 
\bnlem[{\citet{Broderick:2012}}, Lemmas 6.2 and 6.3]
If $\nu$ satisfies \cref{eq:nuassump}, then 
\[
K(t) \uparrow \infty\, \as, \qquad \Phi(t) \uparrow \infty, \qquad \Phi(t)/t \downarrow 0. 
\]
Furthermore, 
\[
|\Phi_{N,j} - \Phi_{j}(N)| < \frac{C_{j}}{N} \max\{\Phi_{j}(N), \Phi_{j+2}(N)\} \to 0. 
\]
\enlem
The final lemma confirms that the asymptotic behaviors of $K_{N}$ and $K_{N,j}$ is 
almost surely the same as the expectations of $K_{N}$ and $K_{N,j}$. 
\bnlem
Assume $\nu$ satisfies \cref{eq:nuassump} and that for some $d \in (0,1)$, $C > 0$, $C_{j} > 0$,
and slowly varying functions $\ell$, $\ell'$, 
$\Phi(t) \sim C \ell(t) t^{d}$ and $\Phi_{j}(t) \sim C_{j}\ell(t)t^{d}$. 
Then for $N \to \infty$, almost surely
\[
K_{n} \sim \Phi_{N} 
\qquad \text{and} \qquad
\sum_{i < j} K_{N,i} \sim \sum_{i < j} \Phi_{N,i}. 
\]
\enlem

\begin{proofof}{\cref{thm:f-rep-power-law-behavior}}
Combining the three lemmas, the result follows as soon as we show that $\nu(\dee \theta)$ satisfies 
\cref{eq:var-nu-asymptotics}. 
$C$ will be a constant that may change from line to line.
We begin by rewriting $\nu(\dee \theta)$ using the change of variable $w = \theta(u^{-1} - 1)$:
\[
\nu(\dee \theta) 
&= C \int_{0}^{1} g(\theta/u)u^{-2-d}(1-u)^{\alpha+d-1}\,\dee u\,\dee\theta \\
&= C \theta^{-1-d} \int_{0}^{\infty} g(w + \theta) \frac{w^{\alpha+d-1}}{(w+\theta)^{\alpha-1}} \,\dee w\,\dee \theta. 
\]
Since $g(x)$ is integrable and continuous, for $x \in [0,1]$, it is upper-bounded by the non-integrable 
function $C_{0}x^{-1}$ for some $C_{0} > 0$.
Combining this upper bound with \cref{eq:g-limiting-behavior} yields
$g(x) \le \phi(x) \defined C_{0}x^{-1}\ind[x \le 1] + C_{1}x^{-1-d-\eps}\ind[x > 1]$ for some $C_{1} > 0$, so
\[
g(w + \theta) \frac{w^{\alpha+d-1}}{(w+\theta)^{\alpha-1}} \le \phi(w) w^{d}.
\]
Since $\phi(w) w^{d}$ is integrable, by dominated convergence the limit
\[
L = \lim_{\theta \to 0} \int_{0}^{\infty} g(w + \theta) \frac{w^{\alpha+d-1}}{(w+\theta)^{\alpha-1}}\,\dee w
\]
exists and is finite. 
Moreover, since $g(x)$ is a continuous density, there exists $M > 0$ and $0 < a < b < \infty$ such that 
$g(x) \ge M$ for all $x \in [a,b]$. 
Hence, for $\theta < a$, 
\[
\int_{0}^{\infty} g(w + \theta) \frac{w^{\alpha+d-1}}{(w+\theta)^{\alpha-1}}\,\dee w
&\ge M\int_{a-\theta}^{b-\theta} \frac{w^{\alpha+d-1}}{(w+\theta)^{\alpha-1}}\,\dee w
> 0, 
\]
so $L > 0$. 
Thus, 
\[
\psi(\theta) &\defined \theta\int_{0}^{1} g(\theta/u)u^{-2-d}(1-u)^{\alpha+d-1}\,\dee u\, \to C \theta^{-d},  & \theta \to 0
\]
and hence for $\delta > 0$ and $\theta$ sufficiently small, $|\psi(\theta) - C\theta^{-d}| < \delta$.
Thus, for $x$ sufficiently small, 
\[
\int_{0}^{x}\psi(\theta)\,\dee \theta 
&\le \int_{0}^{x} C\theta^{-d}\,\dee\theta + \int_{0}^{x}|\psi(\theta) - C\theta^{-d}|\,\dee\theta \\
&\le \frac{C\,x^{1-d}}{1-d} + \delta x \\
&\sim \frac{C\,x^{1-d}}{1-d}, \qquad x \to 0,
\]
which shows that \cref{eq:var-nu-asymptotics} holds. 

\end{proofof}

\section{Proofs of CRM truncation bounds}\label{app:truncproofs}

\subsection{Protobound}\label{app:protoboundproofs}

\begin{nlem}[Protobound]\label{lem:protobound}
  Let $\Theta$ and $\Theta'$ be two discrete random measures. 
  Let $X_{1:N}$   be a collection of random measures generated
  \iid from $\Theta$ with 
  $\supp(X_n) \subseteq \supp(\Theta)$, and
  let $Y_{1:N}$ be a collection of random variables
  where $Y_n$ is generated from $X_n$ via  $Y_n \given X_n \dist f(\cdot \given X_n)$.
  Define $Z_{1:N}$ and $W_{1:N}$ analogously for $\Theta'$.
  Finally, define $Q \defined \mathbbm{1}\left[ \supp(X_{1:N}) \subseteq
  \supp(\Theta')\right]$.
  If $(X_{1:N} | \Theta, \Theta', Q=1) \eqD (Z_{1:N} | \Theta', \Theta)$ almost surely under
  the joint distribution of $\Theta, \Theta'$, then
    \[
    \frac{1}{2}\|p_{Y} - p_{W}\|_{1}
    &\leq 1- \Pr(Q=1),
  \]
  where $p_Y, p_W$ are the marginal densities of $Y_{1:N}$ and $W_{1:N}$. 
\end{nlem}

\begin{proofof}{\cref{lem:protobound-crm}}
This is the direct application of \cref{lem:protobound} to CRMs, 
where 
$\Theta \dist\distCRM(\nu)$,
  and $\Theta'$ is a truncation $\Theta'=\Theta_K$. %
The technical condition is satisfied because the weights in $X_{1:N}$ are sampled independently for each atom in $\Theta$.
\end{proofof}

\begin{proofof}{\cref{lem:protobound-ncrm}}
This is the direct application of \cref{lem:protobound} to NCRMs, 
where $\Theta$ is the normalization of a CRM with distribution 
$\distCRM(\nu)$,
  and $\Theta'$ is the normalization of its truncation. %
The technical condition is satisfied because the conditioning on $X_{1:N} \subseteq \supp(\Theta')$ is equivalent to normalization of $\Theta'$.
\end{proofof}

\begin{proofof}{\cref{lem:protobound}}
We begin by expanding the 1-norm and conditioning on both $\Theta$ and
$\Theta'$ (denoted by conditioning on $\tTheta \defined (\Theta, \Theta')$ for brevity):
\[
\begin{aligned}
  \|p_Y-p_W\|_1 &= \int \left|
  \EE \left[\prod_{n=1}^N f(y_n|Z_n)\right] -
  \EE \left[\prod_{n=1}^N f(y_n|X_n)\right]\right|\dee y\\
  &= \int \left|
  \EE \left[\EE \left[\prod_{n=1}^N f(y_n|Z_n) |
  \tTheta\right] -
  \EE \left[\prod_{n=1}^N f(y_n|X_n) |
  \tTheta\right]\right]\right|\dee y.
\end{aligned}
\]
Then conditioning on $Q$,
\[
\lefteqn{\EE \left[\prod_{n=1}^N f(y_n|X_n) |
\tTheta\right] = \EE \left[\EE \left[\prod_{n=1}^N f(y_n|X_n) |
\tTheta, Q\right] | \tTheta\right]} \\
& \overset{\as}{=} \Pr(Q=1|\tTheta)\EE \left[\prod_{n=1}^N f(y_n|Z_n) |
\tTheta\right] +\Pr(Q=0|\tTheta) \EE \left[\prod_{n=1}^N f(y_n|X_n) |
\tTheta, Q=0\right], \nonumber
\]
where the first term arises from the fact that for any function $\phi$, 
\[
\EE \left[\phi(Z_{1:N}) | \tTheta\right] &\overset{\as}{=} \EE \left[ \phi(X_{1:N}) | \tTheta, Q=1\right],
\]
because $X_{1:N} | \tTheta, Q=1$ is equal in distribution to $Z_{1:N} |
\tTheta$ \as~by assumption.  
Substituting this back in above,
\[
  \lefteqn{\|p_Y-p_W\|_1}  \label{eq:crm-l1-upper-bounding} \\
  &= \int \left|
  \EE \left[\Pr(Q=0|\tTheta)\left(\EE \left[\prod_{n=1}^N f(y_n|Z_n) |
  \tTheta\right] -
  \EE \left[\prod_{n=1}^N f(y_n| X_n) |
  \tTheta, Q=0\right]\right)\right]\right|\dee y   \\
  &\leq \int 
  \EE \left[\Pr(Q=0|\tTheta)\left|\EE \left[\prod_{n=1}^N f(y_n|Z_n) |
  \tTheta\right] -
  \EE \left[\prod_{n=1}^N f(y_n|X_n) |
  \tTheta, Q=0\right]\right|\right]\dee y  \\
  &\leq\int 
  \EE \left[\Pr(Q=0|\tTheta)\left(\EE \left[\prod_{n=1}^N f(y_n|Z_n) |
  \tTheta\right] +
  \EE \left[\prod_{n=1}^N f(y_n|X_n) |
  \tTheta, Q=0\right]\right)\right]\dee y,  
\]
and finally by Fubini's Theorem,
\[
  \lefteqn{\|p_Y-p_W\|_1}    \\
  &\leq  
  \EE \left[\Pr(Q=0|\tTheta)\left(\EE \left[\int\prod_{n=1}^N f(y_n|Z_n)\dee y | \tTheta\right] +
  \EE \left[\int \prod_{n=1}^Nf(y_n|X_n)\dee y | \tTheta, Q=0\right]\right)\right]   \\
  &= \EE \left[\Pr(Q=0|\tTheta)\left(\EE \left[1 | \tTheta\right] + \EE
  \left[1| \tTheta, Q=0\right]\right)\right] \\
 & = 2\,\Pr(Q=0) = 2(1 - \Pr(Q=1)). 
\]
\end{proofof}

\subsection{Series representation truncation}

Recall from \cref{sec:seriestrunc} that a series representation generally has the form
\[
\Theta = \sum_{k=1}^\infty \theta_k \delta_{\psi_k} \qquad \theta_k = \tau(V_k, \Gamma_k) \qquad V_k \distiid g,
\]
where $\Gamma_k=\sum_{\ell=1}^k E_\ell$, $E_{\ell}\distiid\distExp(1)$, are the jumps in a unit-rate homogeneous Poisson point process, $\tau:\reals_+\times\reals_+\to\reals_+$ is a measurable function, $g$ is a distribution on $\reals_+$,
and $\lim_{u\to\infty}\tau(v, u) = 0$ for $g$-almost every $v$.
Note that by \cref{lem:behavior-at-zero} $\bpi(0) = 0$, where $\bpi(x) \defined 1 - \pi(x)$.
This fact will repeatedly prove useful for the proofs in this section. 

The proof of \cref{thm:series-rep-trunc} is based on the following lemma. 
\bnlem  \label{lem:series-rep-trunc}
Under the same hypotheses as \cref{thm:series-rep-trunc},
\[
\begin{split}
\lefteqn{\Pr(\supp(X_{1:N}) \subseteq \supp(\Theta_K))} \\
&= \EE\left[\exp\left\{-\int_{0}^{\infty}\left(1 - \int_{0}^{\infty}\pi(\tau(v, u + G_{K}))^{N}g(\dee v)\right)\,\dee u\right\}\right] \label{eq:series-rep-trunc-equality}.
\end{split}
\] 
\enlem 
\bprf
Let 
\[
p(t,K) \defined \EE\left[\prod_{k=K+1}^{\infty}\pi\left(\tau(V_{k}, \Gamma_{k}+t)\right)^{N}\right],
\]
so $p(0,K) = \Pr(\supp(X_{1:N}) \subseteq \supp(\Theta_K))$.
We use the proof strategy from \citet{Banjevic:2002} and induction in $K$.
For $K=0$,
\[
p(t,0)
&= \EE\left[\EE\left[\pi\left(\tau(V_{1}, u + t)\right)^{N}\prod_{k=2}^{\infty}\pi\Bigg(\tau\Bigg(V_{k}, \sum_{1 < j \le k}E_{j}+u+t\Bigg)\Bigg)^{N} \Bigg| \Gamma_{1} = u \right]\right] \\
&= \int_{0}^{\infty}\int\pi\left(\tau(v, u+t)\right)^{N} p(u+t,0)e^{-u}\,g(\dee v)\,\dee u
\]
since the $V_{k}$ are \iid
Multiplying both sides by $e^{-t}$ and making the change of variable $w = u + t$ yields
\[
e^{-t}p(t,0)
&= \int_{t}^{\infty}\int\pi\left(\tau(v, w)\right)^{N} p(w,0)e^{-w}\,g(\dee v)\,\dee w.
\]
Differentiating both sides with respect to $t$ and rearranging yields
\[
\D{p(t,0)}{t} 
&= p(t,0)\left(1 - \int\pi\left(\tau(v, t)\right)^{N}\,g(\dee v)\right). \label{eq:pt-diffeq}
\]
Since $\lim_{u\to\infty}\tau(v, u) = 0$ and  $\pi(0) = 1$ by \cref{lem:behavior-at-zero},
we can solve \cref{eq:pt-diffeq} and conclude that 
\[
p(t,0) 
&= \exp\left\{-\int_{t}^{\infty}\left(1 - \int\pi\left(\tau(v, u)\right)^{N}\,g(\dee v)\right)\,\dee u\right\} \\
&= \exp\left\{-\int_{0}^{\infty}\left(1 - \int\pi\left(\tau(v, u + t)\right)^{N}\,g(\dee v)\right)\,\dee u\right\}.
\]
We use the inductive hypothesis that
\[
p(t,K) &= \EE[p(t + G_{K},0)], & G_{K} \dist \distGam(K, 1), \quad G_{0} = 0,
\label{eq:il-ih}
\]
which trivially holds for $K = 0$. 
If the inductive hypothesis holds for some $K \ge 0$, then using the tower property,
\[
p(t,K+1)
&= \EE\left[\EE\left[\prod_{k=K+2}^{\infty}\pi\Bigg(\tau\Bigg((V_{k}, \sum_{1 < j \le k}E_{1j}+u+t\Bigg)\Bigg)^{N} \Bigg| \Gamma_{1} = u \right]\right] \\
&= \EE[p(t+E_{1},K)], \qquad E_{1} \dist \distExp(1) \\
&= \EE[p(t+G_{K}+E_{1},0)] \\
&= \EE[p(t+G_{K+1},0)].
\]
\cref{eq:series-rep-trunc-equality} follows by setting $t = 0$.
\eprf
\begin{proofof}{\cref{thm:series-rep-trunc}}
The main result follows by combining \cref{lem:protobound-crm,lem:series-rep-trunc},
applying Jensen's inequality, then using monotone convergence. 
The upper bound $1 - e^{-B_{N,K}} \le 1$ follows immediately from the fact that the integral
\[
\int_{0}^{\infty}\left(1 - \EE\left[\int_{0}^{\infty}\pi(\tau(v, u + G_{K}))^{N}g(\dee v)\right]\right)\,\dee u 
\]
is non-negative. 

Fix $N$. It follows from \cref{eq:hassump} that 
\[
\lim_{K \to \infty}1 - \Pr(\supp(X_{1:N}) \subseteq \supp(\Theta_K)) \to 0.
\] 
It then follows from \cref{lem:convergence-to-zero} that $1 - e^{-\omega_{N,K}} \convP 0$, where 
$\omega_{N,K} \defined \int_{0}^{\infty}\left(1 - \int_{0}^{\infty}\pi(\tau(v, u + G_{K}))^{N}g(\dee v)\right)\,\dee u$. 
By the continuous mapping theorem, conclude that $\omega_{N,K} \convP 0$ as $K \to \infty$ and hence 
$B_{N,K} = \EE[\omega_{N,K}] \to 0$ as $K \to \infty$. 
\end{proofof}

\bnthm[Inverse-L\'evy representation truncation error] \label{thm:il-rep-truncation-error}
For $\Theta \isILrep(\nu)$, the conclusions of \cref{thm:series-rep-trunc} hold with
\[
B_{N,K} = N\int_{0}^{\infty}F_{K}(\nu[x, \infty))(1 - \pi(x))\,\nu(\dee x).
\]
\enthm
\bprf
We have from \cref{thm:series-rep-trunc} that
\[
B_{N,K} = N \EE\left[\int_{0}^{\infty}(1 - \pi(\nuinv(u + G_{K}))\,\dee u\right].
\]
We first make the change of variables $x = \nuinv(u)$ to obtain
\[
\int_{0}^{\infty}\bpi(\nuinv(u + G_{K}))\,\dee u
= \int_{G_{K}}^{\infty} \bpi(\nuinv(u))\,\dee u 
= \int_{0}^{\nuinv(G_{K})} \bpi(x)\nu(\dee x)
\]
Finally, use the fact that for all $a, b \geq 0$, $\nuinv(a) \ge b \iff a \le \nu\left(\left[b, \infty\right)\right)$
and monotone convergence:
\[
B_{N,K} 
&= N \EE\left[\int_{0}^{\nuinv(G_{K})}\bpi(x)\,\nu(\dee x)\right] \\
&= N \EE\left[\int_{0}^{\infty}\ind[x \le \nuinv(G_{K})]\bpi(x)\,\nu(\dee x)\right] \\
&= N \EE\left[\int_{0}^{\infty}\ind[G_{K} \le \nu[x, \infty)]\bpi(x)\,\nu(\dee x)\right] \\
&= N \int_{0}^{\infty}F_{K}(\nu[x, \infty))\bpi(x)\,\nu(\dee x).
\]
\eprf

\bnthm[Thinning representation truncation error] 
For $\Theta \isTrep(\nu)$, the conclusions of \cref{thm:series-rep-trunc} hold with
\[
B_{N,K} = N\int_{0}^{\infty}(1 - \pi(v))\int_{0}^{\frac{\dee \nu}{\dee g}(v)}F_{K}\left(\frac{\dee \nu}{\dee g}(v) - u\right)\dee u\,g(v)\dee v.
\]
\enthm
\bprf
We have from \cref{thm:series-rep-trunc} that
\[
B_{N,K} = N\EE\left[\int_{0}^{\infty}\int_{0}^{\infty}\bpi\left(v\ind\left[\frac{\dee \nu}{\dee g}(v) \ge u + G_{K}\right]\right)\dee u\,g(v)\dee v\right]
\]
Since $\bpi(0) = 0$, using monotone convergence we have
\[
B_{N,K} 
&= N\int_{0}^{\infty}\bpi(v)\int_{0}^{\infty}\EE\left[\ind\left[\frac{\dee \nu}{\dee g}(v) \ge u + G_{K}\right]\right]\dee u\,g(v)\dee v \\
&= N\int_{0}^{\infty}\bpi(v)\int_{0}^{\infty}F_{K}\left(\frac{\dee \nu}{\dee g}(v) - u\right)\dee u\,g(v)\dee v \\
&= N\int_{0}^{\infty}\bpi(v)\int_{0}^{\frac{\dee \nu}{\dee g}(v)}F_{K}\left(\frac{\dee \nu}{\dee g}(v) - u\right)\dee u\,g(v)\dee v\\
&= N\int_{0}^{\infty}\bpi(v)\int_{0}^{\frac{\dee \nu}{\dee g}(v)}F_{K}\left(u\right)\dee u\,g(v)\dee v.
\]
\eprf

\bnthm[Rejection representation truncation error] 
For $\Theta \isRrep(\nu)$, the conclusions of \cref{thm:series-rep-trunc} hold with
\[
B_{N,K} = N\int_{0}^{\infty}F_{K}(\mu[x, \infty))(1 - \pi(x))\,\nu(\dee x).
\]
\enthm
\bprf
We have from \cref{thm:series-rep-trunc} that 
\[
B_{N,K} = N\EE\left[\int_{0}^{\infty}\int_{0}^{1}\bpi\left(\muinv(u + G_{K})\ind\left[\frac{\dee \nu}{\dee \mu}\left(\muinv(u + G_{K})\right)\ge v\right]\right)\dee v\,\dee u\right]
\]
Since $\bpi(0) = 0$, we can eliminate the innermost integral:
\[
\lefteqn{\int_{0}^{1}\bpi\left(\muinv(u + G_{K})\ind\left[\frac{\dee \nu}{\dee \mu}\left(\muinv(u + G_{K})\right)\ge v\right]\right)\dee v} \\
&= \frac{\dee \nu}{\dee \mu}(\muinv(u + G_{K}))\bpi(\muinv(u + G_{K})).
\]
Making the change of variable $x = \muinv(u + G_{K})$ and reasoning analogously to the 
proof of \cref{thm:il-rep-truncation-error}, we obtain
\[
B_{N,K}
&= N\EE\left[\int_{0}^{\muinv(G_{K})}\frac{\dee \nu}{\dee \mu}(x)\bpi(x)\mu(\dee x)\right] 
= N\int_{0}^{\infty}F_{K}(\mu[x,\infty))\bpi(x)\nu(\dee x).
\]
\eprf

\bnthm[Bondesson representation truncation error] \label{thm:bondesson-trunc}
For $\Theta \isRrep(\nu)$, the hypotheses of \cref{thm:series-rep-trunc} are satisfied and its conclusions hold with
\[
B_{N,K} = N \int_{0}^{\infty} \left(1 - \EE\left[\pi(ve^{-G_{K}})\right]\right)\nu(\dee v).
\]
\enthm
\bprf
While \cref{thm:bondesson-trunc} can be proved using \cref{thm:series-rep-trunc},
we take an alternative approach using more direct Poisson process arguments.

\begin{nlem}\label{lem:intermediate-ber-error-bound}
For $K \ge 0$, 
 $G_{K} \dist \distGam(K, c)$, and $G_{0} = 0$,
\[
\Pr(\supp(X_{1:N}) \subseteq \supp(\Theta_K)) 
= \EE\left[\exp\left\{ -\int_{0}^{\infty} \left(1 - \pi(ve^{-G_{K}})^{N}\right)\nu(\dee v)\right\}\right]. \label{eq:ber-event}
\]
\end{nlem}
\begin{proofof}{\cref{lem:intermediate-ber-error-bound}}
For $t\geq0$, the measure $t\Theta$ has distribution $t\Theta\dist\distCRM(\nu_t)$ where $\nu_t(\dee\theta)\defined \nu(\dee\theta/t)$. %
Further, define $\tX_n\given \Theta \distiid \distLP(h, t\Theta)$, and
\[
p(t,K) \defined \Pr(\supp(\tX_{1:N}) \subseteq \supp(t\Theta_K)) = \EE\left[\prod_{k=K+1}^{\infty}\pi\left(t V_{k}e^{-\Gamma_{k}/c}\right)^{N}\right].
\]
We will prove that
\[
p(t,K) = \EE\left[\exp\left\{ -\int_{0}^{\infty} \left(1 - \pi(t ve^{-G_{K}})^{N}\right)\nu(\dee v)\right\}\right],
\]
and then set $t=1$ to obtain the desired result.
The proof proceeds by induction. For $K=0$, the event $\supp(\tX_{1:N}) \subseteq \supp(t\Theta_K)$ 
is equal to $\supp(\tX_{1:N}) \subseteq \emptyset$ and thus $\supp(\tX_{1:N}) = \emptyset$. 
This is in turn equivalent to the probability that after thinning a $\distCRM(\nu_t)$
by $\pi(\theta)$ $N$ times, the remaining process has no atoms, i.e.~the probability 
that $\distCRM\left( \left(1-\pi(\theta)^N\right)\nu_t(\dee\theta)\right)$ has no atoms.
Since a Poisson process with measure $\mu(\dee\theta)$ has no atoms with probability $e^{-\int\mu(\dee\theta)}$,
\[
p(t, 0)  &= \exp\left(-\int_0^\infty\left(1-\pi(\theta)^N\right)\nu_t(\dee\theta)\right)\\
  &= \EE\left[\exp\left(-\int_0^\infty\left(1-\pi(t v e^{-G_0})^N\right)\nu(\dee v)\right)\right].
\]
The second equality follows by the change of variables $v = \theta/t$ and because $G_0 = 0$ with probability 1.
The inductive hypothesis is that for $K \geq 0$, 
\[
p(t, K) = \EE\left[p\left(t e^{-G_K}, 0\right)\right], \qquad G_K\dist\distGam(K, c).
\]
Using the tower property to condition on $\Gamma_{1}/c$ and the fact that the $V_k$ are \iid,
\[
p(t, K+1)
&=\EE\left[ \prod_{k=K+2}^\infty \pi\left(t V_k e^{-\Gamma_{k}/c}\right)^N\right]\\
&=\EE\left[\EE\left[ \prod_{k=K+1}^\infty \pi\left(t e^{-\Gamma_{1}/c}V_k e^{-\Gamma_{k}/c}\right)^N\given \Gamma_{1}/c\right]\right]\\
&= \EE\left[p\left(t e^{-\Gamma_{1}/c}, K\right)\right]
= \EE\left[p\left(t e^{-(\Gamma_{1}/c+G_K)}, 0\right)\right]
= \EE\left[p\left(t e^{-G_{K+1}}, 0\right)\right],
\]
since $\Gamma_{1}/c \dist \distExp(c)$. The desired result follows by setting $t=1$.
\end{proofof}

First combine \cref{lem:protobound-crm,lem:intermediate-ber-error-bound},
then apply Jensen's inequality. %
The bounds on $B_{N,K}$ and fact that $\lim_{K \to \infty} B_{N,K} = 0$ follows by 
the same arguments as in the proof of \cref{thm:series-rep-trunc}.
\eprf

\subsection{Superposition representation truncation}
\begin{proofof}{\cref{thm:superpositiontrunc}}
We begin with \cref{lem:protobound-crm}, and note that 
\[
\Pr\left(\supp(X_{1:N}) \subseteq \supp(\Theta_K)\right) = \Pr\left(\supp(X_{1:N}) \cap \supp(\Theta_K^+) = \emptyset\right).
\]
After generating $X_{1:N}$ from $\Theta$, we can view the point process representing the atoms in $\Theta_K^+$
not contained in any $X_n$ as $\Theta_K^+$ thinned by $\pi(\theta)^N$ (i.e.~the Bernoulli trial with success probability $1-\pi(\theta)$ 
to generate an atom failed $N$ times), and thus the remaining process is $\Theta_K^+$ thinned by $1-\pi(\theta)^N$. Therefore, the above event
is equivalent to the event that $\Theta_K^+$ thinned by $1-\pi(\theta)^N$ 
has no atoms. Using the fact that a Poisson process
with measure $\mu(\dee\theta)$ has no atoms with probability $e^{-\int\mu(\dee\theta)}$,
we have the formula for $B_{N,K}$,
\[
\Pr\left(\supp(X_{1:N}) \subseteq \supp(\Theta_K)\right) &= e^{-\int \left(1-\pi(\theta)^N\right)\nu_K^+(\dee\theta)}.
\]
Since $B_{N,K} \defined \int \left(1-\pi(\theta)^N\right)\nu_K^+(\dee\theta)$ is nonnegative, the error bound lies in the interval $[0, 1]$.
To show that $\lim_{K\to\infty} B_{N,K} = 0$, first note that
\[
 \hspace{-.3cm}\int\!\! \left(1-\pi(\theta)^N\right)\nu(\dee\theta) 
\leq N\!\!\int\!\! \left(1-\pi(\theta)\right)\nu(\dee\theta)
= N\sum_{x=1}^\infty\int\!\! h(x\given\theta)\nu(\dee\theta) < \infty, \label{eq:btotalfinite}
\]
by \cref{eq:hassump}. Further, splitting $\nu$ into its individual summed components, we have that
\[
\int (1-\pi(\theta))\nu(\dee\theta) 
&= \sum_{k=1}^K \int (1-\pi(\theta))\nu_k(\dee\theta) + \int(1-\pi(\theta))\nu_K^+(\dee\theta).\label{eq:splitnukkp}
\]
Combining the results from \cref{eq:btotalfinite,eq:splitnukkp} yields
\[
\lim_{K\to\infty} \int(1-\pi(\theta))\nu_K^+(\dee\theta) = 0.
\]
\end{proofof}

\subsection{Stochastic mapping truncation}

\begin{proofof}{\cref{prop:stochastic-mapping-trunc}}
Let $\tpi(u) = \tilde{h}(0 \given u)$. For notational brevity, define 
$Q$ to be the event where $\supp(X_{1:N}) \subseteq \supp(\Theta_K)$,
and $\tQ$ to be the corresponding transformed event where $\supp(\tX_{1:N}) \subseteq \supp(\tTheta_K)$.
Then  we have
\[
\Pr(\tQ)
&= \EE\left[{\textstyle\prod_{k=K+1}^{\infty}} \tpi(u_{k})^{\tN}\right].
\]
If $h(x \given \theta) = \distBern(x; 1-\pi_{\kappa,\tN}(\theta))$
and $N = 1$, then
\[
\Pr(\tQ)
&= \EE\left[{\textstyle\prod_{k=K+1}^{\infty} \int} \tpi(u)^{\tN}\kappa(\theta_{k}, \dee u)\right] 
= \Pr(Q).
\]
\end{proofof}

\subsection{Truncation with hyperpriors}\label{sec:crmhyperproof}

\begin{proofof}{\cref{prop:trunccrm-hyperprior}}
By repeating the proof of \cref{lem:protobound-crm} in \cref{app:protoboundproofs},
except with an additional use of the tower property to condition on the
hyperparameters $\Phi$, an additional use of Fubini's theorem to swap
integration and expectation, and Jensen's inequality, we have
\begin{align}
\frac{1}{2}\left\|p_Y-p_W\right\|_1 
&\leq \EE\left[1-\Pr\left(\supp(X_{1:N}) \subseteq \supp(\Theta_K) \given \Phi\right)\right]\\
&\leq \EE\left[1-e^{-B_{N,K}(\Phi)}\right]\\
&\leq 1-e^{-\EE\left[B_{N,K}(\Phi)\right]}.
\end{align} 
\end{proofof}

\section{Proofs of normalized truncation bounds} \label{app:normproofs}

\begin{proofof}{\cref{lem:gumbelmax}}
First, we demonstrate that the $\argmax$ is well-defined. Note that
\begin{align}
  \argmax_{i\in\nats} T_i+\log p_i = \argmax_{i\in\nats} \exp\left(T_i + \log
  p_i\right)
\end{align}
if it exists, due to the monotonicity of $\exp$. Similarly, existence of
either proves the existence of the other. Since $T_i$ are
i.i.d.~$\distGumbel(0, 1)$,
\begin{align}
  \Pr\left(\exp\left(T_i + \log p_i\right) > \epsilon\right)
 &= 1-\exp\left(-e^{-(\log\epsilon - \log p_i)}\right)\\
  &= 1-\exp\left(-p_i/\epsilon\right)
\leq 1-(1-p_i/\epsilon).
\end{align}
Therefore,
\begin{align}
  \sum_{i=1}^\infty \Pr\left(\exp\left(T_i + \log p_i\right) > \epsilon\right)
  &\leq \sum_{i=1}^\infty 1-(1-p_i/\epsilon)
  = \epsilon^{-1}\sum_{i=1}^\infty p_i < \infty.
\end{align}
This is sufficient to demonstrate that 
\begin{align}
  \exp\left(T_i + \log p_i\right) \convas 0 \qquad \text{as $i \to \infty$.}
\end{align}
Finally, since any positive sequence converging to $0$ can have only a finite number of 
elements greater than any $\epsilon > 0$, set $\epsilon = \exp(T_1+\log p_1)$, 
and thus 
\begin{align}
\argmax_{i\in\nats} \exp\left(T_i + \log p_i\right) &= \argmax_{i : T_i+\log p_i
\geq \epsilon} \exp\left(T_i + \log p_i\right)
\end{align}
where the right hand side exists because it computes the maximum of a finite,
nonempty set of numbers. Note that the $\argmax$ is guaranteed to be a single
element, since $T_i+\log p_i$ has a purely diffuse distribution on
$\mathbb{R}$.

Now that the \as~existence and uniqueness of the $\argmax$ has been demonstrated, we can compute
its distribution. First, note that
\begin{align}
  \Pr\left(T_i + \log p_i \leq x \, \forall i \in \nats, i\neq j\right) &=
   \hspace{-.3cm}\prod_{i=1, i\neq j}^\infty \exp\left(-e^{-(x - \log p_i)}\right)
  = \exp\left(-e^{-x}(s-p_j)\right),
\end{align}
where $s \defined \sum_i p_i$. So then
\begin{align}
  \Pr\left(j = \argmax_{i\in\nats} T_i + \log p_i\right) &= \Pr\left(T_i + \log
  p_i \leq T_j + \log p_j \, \forall i \in \nats\right)\\
  &= \int e^{-e^{-x}(s-p_j)} e^{-\left(x-\log p_j + e^{-(x- \log p_j)}\right)} \dee x\\
  &= p_j\int e^{-se^{-x}} e^{-x} \dee x\\
  &= \frac{p_j}{\sum_i p_i}\int e^{-e^{-(x-\log s)}} e^{-(x-\log s)} \dee x\\
  &= \frac{p_j}{\sum_i p_i},
\end{align}
where the last integral is 1 since its integrand is the $\distGumbel(s, 1)$ density. 
\end{proofof}

\subsection{Normalized series representation truncation}
\begin{proofof}{\cref{thm:ndbound}}
First, we apply \cref{lem:protobound-ncrm},
  \begin{align}
    \frac{1}{2}\left\|p_Y - p_W\right\|_1 &\leq 1- \Pr\left(X_{1:N} \subseteq \supp(\Xi_K)\right).
  \end{align}
  Next, by Jensen's
  inequality, 
  \begin{align}
    \Pr\left(X_{1:N} \subseteq \supp(\Xi_K)\right)= \EE\left[
    \left(\frac{\Theta_K\left(\Psi\right)}{\Theta\left(\Psi\right)}\right)^N\right]
    &\geq
    \EE\left[\frac{\Theta_K\left(\Psi\right)}{\Theta\left(\Psi\right)}\right]^N\\
    &= \Pr(X_1 \in
    \text{supp}(\Xi_K))^N.\label{eq:intermediate_step_innptb}
  \end{align}
The remaining
  part of this proof quantifies the probability that sampling $X_1$ from $\Xi$ 
  generates an atom in the support of $\Xi_K$ (equal to the support of $\Theta_K$, since $\Xi$ is just the normalization of $\Theta$). 
To do this, we use the trick based on \cref{lem:gumbelmax}:
we log-transform the rates in $\Theta$, perturb them all by \iid Gumbel random variables, and quantify
the probability that the max occurs within the atoms of $\Theta_K$.

First, we split the sequential representation of $\Theta$ into the truncation $\Theta_K$ and its tail $\Theta_K^+$, 
using the form from \cref{eq:general-tau-series-rep},
\[
\Theta = \sum_{k=1}^\infty \tau(V_k, \Gamma_k)\delta_{\psi_{k}}= \sum_{k=1}^K \tau(V_k, \Gamma_k)\delta_{\psi_{k}}+\sum_{k=K+1}^\infty \tau(V_k, \Gamma_k)\delta_{\psi_{k}}\defined \Theta_K+\Theta_K^+.
\]
Next, we define the maximum of the log-transformed, Gumbel perturbed rates in $\Theta_K$ as
\[
M_K &\defined \max_{1\leq k\leq K} \log\tau(V_k, \Gamma_k)+W_k,
\]
where $W_k \distiid \distGumbel(0, 1)$. 
Since $\Gamma_k$ are from a unit-rate homogeneous Poisson process, if we condition on the value of $\Gamma_K$, this is equivalent
to conditioning on the event that the Poisson process has exactly $K-1$ atoms on $[0, \Gamma_K]$, and the $K^\text{th}$ atom is at $\Gamma_K$.
And since $M_K$ does not depend on the ordering of $(\Gamma_k)_{k=1}^K$, 
\[
M_K \given \Gamma_K \eqD \max \left\{\log\tau(V_K, \Gamma_K)+W_K, \max_{1\leq k\leq K-1} \log\tau(V_k, U_k)+W_k \right\},
\]
where $U_k \distiid \distUnif(0, \Gamma_K)$. Therefore, $M_K \given \Gamma_K$ is the maximum of a collection of independent random variables,
so we can compute its CDF by simply taking the product of the CDFs of each of those random variables.
Using standard techniques for transformation of independent random variables,
we have that
\[
\Pr\left(\log\tau(V_k, U_k)+W_k \leq x \given\Gamma_K\right) &= \int_0^\infty \int_0^1 g(v) e^{-\tau\left(v, \Gamma_K u\right)e^{-x}}\dee u\dee v, \quad k< K\\
\Pr\left(\log\tau(V_K, \Gamma_K)+W_K \leq x\given \Gamma_K\right) &= \int_0^\infty g(v)e^{-\tau\left(v, \Gamma_K\right)e^{-x}}\dee v,
\]
so
\[
\Pr\left(M_K \leq x \given \Gamma_K\right) &= \left(\int_0^\infty \int_0^1 g(v) e^{-\tau\left(v, \Gamma_K u\right)e^{-x}}\dee u\dee v\right)^{K-1}\left(\int_0^\infty g(v)e^{-\tau\left(v, \Gamma_K\right)e^{-x}}\dee v\right).
\]
Defining the function
\[
J(u, t) \defined \EE\left[e^{-t \cdot \tau\left(V, u\right)}\right], \qquad V\dist g,
\]
we have
\[
\Pr\left(M_K \leq x \given \Gamma_K\right) &= \left(\int_0^1 J\left(\Gamma_K u, e^{-x}\right) \dee u\right)^{K-1}J\left(\Gamma_K, e^{-x}\right).
\]
Next, we define an analogous maximum for the tail process rates in $\Theta_K^+$,
\[
M_K^+ \defined \sup_{k>K} \log\tau(V_k, \Gamma_k)+W_k.
\]
Conditioning on $\Gamma_K$, we have
\[
M_K^+ \given \Gamma_K \eqD \sup_{k\geq 1} \log\tau(V_k, \Gamma'_k+\Gamma_K)+W_k,
\]
where $\Gamma'_k$ is a unit-rate homogeneous Poisson process on $\reals_+$. Now note that since $\Gamma'_k$ is a Poisson point process on $\reals_+$,
so is $\log \tau(V_k, \Gamma'_k+\Gamma_K)+W_k$ (using Poisson process stochastic mapping), with rate measure 
\[
\left(\int_0^\infty\int_0^\infty e^{-(t-\log\tau(v, u+\Gamma_K)) - e^{-(t-\log\tau(v, u+\Gamma_K))}} g(v)\dee u\dee v\right) \dee t.
\]
Therefore, $\Pr\left(M_K^+ \leq x \given \Gamma_K\right)$ is equal to the probability that the above Poisson point process
has no atoms with position greater than $x$. Since a Poisson process on $\reals_+$ with measure $\mu$ has no atoms above $x$ 
with probability $e^{-\int_x^\infty \mu(\dee t)}$,
\[
\Pr\left(M_K^+ \leq x \given \Gamma_K\right) &= e^{-\int_{x}^\infty \left(\int_0^\infty\int_0^\infty e^{-(t-\log\tau(v, u+\Gamma_K)) - e^{-(t-\log\tau(v, u+\Gamma_K))}} g(v)\dee u\dee v\right) \dee t}.
\]
Noticing that the integrand in $t$ is a Gumbel density, we can use Fubini's theorem to swap integrals and evaluate:
\[
\Pr\left(M_K^+ \leq x \given \Gamma_K\right) &= e^{\int_0^\infty\int_0^\infty\left(e^{-\tau(v, u+\Gamma_K)e^{-x}}-1\right) g(v)\dee u\dee v}\\
&=e^{\int_0^\infty\left( J\left(u+\Gamma_K, e^{-x}\right) - 1 \right)\dee u}.
\]
Taking the derivative yields the density of $M_K^+ \given \Gamma_K$ with respect to the Lebesgue measure.
Therefore, using the Gumbel-max trick from \cref{lem:gumbelmax}, we substitute the results for $M_K \given \Gamma_K$ and $M_K^+ \given \Gamma_K$ into the original bound yielding
\[
    \frac{1}{2}\left\|p_Y - p_W\right\|_1 &\leq 1- \Pr\left(X_{1} \subseteq \supp(\Xi_K)\right)^N\\
&= 1-\left(1-\EE\left[\Pr\left(M_K < M_K^+ \given \Gamma_K\right)\right]\right)^N,
\]
where (using the substitution $t = e^{-x}$)
\[
&\Pr\left(M_K < M_K^+ \given \Gamma_K\right) = \int_{-\infty}^\infty \Pr\left(M_K \leq x \given \Gamma_K\right)\frac{\dee}{\dee x}\Pr\left(M_K^+\leq x\given\Gamma_K\right)\dee x\\
= &\int_{0}^{\infty} J\left(\Gamma_K, t\right)\left(\int_0^1 J\left(\Gamma_K u, t\right) \dee u\right)^{K-1} 
\left(-\frac{\dee}{\dee t}e^{\int_0^\infty\left( J\left(u+\Gamma_K, t\right) - 1 \right)\dee u}\right) \dee t.
\]
The fact that the bound is between 0 and 1 is a simple consequence of the fact that $\Pr\left(X_1 \in \supp(\Xi_K)\right) \in [0, 1]$, and the asymptotic result follows from the fact that $\Pr\left(X_1\in\supp\left(\Xi_K\right)\right) \to 1$ as $K\to\infty$.
\end{proofof}

\subsection{Normalized superposition representation truncation}\label{sec:nrmsupertruncprf}
\begin{proofof}{\cref{thm:nfvbound}}
The same initial technique as in the proof of \cref{thm:nfvbound}
yields
\begin{align}
\frac{1}{2} \left\|p_Y-p_W\right\|_1 \leq 1- \Pr\left(X_1 \in \supp\left(\Xi_K\right)\right)^N.
\end{align}
  The remaining
  part of this proof quantifies the probability that sampling $X_1$ from $\Xi$ 
  generates an atom in the support of $\Xi_K$. Since most of the following developments are similar 
  for $\Theta_K, \nu_K$ and
  $\Theta_{K}^+, \nu_K^+$, we will focus the discussion on $\Theta_K, \nu_K$ and reintroduce
  the tail quantities when necessary. 
  First, we transform the rates of $\Theta_K$ under the stochastic mapping $w = \log \theta + W$, where $W\distiid \distGumbel(0, 1)$,
resulting in a new Poisson point process with rate measure
\begin{align}
  \left(\int e^{-(t-w) -e^{-(t-w)}} e^w\nu_K(e^w)\dee w\right) \dee t.
\end{align}
The probability that all points in
this Poisson point process are less than a value $x$ is equal 
to the probability that there are no atoms above $x$. Defining $M_K$ to be
the supremum of the points in this process, combined with the basic
properties of Poisson point processes, we have
\begin{align}
  \Pr(M_K \leq x) &= \exp\left(-\int_{x}^\infty\int e^{-(t-w) -e^{-(t-w)}}
  e^w\nu_K(e^w)\dee w \dee t\right)
\end{align}
Using Fubini's theorem to swap the integrals, we can evaluate the inner integral
analytically by noting the integrand is a $\distGumbel(w, 1)$ density,
\begin{align}
  \Pr(M_K \leq x) &= \exp\left(-\int (1-e^{-e^{-(x-w)}})
  e^w\nu_K(e^w)\dee w \right)\\
  &= \exp\left(\int \left(e^{-\theta e^{-x}}-1\right) \nu_K(\dee\theta)
  \right).
\end{align}
We can take the derivative with respect to $x$ to obtain its
density with respect to the Lebesgue measure.
The above derivation holds true for $\Theta_K^+$, replacing $\nu_K$ with $\nu_K^+$ and $M_K$ with $M_K^+$.
Therefore, using the Gumbel-max trick from \cref{lem:gumbelmax},
we substitute the results for $M_K$ and $M_K^+$ into the original bound, yielding 
\[
\frac{1}{2} \left\|p_Y-p_W\right\|_1 &\leq 1- \Pr\left(X_1 \in \supp\left(\Theta_K\right)\right)^N\\
& = 1-\left(1- \Pr\left(M_K < M_K^+\right) \right)^N,
\]
where (using the substitution $t=e^{-x}$, and the fact that $\theta e^{-\theta t}$ is dominated by $\theta$ to swap integration and differentiation)
\[
\Pr\left(M_K < M_K^+\right) &= \int_{-\infty}^\infty\Pr(M_K \leq x) \frac{\dee}{\dee x} \Pr\left(M_K^+ \leq x\right) \dee x\\
&= \int_{-\infty}^\infty e^{\int \left(e^{-\theta e^{-x}}-1\right) \nu(\dee\theta)} \frac{\dee}{\dee x}\left(\int \left(e^{-\theta e^{-x}}-1\right) \nu_K^+(\dee\theta)\right) \dee x\\
&= -\int_{0}^\infty e^{\int \left(e^{-\theta t}-1\right) \nu(\dee\theta)} \frac{\dee}{\dee t}\left(\int \left(e^{-\theta t}-1\right) \nu_K^+(\dee\theta)\right) \dee t\\
&= \int_{0}^\infty e^{\int \left(e^{-\theta t}-1\right) \nu(\dee\theta)} \left(\int\theta e^{-\theta t}\nu_K^+(\dee\theta)\right) \dee t.
\]
  The fact that the bound lies between 0 and 1 is a simple consequence of the fact that $\Pr\left(X_1 \in \supp(\Xi_K)\right) \in [0,1]$. 
  The fact that the error bound asymptotically approaches 0  is a consequence of the
  monotone convergence theorem applied to the decreasing sequence of functions 
$\theta\nu_K^+(\dee\theta)$.
\end{proofof}

\subsection{Truncation with hyperpriors}
\begin{proofof}{\cref{prop:truncncrm-hyperprior}}
By repeating the proof of \cref{lem:protobound-ncrm} in \cref{app:protoboundproofs},
except with an additional use of the tower property to condition on the
hyperparameters $\Phi$, an additional use of Fubini's theorem to swap
integration and expectation, and Jensen's inequality, we have
\begin{align}
\frac{1}{2}\left\|p_Y-p_W\right\|_1 &\leq \EE\left[1-\Pr\left(X_{1:N} \subseteq \supp(\Theta_K) \given \Phi\right)\right]\\
&\leq 1-\EE\left[(1-B_K(\Phi))^N\right]\\
&\leq 1-\left(1-\EE\left[B_K(\Phi)\right]\right)^N.
\end{align} 
\end{proofof}

\bibliographystyle{ba}
\bibliography{mybib,tbib}

\end{document}